\documentclass[a4paper,10pt]{amsart}

\usepackage{amsmath,amsfonts,amsthm, amssymb, epsfig}

\def\C{\mathbb{C}}
\def\Z{\mathbb{Z}}
\def\N{\mathbb{N}}

\def\qed{$\hfill \diamondsuit$}
\def\U{\mathbf{U}}
\def\H{\mathbf{H}}
\def\p{\mathbf{p}}
\def\h{\mathfrak{h}}
\def\g{\mathfrak{g}}
\def\n{\mathfrak{n}}
\def\X{\mathbb{X}}
\def\x{\vec{x}}
\def\c{\vec{c}}
\def\l{\underline{\lambda}}
\def\O{{\mathcal{O}}}
\def\Q{\mathbb{Q}}

\newtheorem{theo}{Theorem}
\newtheorem{prop}{Proposition}
\newtheorem{lem}{Lemma}
\newtheorem{cor}{Corollary}
\newtheorem{claim}{Claim}
\newtheorem*{conj}{Conjecture}
\numberwithin{equation}{section}
\numberwithin{lem}{section}
\numberwithin{prop}{section}
\numberwithin{cor}{section}
\numberwithin{theo}{section}
\numberwithin{claim}{section}

\begin{document}

\title{ Noncommutative Projective Curves and Quantum Loop Algebras}
\author{Olivier Schiffmann}
\date{}
\begin{abstract} 
We show that the Hall algebra of the category of coherent sheaves 
on a weighted projective line over a finite field provides a realization of 
the (quantized) enveloping algebra of a certain nilpotent subalgebra
of the affinization of the 
corresponding
Kac-Moody algebra. In particular, this yields a geometric realization of 
the quantized
enveloping algebra of elliptic (or 2-toroidal) algebras of types 
$D_4^{(1,1)}$, $E^{(1,1)}_6$, $E^{(1,1)}_7$ and $E_{8}^{(1,1)}$
in terms of coherent sheaves on weighted projective lines of genus one,
or equivalently in terms of equivariant coherent sheaves on elliptic curves.
\end{abstract} 

\maketitle
\noindent

{\centerline{\small{\textit{Dedicated to Igor Frenkel 
on the occasion of his 50th birthday}}}

\vspace{.1in}

\paragraph{} \textbf{Introduction.} The geometric approach to quantum groups
developed in the past 15 years is based on a deep relationship between simple
or affine (and to a lesser extent Kac-Moody) Lie algebras and certain 
finite-dimensional hereditary algebras. More precisely, let $\mathfrak{g}$
be a Kac-Moody algebra and $\Gamma$ its Dynkin diagram. 
Ringel proved in \cite{Ri1}
that the Hall algebra of the category of representations, over a finite field
$\mathbb{F}_q$, of a quiver whose underlying graph is $\Gamma$ provides
a realization of the ``positive part'' $\U^+_q(\mathfrak{g})$ of the 
Drinfeld-Jimbo quantum group associated to $\mathfrak{g}$. This result was the
starting point of Lusztig's geometric construction of the canonical basis
of $\U_q^+(\mathfrak{g})$ (see \cite{Lu}).

\vspace{.1in}

Another natural example of a hereditary category is provided by the category
$Coh(X)$ of coherent sheaves on a smooth projective curve $X$. In the 
remarkable paper \cite{K},
Kapranov observed a striking analogy between a function field analog
of the algebra of unramified automorphic forms (for $GL(N)$ for all $N$)
and Drinfeld's loop-like
realization of quantum affine algebras. In particular,
his result provides an isomorphism between a natural subalgebra
of the Hall algebra of the category $Coh(\mathbb{P}^1(\mathbb{F}_q))$
and Drinfeld's ``positive part'' $\U^+_q(\widehat{\mathfrak{sl}}_2)$.

\vspace{.1in}

In this paper we extend Kapranov's result. Rather than considering higher genus
smooth projective curves (whose Hall algebras are wild) we study the Hall 
algebra of the category of coherent sheaves on certain ``noncommutative smooth
projective curves''- the so-called \textit{weighted projective lines}, 
introduced by
Geigle and Lenzing (\cite{GL}). These fundamental examples of ``noncommutative
smooth projective curves'' have attracted some attention lately (see e.g
\cite{Happ2} and \cite{CB}). Our main result (Theorem~5.2) provides, when 
$\mathfrak{g}$ is a 
simple Lie algebra  or an affine Lie algebra of type $D_4^{(1)}, E_{6}^{(1)},
E_7^{(1)}$ or $E_8^{(1)}$, a natural embedding
$\U_q(\widehat{\n}) \hookrightarrow \H_{Coh(\X_{\mathfrak{g}})}$
of a certain ``positive part'' of $\U_q(\mathfrak{L}\mathfrak{g})$ into the Hall
algebra of a suitable weighted projective line $\X_{\mathfrak{g}}$.
Here $\mathfrak{L}\g$ is the universal central extension of $\g[t,t^{-1}]$.

\vspace{.1in}

In particular, this gives a geometric construction of the quantum toroidal
algebras of type $D_4, E_6, E_7$ and $E_8$. In this case, the two loops have a 
natural
interpretation as the two discrete invariants (rank and degree) of a coherent 
sheaf on $\X_{\mathfrak{g}}$. Note that the exceptional role played
by these four types among all affine Dynkin diagrams is well-known in the 
theory of quadratic forms
(see e.g \cite{Ri3}). 
Our computations extend those of Bauman and Kassel in \cite{BK}, where the
case of the projective line $\mathbb{P}^1$ is considered in details.
Our proof also makes use in an essential way of the 
Harder-Narasimhan filtration on semistable coherent sheaves and on the theory
of mutations. The use of the Harder-Narasimhan filtration in the context
of Hall algebras is also noted in \cite{Re}.

\vspace{.1in}

Recently, Y. Lin and L. Peng independently obtained another construction
of the enveloping algebra of double-loop algebras of type $D_4, E_6,E_7$
and $E_8$ by considering the root category (as in \cite{PX}) of certain
tame algebras $\Lambda_\mathfrak{g}$ (the \textit{tubular algebras}).
Our construction is related to theirs by the existence, in this case, of a 
derived equivalence
$$D^b(Coh(\X_{\mathfrak{g}})) \simeq D^b(Rep(\Lambda_{\mathfrak{g}})).$$

\vspace{.1in}

Finally, we note the following link of our work with the McKay correspondence.
Let $\mathfrak{g}$ be a simple Lie algebra of type $A_{2n}, D_n$ or
$E_k$ for $k=6,7,8$, and let $\Gamma \supset \{\pm 1\}$
be the finite subgroup of $SL(2,\C)$ associated to $\mathfrak{g}$ by the McKay
correspondence, and set $\Gamma'=\Gamma/\{\pm 1\} \in SO(3)$. 
When the base field is 
algebraically closed, the (ramified) quotient $\mathbb{P}^1/\Gamma'$ and 
the (smooth, noncommutative) curve $\X_{\mathfrak{g}}$ are related
by an equivalence of categories $Coh(\X_{\mathfrak{g}}) \simeq Coh_{\Gamma'}
(\mathbb{P}^1)$ and a
 derived equivalence
$$D^b(Coh(\X_{\mathfrak{g}})) \simeq 
D^b(\Lambda_{\mathfrak{g}}) 
\simeq D^b_{\Gamma'}(Coh(\mathbb{P}^1)),$$
where $\Lambda_{\mathfrak{g}}$ is the path algebra of the McKay quiver
of $\Gamma$, and $D^b_{\Gamma'}(Coh(\mathbb{P}^1))$ is the 
$\Gamma'$-equivariant
derived category of $Coh(\mathbb{P}^1)$.
Furthermore, the minimal 
desingularization $\widetilde{T^*(\mathbb{P}^1)/\Gamma'}$ is isomorphic
to the minimal desingularization $\widetilde{\C^2//\Gamma}$ of the Kleinian
singularity $\C^2//\Gamma$ (see \cite{Lam}). Such spaces (the so-called
ALE spaces), and more generally
the moduli space of vector bundles on them, 
have recently been used by Nakajima in his geometric construction of 
representations of $\widehat{\mathfrak{g}}$ (see \cite{Nak}).

\vspace{.1in}

Now let $\mathfrak{g}$ be of type $D_4^{(1)}, E_6^{(1)}, E_7^{(1)}$ or
$E_8^{(1)}$. As shown in \cite{GL}, $\S 5.8$, we now have
$Coh(\X_{\mathfrak{g}})\simeq Coh_G(\mathcal{E})$
where $\mathcal{E}$ is an elliptic curve
and $G$ a finite subgroup of $\mathrm{Aut}(\mathcal{E})$. This suggests
the study of the ``noncommutative space'' $T^*\X_{\mathfrak{g}}$ 
and the moduli space of vector bundles on it in the spirit of Nakajima's
work. 

\vspace{.1in}

More generally, as explained to us by H. Lenzing, there exists,
for \textit{any}
weighted projective line $\X_{\g}$ with at least three marked points, a pair
$(X,G)$ consisting of a smooth projective curve $X$ and a finite group of
automorphisms $G \subset Aut(X)$ such that $Coh_G(X) \simeq 
Coh(\X_{\p,\l})$). Conversely, if $X$ is a smooth projective curve and $G$
a group of automorphisms of $X$ such that $X/G \simeq \mathbb{P}^1$ then
$Coh_G(X)$ is equivalent to the category of coherent sheaves on the weighted
projective line $\mathbb{X}_{\p,\l}$, where $\l$ and $\p$ are the ramification locus
and multiplicities respectively. 

\vspace{.2in}

\centerline{\bf{Acknowledgements.}}

\vspace{.1in}

I would like to thank John Duncan, Igor Frenkel, Michael Kapranov
and Markus Reineke for interesting discussions, as well
as B. Enriquez, L. Hille, D. Kussin, H. Lenzing and C. Ringel.  
I am grateful to Hadi Salmasian
for teaching me all I know about Farey sequences.
  
\section{Loop algebras and quantum groups}
\vspace{.2in}
\paragraph{\textbf{1.1.}}  Let $A=(a_{ij})_{i,j=1}^r$
 be an irreducible symmetric generalized Cartan matrix of simply laced type 
and let $(\h, 
\Pi, \Pi^\vee)$ be a realization of $A$ (see \cite{Kac}), where 
$\Pi=\{\alpha_1,\ldots,\alpha_r\} \subset \h^*$ and $\Pi^\vee=
\{\alpha_1^\vee,\ldots \alpha_r^\vee\} \subset \h$. Let $\h'=\mathrm{Span}_\C
\{\alpha_1^\vee, \ldots \alpha_r^\vee\}$ and let $\g=\n_- 
\oplus \h' \oplus \n_+$ be the associated complex Kac-Moody algebra, 
with generators 
$\{e_i, f_i\;|i=1,\ldots,r\}\cup \h'$
satisfying the following set of relations
\begin{equation*}
[\h',\h']=0,\qquad [e_i,f_j]=\delta_{ij}\alpha^\vee_i,
\end{equation*}
\begin{equation*}
[\alpha_j^\vee,e_i]=a_{ij} e_{\alpha_i}, \qquad [\alpha_j^\vee,f_i]=
-a_{ij} f_i,
\end{equation*}
\begin{equation*}
ad(e_i)^{1-a_{ij}}(e_j)=0, \qquad 
ad(f_i)^{1-a_{ij}}(f_{j})=0.
\end{equation*}
Let $(\,,\,)$ be a standard invariant bilinear form on $\g$ (see \cite{Kac}).
In particular, we have $(e_i,f_i)=1$ for $i=1,\ldots , r$. Let us further
denote by $\Delta\subset \h^*$ the root system of $\g$ and by $Q=\bigoplus_i
\Z\alpha_i$ the root lattice. We also let $\Delta^{im} \subset \Delta$
stand for the subset of imaginary roots.
The $\Z$-module $Q$ is equipped with the Cartan
bilinear form determined by the relation $(\alpha_i, \alpha_j)=a_{ij}$.

\vspace{.2in}

\paragraph{\textbf{1.2.}} Let
$\widehat{\g}=\g[t,t^{-1}] \oplus \C \mathbf{c}$ be the affinization of
$\g$, where $\mathbf{c}$ is a central element, and where the Lie bracket
is given by
$$[xt^n,yt^m]=[x,y]t^{n+m} + n\delta_{n,-m}(x,y)\textbf{c}.$$
The algebra $\widehat{\g}$ is naturally equipped with a $\widehat{Q}=
Q \oplus \Z \delta$-grading, where
$$\widehat{\g}[\alpha+l \delta]=\g[\alpha] t^l, \qquad \widehat{\g}[0]=
\h \oplus \C c, \qquad \widehat{\g}[l \delta]=\h t^l.$$
We extend the Cartan form to $\widehat{Q}$ by 
setting $(\delta, \alpha)=0$ for all $\alpha \in \widehat{Q}$. 
The root system of $\widehat{\g}$ is $\widehat{\Delta}=
\Z^* \delta \cup \{\Delta + \Z \delta \}$. We will say that 
a root
$\alpha \in \widehat{\Delta}$ is \textit{real} if $(\alpha,\alpha)=2$
and \textit{imaginary} if $(\alpha,\alpha)\leq 0$.

\vspace{.2in}

\paragraph{\textbf{1.3.}} We will consider certain extensions of the
Lie algebra $\widehat{\g}$.
Let us denote by $\mathfrak{L}\g$ the Lie algebra generated by
elements $h_{i,k}, e_{i,k}, f_{i,k}$ for $i=1,\ldots r$, $k \in \Z$ and 
$\textbf{c}$ subject to the following relations :
\begin{equation*}
[h_{i,k}, h_{j,l}]=k \delta_{k,-l}a_{ij}\textbf{c},
\qquad [e_{i,k},f_{j,l}]=\delta_{i,j}h_{i,k+l}+ k \delta_{k,-l} \textbf{c},
\end{equation*}
\begin{equation*}
[h_{i,k}, e_{j,l}]=a_{ij}e_{j,l+k},\qquad [h_{i,k}, f_{j,l}]=
-a_{ij}f_{j,k+l},
\end{equation*}
\begin{equation*}
[e_{i,k+1},e_{j,l}]=[e_{i,k},e_{j,l+1}],\qquad
[f_{i,k+1},f_{j,l}]=[f_{i,k},f_{j,l+1}]
\end{equation*}
\begin{equation*}
\begin{split}
&[e_{i,k_1},[e_{i,k_2},[\ldots [e_{i,k_n},e_{j,l}]\cdots ]=0 \qquad 
\mathrm{if\;} n=1-a_{ij},\\
&[f_{i,k_1},[f_{i,k_2},[\ldots [f_{i,k_n},f_{j,l}]\cdots ]=0 \qquad 
\mathrm{if\;} n=1-a_{ij}.
\end{split}
\end{equation*}
It is clear that $e_i \mapsto e_{i,0}, f_i \mapsto f_{i,0}, \alpha_i^\vee
\mapsto h_{i,0}$ defines an embedding $\g \subset \mathfrak{L}\g$. Moreover,
giving generators $e_{i,k}, f_{i,k}, h_{i,k}$ degrees $\alpha_i + k\delta$,
$-\alpha_i + k\delta$ and $k\delta$ respectively endows $\mathfrak{L}\g$
with a $\widehat{Q}$-grading.
\paragraph{} There is
a surjective Lie algebra morphism 
$\varphi:\mathfrak{L}\g \to \widehat{\g}$ given by
by $e_{i,k} \mapsto e_{i}t^k, \; 
f_{i,k} \mapsto f_{i}t^k,\; h_{i,k} \mapsto \alpha_i^\vee t^k,
\mathbf{c} \mapsto \mathbf{c}$. It is proved in \cite{MRY} that
$$Ker\;\varphi \subset \bigoplus_{\underset{k \in \Z}{\alpha \in \Delta^{im}}}
\mathcal{L}\g[\alpha+k\delta].$$
In particular, the root system of $\mathfrak{L}\g$ is also $\widehat{\Delta}$
and the weight spaces corresponding to real roots are one-dimensional.

\paragraph{}When $\g$ is a simple complex Lie algebra, it is known that
$\varphi$ is an isomorphism $\mathfrak{L}\g \simeq \widehat{\g}$. Now suppose
that $\g=\g_0[s,s^{-1}] \oplus \C \mathbf{c}_0=\widehat{\g_0}$ is 
an (untwisted) affine algebra.
According to \cite{MRY}, $\mathfrak{L}\g$ is 
then the universal
central extension of $\g_0[s,s^{-1},t,t^{-1}]$ (the toroidal algebra, or
elliptic algebra). Let $\delta_0 \in \Delta$ stand for the imaginary root
of $\widehat{\g_0}$, and let $\Delta_0 \subset \Delta$ be the root system of
$\g_0$, so that $\Delta=\Z^*\delta_0 \cup (\Delta_0 + \Z\delta_0)$, and
$$\widehat{\Delta}=(\Z \delta_0 + \Z \delta)\cup (\Delta_0 + \Z \delta_0
+\Z\delta) \backslash\{0\}.$$
We have (see \cite{MRY}, Proposition 3.6 )
\begin{equation}\label{E:multcent}
\mathrm{dim}\;\mathfrak{L}\g[\alpha + k_0 \delta_0 + k \delta]
=
\begin{cases}
1 & \qquad \mathrm{if\;} \alpha \neq 0\\
l+1 & \qquad \mathrm{if\;} \alpha=0, (k_0,k) \neq (0,0)\\
l+2 & \qquad \mathrm{if\;} \alpha=k_0=k=0
\end{cases}
\end{equation}

\vspace{.2in}

\paragraph{\textbf{1.4.}} We will now consider certain nilpotent subalgebras
in $\mathfrak{L}\g$ and $\widehat{\g}$. Let $p_1, \ldots, p_n \in \N$
and let us consider the
``star'' Dynkin diagram
$$
\centerline{\epsfbox{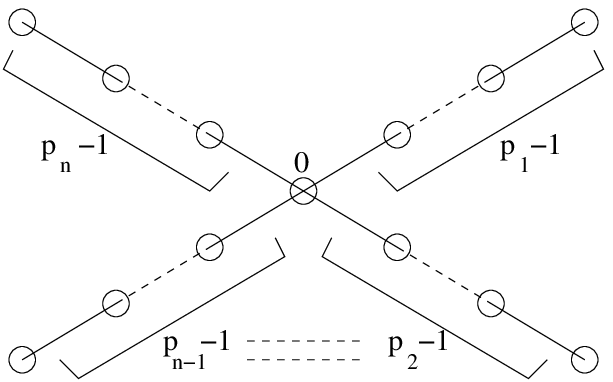}}
$$

Let us denote by $\star$ the central node and let 
$J_1, \ldots, J_{n}$ be the 
branches, which are subdiagrams are of type $A_{p_1-1},
\ldots A_{p_n-1}$ respectively. We can assume that $p_i >1$ for all $i$.
These examples include all finite type Dynkin 
diagrams as well as the affine Dynkin diagrams of types $D_4^{(1)}, E_6^{(1)},
E_7^{(1)}$, and $E_8^{(1)}$.

\vspace{.1in}
Let $\g$ be the Kac-Moody algebra corresponding to $\Gamma$, and
let $\mathfrak{p}$ be the maximal parabolic subalgebra of $\g$ associated
to the node $\star$. Denote by $\mathfrak{l}$ and $\mathfrak{m}$ its standard
Levi and nilpotent radical respectively. Let us write $\mathfrak{l}=
\mathfrak{n}^+_{\mathfrak{l}} \oplus \h \oplus \mathfrak{n}^-_{\mathfrak{l}}$
for the standard Cartan decomposition of $\mathfrak{l}$ and set 
$$\widehat{\mathfrak{t}}':= \mathfrak{n}^+_{\mathfrak{l}}\oplus t 
\mathfrak{l}[t] \subset \mathfrak{l}[t,t^{-1}],
\qquad \widehat{\mathfrak{f}}':= \mathfrak{m}[t,t^{-1}].$$
Finally, we put 
\begin{equation}\label{E:044}
\widehat{\n}':=\widehat{\mathfrak{f}}' \oplus \widehat{\mathfrak{t}}'\subset
\widehat{\g}.
\end{equation}

\vspace{.1in}
The following description of $\widehat{\n}'$ will be convenient for us.
For every $s=1, \ldots, n$ let $\g^{s}$ be the 
subalgebra of $\g$ generated by the elements 
$\{e_j, f_j, \alpha_j^\vee \;|j \in J_s\}$. Let 
$\theta_{s}$ be the highest root of $\g^s$ and let 
$f_{\theta_{s}} \in \g^s_{-\theta_s}$ be 
any lowest root vector. 
Let $\widehat{\n}^s \subset \widehat{\g}$ be generated by
$\{e_j \;|j \in J_s\} \cup \{f_{\theta_{s}}t\}$. Then
$\widehat{\n}'$ is the subalgebra generated by
\begin{equation}\label{E:generators}
\bigcup_{s=1}^n \widehat{\n}^s
\cup \{e_{\star}t^l\;| l \in \Z\} \cup \{h_{\star}t^r\;|r \in \N^*\}.
\end{equation}
and $\widehat{\mathfrak{t}}'$ is the subalgebra
generated by $\widehat{\n}^s$ for $s=1, \ldots, n$ and $h_{\star}t^r$
for $r \geq 1$.

\vspace{.1in} 

\paragraph{}Observe that, for any $i$, the subalgebra of 
$\mathfrak{L}\g$ generated
by $e_{j,l}, f_{j,l}, h_{j,l}$ for $j \in J_i$ and $l \in \Z$ is canonically
isomorphic (via the map $\varphi$) to an affine algebra of type $A_{p_i-1}$,
and we can define the subalgebra $\widehat{\n}^s \subset \mathfrak{L}\g$
in the same way as above. We now define $\mathfrak{L}\n$ to be the subalgebra
of $\mathfrak{L}\g$ generated by
\begin{equation}\label{E:generators2}
\bigcup_{s=1}^n \widehat{\n}^s
\cup \{e_{\star,l}\;| l \in \Z\} \cup \{h_{\star,r}\;|r \in \N^*\}.
\end{equation}


\vspace{.2in}
\paragraph{\textbf{1.5.}} Let $p >0$. The Lie algebra $\mathfrak{gl}_p$
of $p$ by $p$ matrices has a natural basis 
$$\{E_{ij}\;|i,j=1, \ldots p\}$$ such that
$[E_{ij},E_{kl}]=\delta_{jk}E_{il}-\delta_{il}E_{kj}$. Set $h_0=-E_{11}$. Then 
$$\mathfrak{gl}_p=\mathfrak{sl}_p \oplus \C h_0.$$
Let $(x,y)=\mathrm{Tr}(xy)$ be the (normalized) Killing form.  
Define $\widehat{\mathfrak{gl}}_p=\mathfrak{gl}_p[t,t^{-1}] \oplus 
\C\mathbf{c}$ as in Section 1.2. Set $e_i=E_{i,i+1}$ for $i=1, \ldots, 
p-1$ and $e_0=E_{p1}\otimes t$.
Let $\widehat{\mathfrak{sl}}_p^+
\subset
\widehat{\mathfrak{sl}}_p$ (resp. $\widehat{\mathfrak{gl}}_p^+
\subset
\widehat{\mathfrak{gl}}_p$) be the subalgebra generated by
$e_i$, $i=0, \ldots, p-1$ (resp.
generated by $e_i$, $i=1, \ldots, p-1$ and $\mathfrak{gl}_p\otimes 
t\C[t]$).
In particular,
\begin{equation}\label{E:132}
\U(\widehat{\mathfrak{gl}}_p^+)=\U(\widehat{\mathfrak{sl}}_p^+) \otimes
\C[h_1,h_2,\ldots],
\end{equation}
where we set $h_i=h_0t^i$.

\vspace{.2in}

\paragraph{\textbf{1.6.}} Let $\Gamma$ be as in Section 1.4. 
We give a presentation of (a certain extension of) the Lie algebra
$\widehat{\n}'$ by generators and 
relations.
For each $s$ we let $\phi_s: A_{p_s-1} \stackrel{\sim}{\to} J_s$
be the unique bijection such that $\phi_s(1)$ is adjacent to the
central node $\star$.

\vspace{.1in}

\paragraph{\textbf{Definition.}} Let $\widehat{\n}$ be the Lie
algebra generated by $\varepsilon_j$ for $j \in \Gamma \backslash \{\star\}$,
$\varepsilon_0^{(s)}$ for $s=1, \ldots, n$, and $\varepsilon_{\star,t}$, 
$\theta_{\star,r}$,
for $r \geq 1$ and $t \in \Z$, subject to the following set of relations~:
\begin{enumerate}
\item[i)] For all $s$, the assignment $e_j \mapsto \varepsilon_{\phi_s(j)}$,
$e_0 \mapsto \varepsilon_0^{(s)}$ and $h_r \mapsto \theta_{\star,r}$
extends to an embedding $\widehat{\mathfrak{gl}}_{p_s}^+ 
\hookrightarrow \widehat{\n}$, and $[\widehat{\mathfrak{gl}}_{p_s}^+,
\widehat{\mathfrak{gl}}_{p_{s'}}^+]=0$ if $s \neq s'$.
\item[ii)] We have 
\begin{equation}
[\theta_{\star,r},\varepsilon_{\star,t}]=2\varepsilon_{\star,r+t}, \qquad 
[\varepsilon_{\star,t},\varepsilon_{\star,t'}]=0
\end{equation}
\item[iii)] 
\begin{equation}\label{E:134}
a_{\star j}=0 \Rightarrow [\varepsilon_{\star,t},\varepsilon_j]=0,
\qquad [\varepsilon_{\star,t},\varepsilon_0^{(s)}]=0,\;\; s=1, \ldots ,n
\end{equation}
\begin{equation}\label{E:133}
p_s >2 \Rightarrow [\varepsilon_0^{(s)},[\varepsilon_{\star,t},\varepsilon_{\phi_s(1)}]]=0.
\end{equation}
\item[iv)] For $s=1, \ldots, n$ and $j=\phi_s(1)$ set
$\varepsilon_{j,r}=[\varepsilon_j,\theta_{\star,r}]$. Then for any $t, t_1, t_2 \in \Z$
and $r,r_1,r_2 \geq 1$,
\begin{equation}
[\varepsilon_{\star,t_1},[\varepsilon_{\star,t_2},\varepsilon_{j,r}]]=[\varepsilon_{j,r_1},[\varepsilon_{j,r_2},\varepsilon_{\star,t}]]=0
\end{equation}
\begin{equation}\label{E:135}
[\varepsilon_{\star,t},\varepsilon_{j,r+1}]=[\varepsilon_{\star,t+1},\varepsilon_{j,r}]
\end{equation}
\end{enumerate}
\paragraph{}It is easy to check that the assignement $\varepsilon_{\star,l} \mapsto 
e_{\star,l},
\theta_{\star,k} \mapsto h_{\star,k}, \varepsilon_j \mapsto e_j, \varepsilon_0^{(s)} \mapsto 
f_{\theta_s}t$ extends to a surjective homomorphism $\phi: \widehat{\n}
\to \mathfrak{L}{\n}$. Moreover, $\widehat{\n}$ is $\widehat{Q}$-graded.

\begin{prop} We have
 $$Ker \;\phi \subset \bigoplus_{\underset{l \in \Z}{\alpha \in \Delta^{im}}}
\widehat{\n}[\alpha+l\delta]$$
Moreover, $\phi$ is an isomorphism if $\Gamma$ is a finite or an affine 
Dynkin diagram.
\end{prop}
\noindent
\textit{Proof.} For $s=1, \ldots
,n$ we denote by
$\n^s_- \subset \mathfrak{sl}_{p_s}$ the subspace spanned by $E_{k,l}$
for $k>l$. 
\begin{lem} For any $s$, any $t \in \Z$ and any 
$x \otimes t^r \in \n^s_-\otimes t\C[t]
\subset \widehat{\mathfrak{gl}}_{p_s}^+$ we have $[\varepsilon_{\star,t},x\otimes
t^r]=0$.\end{lem}
\noindent
\textit{Proof.}
Set $\varepsilon_{0,r}^s=[\theta_{\star,r},\varepsilon_0^{(s)}]$. If $p_s=2$ then the lemma 
follows from (\ref{E:134}). If $p_s>2$, observe that
$[\varepsilon_{0,r},\varepsilon_{\star,t}]=0$ and that
\begin{equation}\label{E:136}
\begin{split}
[\varepsilon_{\star,t},[\varepsilon_{0,r},\varepsilon_{1}]]&=[\varepsilon_{0,r},[\varepsilon_{\star,t},\varepsilon_1]]\\
&=[[\theta_{\star,r},\varepsilon_0],[\varepsilon_{\star,t},\varepsilon_1]]\\
&=[\theta_{\star,r},[\varepsilon_0,[\varepsilon_{\star,t},\varepsilon_1]]] + [\varepsilon_{0}, [\theta_{\star,r},[\varepsilon_{\star,t},
\varepsilon_1]]]\\
&=2[\varepsilon_0,[\varepsilon_{\star,t+r},\varepsilon_1]] -[\varepsilon_0,[\varepsilon_{\star,t},\varepsilon_{1,r}]]\\
&=0
\end{split}
\end{equation}
where we have used (\ref{E:133}) and (\ref{E:135}). Any 
$x \otimes t^r \in \n^s_-\otimes t\C[t]$ is a linear combination of
terms of the form $[\varepsilon_{0,r},[\varepsilon_{i_1},[\varepsilon_{i_2}, \cdots ]]$ where
$i_1,i_2 \ldots \in \{1, \ldots, p_s-1\}$ are distinct. The lemma
is now a consequence of (\ref{E:136}) . \qed

\vspace{.1in}

\paragraph{}The algebra $\widehat{\n}$ is linearly spanned by
iterated commutators
$$[\alpha_1, \ldots, \alpha_r]=[\alpha_1, [\alpha_2,[ \ldots , \alpha_r] \cdots
]]$$
with 
$\alpha_i \in \{e_{\star,l}, e_j t^k, h_jt^k,xt^r\;|\; x \in \n^s_-; s=1, \ldots, n; l \in \Z; k \geq 0; r \geq 1\}$. Set 
$$\widehat{\n}_{\pm}=\bigoplus_{\alpha \in \Delta^{\pm}, l \in \Z} 
\widehat{\n}[\alpha + l \delta], \qquad (\mathfrak{L}\n)_{\pm}
=\bigoplus_{\alpha \in \Delta^{\pm}, l \in \Z} 
\mathfrak{L}\n[\alpha + l \delta].$$

From the defining relations and the Lemma above
one easily deduces that $\widehat{\n}_+$ is generated by 
$\varepsilon_{\star, l}, \varepsilon_jt^k$ for $l \in \Z$, 
$j \in \Gamma \backslash \{\star\}$ and $k \geq 0$. Similarly,
$\widehat{\n}_-$ is generated by $\{xt^r|\; x \in \n^s_-, r \geq 1, s=1, \ldots
, n\}$. It is clear that $Ker\;\phi_{|\widehat{\n}_-}=0$. On the other hand,
define $\widehat{\n}_{E,+}$ to be the Lie algebra generated by elements
$\epsilon_{\star,l}, \epsilon_{j,k}$ for $l \in \Z$, $j \in \Gamma \backslash
\{\star\}$ and $k \geq 0$ subject to the relations
$$
[\epsilon_{i,k+1},\epsilon_{j,l}]=[\epsilon_{i,k},\epsilon_{j,l+1}] \qquad \forall \;k,j, \; \forall \;
i,j \in \Gamma,$$
$$[\epsilon_{i,k_1},[\epsilon_{i,k_2},[\ldots [\epsilon_{i,k_n},\epsilon_{j,l}]\cdots ]=0 \qquad 
\mathrm{if\;} n=1-a_{ij}, \; \forall\; k_1, \ldots, k_n, l.$$
There is an obvious surjective map 
$\psi:\widehat{\n}_{E,+} \to \widehat{\n}_+$. But from \cite{E}, 
Proposition 1.5 it follows that 
$$Ker\;\phi \psi \subset \bigoplus_{\alpha \in \Delta^{im}, l \in \Z}
\widehat{\n}_{E,+}[\alpha+l\delta],$$
and from \cite{E}, Proposition 1.6 we deduce that $\phi \psi$ is injective 
when $\Gamma$ is finite or affine. This proves the Proposition. \qed

\vspace{.1in}

\paragraph{}Observe that if $\Gamma$ is finite then 
$\widehat{\n}=\mathfrak{L}\n=\widehat{\n}'$. 
Let us now assume that $\Gamma$ is affine. Let $\widehat{\mathfrak{t}} \subset
\widehat{\n}$ be the subalgebra generated by $\mathfrak{gl}_{p_i}^+$
for $i=1, \ldots n$, and let $\widehat{\mathfrak{f}} \subset \widehat{\n}$ be 
its (unique) $\widehat{Q}$-graded complement. It is clear that 
$\widehat{\mathfrak{t}}
\simeq \widehat{\mathfrak{t}}'$. We will say that a root 
$\alpha =\sum_{i \in \Gamma}
\lambda_i \alpha_i$ is $\star$-positive if $\lambda_{\star} >0$.

\begin{cor} Assume that $\Gamma$ is affine. Then 
\begin{equation}\label{E:multcent2}
\mathrm{dim}\;\widehat{\mathfrak{f}}[\alpha + k_0 \delta_0 + k \delta]
=
\begin{cases}
1 & \qquad \mathrm{if\;} \alpha \in \Delta_0^+ \text{\;and\;}
 \alpha + k_0 \delta_0\text{\;is\;}\star-\text{positive}\\
l+1 & \qquad \mathrm{if\;} \alpha=0,\;k_0 >0\\
0 & \qquad \mathrm{otherwise}
\end{cases}
\end{equation}
where $l=|\Gamma|-1$ is the rank of $\g_0$.
\end{cor}

\vspace{.2in}

\paragraph{\textbf{1.7.}} For any $n \in \N$ we set
$$[n]=\frac{v^n-v^{-n}}{v-v^{-1}},\qquad [n]!=\prod_{i=1}^n [i]!, \qquad
\bmatrix n\\ m \endbmatrix = \frac{[n]!}{[n-m]![m]!},$$
and put $\C_v=\C[v,v^{-1}]$. Let $\mathbf{U}$ be the $\C(v)$-algebra
generated by 
$\{E_i, F_i, K_i, K_i^{-1}\;|i=1,\ldots , r\}$ 
satisfying the following relations
$$K_i K_i^{-1}=K_i^{-1}K_i=1,\qquad [K_i, K_j]=0,$$
$$K_iE_{j}K_i^{-1}=v^{a_{ij}}E_{j},\quad
K_iF_{j}K_i^{-1}=v^{-a_{ij}}F_{j}$$
$$E_{i}F_{j}-F_{j}E_{i}=\delta_{ij}\frac{K_i-K_i^{-1}}{v -v^{-1}},$$
$$\sum_{k=0}^{1-a_{ij}} (-1)^k \bmatrix 1-a_{ij} \\ k
\endbmatrix E_{i}^{1-a_{ij}-k}E_{j}E_{i}^k=0,\quad i\ne j,$$
$$\sum_{k=0}^{1-a_{ij}} (-1)^k \bmatrix 1-a_{ij} \\ k
\endbmatrix F_{i}^{1-a_{ij}-k}F_{j}F_{i}^k=0,\quad i\ne j.$$
The quantum enveloping algebra $\U_v(\g)$
is by definition the $\C_v$-subalgebra of $\mathbf{U}$ generated by
$K_i, K_i^{-1}$ for $i =1, \ldots, r$ and the divided powers
$$E_i^{(n)}:=\frac{E_i^n}{[n]!},\qquad F_i^{(n)}:=\frac{F_i^n}{[n]!}$$
for $i=1, \ldots ,r$ and $n \geq 1$. As usual, we let 
$\mathbf{U}^+_v(\mathfrak{g})$ be the $\C_v$-subalgebra of $\U_v(\g)$
generated by $E_i^{(n)}$ for $i=1, \ldots,r$ and $n \geq 1$.
\vspace{.2in}
\paragraph{\textbf{1.8.}} The quantum loop algebra (with zero central charge)
$\U_v(L\g)$ is the 
$\C(v)$-algebra generated by $x_{i,k}^{\pm}, h_{i,l}$ and $K_i^{\pm}$ 
for $i=1, \ldots r$,
$k \in \Z$ and $l \in \Z^*$ subject to the following relations
$$K_i K_i^{-1}=K_i^{-1}K_i=1,\qquad [K_i,K_j]=[K_i,h_{j,l}]=
[h_{i,l},h_{j,k}]=0 ,$$
$$K_i x_{jk}^{\pm}K_{i}^{-1}=v^{\pm a_{ij}}x_{jk}^{\pm},$$
$$[h_{i,l},x_{j,k}^{\pm}]=\pm \frac{1}{l}[la_{ij}]x^{\pm}_{j,k+l},$$
$$x_{i,k+1}^{\pm}x_{j,l}^{\pm}-v^{\pm a_{ij}}x^{\pm}_{j,l}x^{\pm}_{i,k+1}=
v^{\pm a_{ij}}x_{i,k}^{\pm}x_{j,l+1}^{\pm}-x_{j,l+1}^{\pm}x_{i,k}^{\pm}$$
$$[x_{i,k}^+,x_{j,l}^-]=\delta_{ij}
\frac{\psi_{i,k+l}-\varphi_{i,k+l}}{v-v^{-1}},$$
\begin{equation*}
\begin{split}
\mathrm{For\;}i \neq j&\;\mathrm{and\;}n=1-a_{ij},\\
&\mathrm{Sym}_{k_1, \ldots, k_n} \sum_{t=0}^{1-a_{ij}}(-1)^t 
\bmatrix n & t \endbmatrix x_{i,k_1}^{\pm} \cdots x_{i,k_t}^{\pm}x_{j,l}^{\pm}
x^{\pm}_{i,k_{t+1}} \cdots x^{\pm}_{i,k_n}=0
\end{split}
\end{equation*}
where $\mathrm{Sym}_{k_1, \ldots, k_n}$ denotes symmetrization with respect
to the indices $k_1, \ldots , k_n$, and where $\psi_{i,k}$ and $\varphi_{i,k}$
are defined by the following equations :
$$\sum_{k \geq 0} \psi_{i,k}u^k=K_i \mathrm{exp}\big( (v-v^{-1})
\sum_{k=1}^\infty h_{i,k}u^k\big),$$ 
$$\sum_{k \geq 0} \varphi_{i,k}u^k=K_i^{-1} \mathrm{exp}\big( -(v-v^{-1})
\sum_{k=1}^\infty h_{i,-k}u^{-k}\big).$$
\vspace{.2in}
\paragraph{\textbf{1.9.}} The quantum loop algebra $\U_v(L\mathfrak{gl}_n)$
is the $\C(v)$-algebra generated by $U_v(L\mathfrak{sl}_n)$ and generators
$K_0, K_0^{-1}, h_{0,l}$ for $l \in \Z^*$ with relations
$$K_0 K_0^{-1}=K_0^{-1}K_0=1,$$
$$[K_0,K_i]=[K_0,h_{i,l}]=[K_0,h_{0,l}]=
[h_{0,l},h_{i,k}]=[h_{0,l},h_{0,k}]=0,$$
$$K_0x_{i,k}^{\pm}K_0^{-1}=v^{\mp\delta_{i1}}x_{i,k}^{\pm},$$
$$[h_{0,l},x^{\pm}_{i,k}]=\mp \delta_{i,1}\frac{[l]}{l}x^{\pm}_{i,l+k}.$$
\vspace{.1in}
\section{Weighted projective lines}
\vspace{.2in}

\paragraph{}We recall in this section the definition and main properties of
the weighted projective lines $\X_{\p,\l}$ and the categories of coherent 
sheaves on them, as studied in \cite{GL}. 
\vspace{.2in}
\paragraph{\textbf{2.1.}} Let $k$ be an arbitrary field and let 
$\p=\{p_1,\ldots, p_n\} \in (\N^*)^n$. Consider the $\Z$-module
$$L(\p)=\Z \x_1 \oplus \cdots \oplus \Z \x_n / J$$
where $J$ is the $\Z$-submodule generated by $\{p_1 \x_1 -p_s \x_s\;|s =2, 
\ldots ,n\}$. Set $\c=p_1\x_1 = \cdots = p_n \x_n$. Then $L(\p)/\Z \c \simeq
\prod_{s=1}^n \Z/p_s\Z$ and we can consider $L(\p)$ as a rank 
one abelian group. We set $L^+(\p)=\{ l_1 \x_1 + \cdots + l_n \x_n\;|l_i \geq 0
\}$ and we write $\x \geq \vec{y}$ if $\x-\vec{y} \in L^+(\p)$.

Let $k[L(\p)]$ be the group Hopf algebra of $L(\p)$ and let
$$G(\p)=\mathrm{Spec}\;k[L(\p)] \subset (k^*)^n$$
be the associated affine algebraic group. Hence the $\overline{k}$-points of
$G(\p)$ are 
$$G(\p)(\overline{k})= \{(t_1, \ldots , t_n) \in (\overline{k}^*)^n\;|
t_1^{p_1}=\cdots = t_n^{p_n}\}.$$

By definition $G(\p)$ acts on $k^n$ and induces on $S(\p):=k[X_1,\ldots , X_n]$
the structure of an $L(\p)$-graded algebra such that $\mathrm{deg}\;X_s=\x_s$.
\paragraph{}Let $\l=\{\lambda_1,\ldots , \lambda_n\}$ be a collection of 
distinct closed points (of degree $1$) of $\mathbb{P}^1(k)$, normalized in 
such a way that
$\lambda_1=\infty$, $\lambda_2=0$ and $\lambda_3=1$. Let $I(\p,\l)$ be the
$L(\p)$-graded ideal in $S(\p)$ generated by the polynomials $X_s^{p_s}-
(X_2^{p_2}-
\lambda_s X_1^{p_1})$ for $s=3, \ldots , n$ and set $S(\p,\l)=S(\p)/I(\p,\l)$.
We will denote by $x_1,\ldots , x_n$ the images of $X_1, \ldots , X_n$ in 
$S(\p,\l)$. It is known that $S(\p,\l)$ is a graded factorial domain of 
Krull dimension equal to $2$. We let
$$\X_{\p,\l}=\mathrm{Specgr}\;S(\p,\l)$$
be the set of prime homogeneous ideals in $S(\p,\l)$.
\begin{prop}[\cite{GL}] The nonzero homogeneous prime elements in $S(\p,\l)$
are of the form
\begin{enumerate}
\item[i)] $x_1,\ldots, x_n$ (\textit{the\;exceptional\;points}), or
\item[ii)] $F(x_1^{p_1},x_2^{p_2})$ where $F \in k[T_1,T_2]$ is a prime
homogeneous polynomial (\textit{ordinary\;points}).
\end{enumerate}
\end{prop}
We will denote by $\sigma_s$, $s=1, \ldots, n$ the closed points of 
$\X_{\p,\l}$ corresponding to the primes $x_1, \ldots, x_n$. 
Let $p=\mathrm{l.c.m}(p_1, \ldots, p_n)$. We define the \textit{degree} of a 
closed point $x$ by $deg(\sigma_i)=\frac{p}{p_i}$ for exceptional points
and $deg(x)=pd$ if $x$ is associated to a prime homogeneous polynomial
$P(x_1^{p_1},x_2^{p_2})$ of degree $d$.

For any homogeneous $f \in S(\p,\l)$ we set $V_f=\{\mathfrak{p}
\in \X_{\p,\l}\;|f \in \mathfrak{p}\}$ and $D_f=\X_{\p,\l} \backslash V_f$.
The sets $D_f$ form a basis of the Zariski topology on $\X_{\p,\l}$.
\vspace{.2in}
\paragraph{\textbf{2.2.}} Let $\O_\X$ be the sheaf of \textit{$L(\p)$-graded}
algebras
on $\X_{\p,\l}$ associated to the presheaf $D_f \mapsto S(\p,\l)_f$ where
$$S(\p,\l)_f=\{\frac{g}{f^l}\;|g \in S(\p,\l),\; l \in \N\}$$
is the localization of $S(\p,\l)$ at $f$. Let us denote by $\O_\X$-Mod the
category of sheaves of \textit{$L(\p)$-graded} $\O_\X$-modules on $\X_{\p,\l}$.

\paragraph{}For any $\vec{x} \in L(\p)$ and any object 
$\mathcal{M}=\bigoplus_{\vec{l} \in L(\p)} \mathcal{M}[\vec{l}]$ in 
$\O_\X$-Mod we denote by $\mathcal{M}(\vec{x})$ the shift of $\mathcal{M}$ by
$\vec{x}$ (i.e $\mathcal{M}(\vec{x})[\vec{l}]:=\mathcal{M}[\vec{x}+\vec{l}]$). The category $Coh(\X_{\p,\l})$ is now defined as the full subcategory of
$\O_\X$-Mod
consisting of sheaves $\mathcal{M}$ for which there exists an open covering
$\{U_i\}$ of $\X_{\p,\l}$ and, for each $i$, an exact sequence
$$\bigoplus_{s=1}^N \O_\X(\vec{l}_s)_{|U_i} \to \bigoplus_{t=1}^M \O_\X
(\vec{k}_t)_{|U_i} \to \mathcal{M}_{|U_i} \to 0.$$

Note that for any $\mathcal{M} \in Coh(\X_{\p,\l})$ the space 
$\Gamma(\X_{\p,\l}, \mathcal{M})$ of global sections of $\mathcal{M}$ is an
$L(\p)$-graded $S(\p,\l)$-module.

\begin{prop}[\cite{GL}] The category $Coh(\X_{\p,\l})$ is
hereditary, i.e
$$\mathrm{Ext}^2(\mathcal{M},\mathcal{N})=0$$
for any $\mathcal{M},\mathcal{N}$.
Moreover, for any $\mathcal{M},\mathcal{N}$, $Hom(\mathcal{M},\mathcal{N})$
and $Ext^1(\mathcal{M},\mathcal{N})$ are finite-dimensional.\end{prop}

We now define the \textit{support} of a sheaf. Let $x$ be any closed point
of $\X_{\p,\l}$. Let us denote by $\mathfrak{m}_x$ the corresponding maximal
ideal and set
$$\O_{\X ,x}=\{\frac{f}{g}\;|f, g \in S(\p,\l), g\;\mathrm{homogeneous},
g \not\in \mathfrak{m}_x\}.$$
This is an $L(\p)$-graded discrete valuation ring. If $\mathcal{M} \in 
Coh(\X_{\p,\l})$ we define $\mathcal{M}_x=\mathcal{M}(D_f) \otimes_{\O_\X(D_f)}
\O_{\X ,x}$ where $D_f$ is any neighborhood of $x$. The support of 
$\mathcal{M}$ is the Zariski closure of the set 
$\{x \in \X_{\p,\l}\;|\mathcal{M}_x\neq 0\}$.

Let $\xi \in \X_{\p,\l}$ denote the generic point, associated to the prime
ideal $(0)$. Set
$$\O_{\X,\xi}=\{\frac{f}{g}\;|f,g \in S(\p,\l), g \; homogeneous, g \neq 0\}.$$
This $L(\p)$-graded ring is Morita equivalent to its degree zero component
$\O_{\X,\xi}[0] \simeq k(\frac{x_1^{p_1}}{x_2^{p_2}})$. The \textit{rank}
of a coherent sheaf $\mathcal{M}$ is defined as
$$\mathbf{r}(\mathcal{M})=\mathrm{dim}\; \mathcal{M}(D_f) \otimes_{\O_\X(D_f)}
\O_{\X,\xi}$$
for any $D_f$. 
\vspace{.2in}
\paragraph{\textbf{2.3.}} Let $S(\p,\l)$-modgr be the category of finitely
generated $L(\p)$-graded $S(\p,\l)$-modules. As in the classical case, one
can introduce a sheafification functor 
\begin{align*}
Sh:\; S(\p,\l)-\mathrm{modgr} &\to Coh(\X_{\p,\l})\\
M &\mapsto \tilde{M}
\end{align*}
where $\tilde{M}$ is the sheaf associated to the presheaf
$$D_f \mapsto M_f=\{\frac{m}{f^l}\;|m \in M, l \in \N\}.$$

Let $L^+(\p) \subset L(\p)$ be the monoid generated by $\x_1, \ldots \x_n$.
Denote by $S^+(\p,\l)$-modgr be the full subcategory of $S(\p,\l)$-modgr
consiting of $L^+(\p)$-graded modules. Consider the global sections functor
\begin{align*}
\Gamma_+:\;Coh(\X_{\p,\l}) &\to S^+(\p,\l)-\mathrm{modgr}\\
\mathcal{M} &\mapsto \bigoplus_{l \geq 0} \Gamma(X,\mathcal{M})[\vec{l}].
\end{align*}

The following analogue of Serre's theorem holds :
\begin{prop}[\cite{GL}] The functor $\Gamma_+$ defines a full embedding and
for any $\mathcal{M} \in Coh(\X_{\p,\l})$ we
have $Sh(\Gamma_+(\mathcal{M})) \simeq \mathcal{M}$.
 The functor $Sh$ is exact and induces an equivalence
$$Sh:\;S(\p,\lambda)-\mathrm{modgr}/S(\p,\l)-\mathrm{modgr}_0 
\stackrel{\sim}{\to} Coh(\X_{\p,\l}),$$
where $S(\p,\l)-\mathrm{modgr}_0$ denotes the quotient of $S(\p,\l)$-modgr
by the Serre subcategory of finite-length modules.
\end{prop}
In particular, we have
\begin{equation}\label{E:23}
\mathrm{Hom}(\O_{\X}(\vec{x}), \O_{\X}(\vec{y}))=S(\p,\l)[\vec{y}-\vec{x}].
\end{equation}
\vspace{.1in}
\paragraph{\textbf{2.4.}} A sheaf $\mathcal{M}$ is \textit{locally free}
or is a \textit{vector bundle}
if there exists an open covering $\{U_i\}$ of $\X_{\p,\l}$ and, 
for each $i$, an isomorphism $\bigoplus_{i=1}^N \O_{\X}(\vec{l}_i)_{|U_i} 
\stackrel{\sim}{\to}\mathcal{M}_{|U_i}$ for some suitable $\vec{l}_i$. Every
rank one locally free sheaf on $\X_{\p,\l}$ is isomorphic to $\O_{\X}(\vec{x})$
for
some $\vec{x} \in L(\p)$. The full subcategory consisting of locally free 
sheaves is stable under extensions (but not under quotient, and is not 
abelian).

\vspace{.1in}

A sheaf $\mathcal{M}$ is a \textit{torsion sheaf} if it 
is a finite length 
object in $Coh(\X_{\p,\l})$, or, equivalently, if it is of rank zero.
By definition, the full subcategory $\mathcal{T}$
consisting of torsion sheaves is a Serre subcategory. The 
following properties can be found in \cite{GL} :
\begin{enumerate}
\item[i)] Every sheaf $\mathcal{M} \in Coh(\X_{\p,\l})$ can be decomposed
as $\mathcal{M} = \mathcal{M}_t \oplus \mathcal{M}'$ for a (unique) 
$\mathcal{M}_t \in \mathcal{T}$ and some $\mathcal{M}' \in \mathcal{F}$.
\item[ii)] $\mathrm{Hom}\;(\mathcal{M}_t, \mathcal{M}_f)=
\mathrm{Ext}^1\;(\mathcal{M}_f, \mathcal{M}_t)=0$ for 
every $\mathcal{M}_t
\in \mathcal{T}$ and $\mathcal{M}_f \in \mathcal{F}$.
\end{enumerate}
\vspace{.1in}

The category $\mathcal{T}$
decomposes as a product $\mathcal{T}=\sqcup_{x \in \X_{\p,\l}} \mathcal{T}_x$
where $\mathcal{T}_x$ is the category of torsion sheaves with support at
the closed point $x \in \X_{\p,\l}$. Moreover, $\mathcal{T}_x$ is
Morita-equivalent to the category $\mathrm{Mod}_0\O_{\X,x}$ of finite-length
$L(\p)$-graded $\O_{\X,x}$-modules.

These can be explicitely described.
Let $\mathcal{C}_l$ denote quiver of type $A^{(1)}_{l-1}$ with cyclic
orientation. We let $\mathcal{C}_1$ be the quiver with one vertex and one
arrow. A representation of $\mathcal{C}_p$ over a field $k$ is by definition a
collection of $k$-vector spaces $V_i$ with $i \in \Z/p\Z$ together with a 
collection of linear maps $x_i : V_i \to V_{i-1}$. Morphisms between two 
representations $(V_i, x_i)$ and $(W_i, y_i)$ are $k$-linear maps $\phi_i: V_i
\to W_i$ such that $\phi_ix_i=y_i\phi_i$. Finally, a representation $(V_i,x_i)$
is called nilpotent if there exists $N \gg 0$ such that $x_{i+N} \cdots x_i=0$
for all $i$.

\begin{lem} The following holds :
\begin{enumerate}
\item[i)] Let $x$ be an ordinary closed point
of degree $d$ (see Proposition 2.1 ) and let $k_x$ denote the residue field
at $x$. Then $\mathcal{T}_x$ is Morita-equivalent to the category
of nilpotent representation of the cyclic quiver $\mathcal{C}_1$ over $k_x$.\\
\item[ii)]  Let $1 \leq s \leq n$. The category 
$\mathcal{T}_{\sigma_s}$ is Morita equivalent to the category of nilpotent 
representations of $\mathcal{C}_{p_s}$ over $k$.
\end{enumerate}
\end{lem}

In case i) of the lemma above, the simple module $S_x$ can be realized 
in the following way : let $\pi_x \in S(\p,\l)$ denote the prime corresponding
to $x$. Multiplication by $\pi_x$ leads to an exact sequence
$$0 \to \O_\X \to \O_{\X}(d\vec{c}) \to S_x \to 0.$$
Moreover, for any $\vec{k} \in L(\mathbf{p})$ we have $S_x(\vec{k})\simeq S_x$.
In case ii), the simple sheaves $S^{(s)}_{j}$, $j =1, \ldots ,p_s$ with support
at $\sigma_s$ can be realized as follows :
multiplication by $x_s$ leads to exact sequences
$$0 \to \O_\X((j-1)\vec{x}_s) \to \O_{\X}(j\vec{x}_s) \to S^{(s)}_{j} 
\to 0.$$
Moreover, for any $\vec{k}=\sum_i k_i \x_i \in L(\mathbf{p})$ we have
\begin{equation}\label{E:299}
S_j^{(s)} (\vec{k}) \simeq S^{(s)}_{j + k_s\;\mathrm{mod}\;p_s}.
\end{equation}
In particular, $S^{(s)}_{j}(t\vec{c}) \simeq S_j^{(s)}$ for any $t \in \Z$.
Finally, we note the following result :
\begin{prop} Let $K_0(Coh(\X_{\p,\l}))$ be the Grothendieck group of
$Coh(\X_{\p,\l})$ and let $[\mathcal{M}]$ denote the class of a coherent
sheaf $\mathcal{M}$. Then
$$K_0(Coh(\X_{\p,\l}))\simeq \big(\Z[\O_{\X}] \oplus \Z[\O_{\X}(\vec{c})]
\bigoplus_{s,j} \Z[S^{(s)}_{j}]\big) / I$$
where $I$ is the subgroup generated by $\{\sum_{j} [S^{(s)}_{j}] + [\O_\X] -
[\O_\X(\vec{c})]\}_{s=1,\ldots, n}.$
\end{prop}

\vspace{.2in}

\paragraph{\textbf{2.5.}} Let $\vec{\omega}=(n-1)\vec{c} -\sum_i \vec{x}_i$.

\begin{theo}[Serre duality, \cite{GL}] Let $\mathcal{M}, \mathcal{N} \in 
Coh(\X_{\p,\l})$.
There is a functorial isomorphism
$$(\mathrm{Ext}^1(\mathcal{M}, \mathcal{N}))^* \stackrel{\sim}{\to}
\mathrm{Hom}(\mathcal{N},\mathcal{M}(\vec{\omega})).$$
\end{theo}
\vspace{.2in}

\paragraph{\textbf{2.6.}} Define a
group homomorphism $\partial:\;L(\p) \to \Z$ by setting $\partial(\vec{x}_s)
=\frac{p}{p_s}$. There is a unique \textit{degree map} $\mathbf{d}:\;
K_0(Coh(\X_{\p,\l})) \to \Z$ satisfying 
$\mathbf{d}([\O_\X(\vec{x})]) = \partial(\vec{x})$ and
$\mathbf{d}([S])=deg(x)$ if $S$ is a simple torsion sheaf supported at a
point $x$.

Now for any $\mathcal{M} \in Coh(\X_{\p,\l})$ set $\chi(\mathcal{M})=
\mathrm{dim\;Hom}(\O_\X,\mathcal{M})-\mathrm{dim\;Ext}^1(\O_X,\mathcal{M})$
(the \textit{Euler characteristic} of $\mathcal{M}$), and put
$$\overline{\chi}(\mathcal{M}) = \sum_{k=0}^{p-1} \chi(\mathcal{M}
(-k \vec{\omega})).$$

\begin{theo}[Riemann-Roch theorem, \cite{GL}] For any $\mathcal{M}\in
Coh(\X_{\p,\l})$ we have 
$$\overline{\chi}(\mathcal{M})=\mathbf{r}(\mathcal{M}) 
\overline{\chi}(\O_{\X}) + \mathbf{d}(\mathcal{M}),$$
and 
$\overline{\chi}(\O_\X)=-\frac{p}{2}\partial(\vec{\omega}).$
\end{theo}

The previous theorem motivates the following definition of the genus of 
$\X_{\p,\l}$ : 
$$g(\X_{\p,\l})=1+ \frac{1}{2} \partial(\vec{\omega}).$$
 Using
this notation, the Riemann-Roch theorem admits the following corollary.
We set $\langle M, N \rangle = \mathrm{dim\;Hom}\;(M,N)-
\mathrm{dim\;Ext}^1(M,N)$. Then, for all coherent sheaves $M$ and $N$,
\begin{equation}\label{E:RR}
\sum_{k=0}^{p-1} \langle M(k\vec{\omega}), N \rangle = p(1-g_{\X_{\p,\l}})
\mathbf{r}(M)\mathbf{r}(N) + \mathbf{r}(M)\mathbf{d}(N)-\mathbf{r}(N)
\mathbf{d}(M).
\end{equation}

Let us 
associate to the sequence $(p_1,\ldots, p_l)$ the ``star Dynkin diagram''
$\mathbb{T}_{p_1,\ldots, p_n}$ (see Section 1.3). The genus of $\X_{\p,\l}$
and the properties of $\mathbb{T}_{p_1,\ldots,p_n}$ are related in the 
following way :
\begin{enumerate}
\item[i)] $g(\X_{\p,\l}) <1$ if and only if $\mathbb{T}_{p_1,\ldots ,p_n}$
is a finite Dynkin diagram,
\item[ii)] $g(\X_{\p,\l}) =1$ if and only if $\mathbb{T}_{p_1,\ldots ,p_n}$
is an affine Dynkin diagram.
\end{enumerate}
\vspace{.2in}
\paragraph{\textbf{2.7.}} One of our main tools in the study of the category
$Coh(\X_{\p,\l})$ will be the Harder-Narasimhan (HN) filtration, introduced
in \cite{GL}.
 Define the slope of a nonzero coherent sheaf
$\mathcal{F} \in Coh(\X_{\p,\l})$ as $\mu(\mathcal{F})=
\frac{\mathbf{d}(\mathcal{F})}{\mathbf{r}(\mathcal{F})} \in \mathbb{Q} \cup 
\{\infty\}$ (where
$\mathbf{r}(\mathcal{F})$ and $\mathbf{d}(\mathcal{F})$ are the rank and degree
of $\mathcal{F}$). Note that $\mu(\mathcal{F})=\infty$ if and only if 
$\mathcal{F}$ is a torsion sheaf.
The sheaf $\mathcal{F}$ is called stable (resp. semistable)
if for every subsheaf $\mathcal{G} \subset \mathcal{F}$ we have 
$\mu(\mathcal{G}) < \mu(\mathcal{F})$ (resp. $\mu(\mathcal{G}) \leq
 \mu(\mathcal{F})$). For any $q \in \Q$ let $\mathcal{C}_q$ be the full
subcategory of $Coh(\X_{\p,\l})$ consisting the zero sheaf together with
all sheaves of slope $q$. 

\begin{prop}[\cite{GL}]\label{P:HN1}
For each $q \in \Q$ the subcategory $\mathcal{C}_q$
is abelian and closed under extensions. Moreover, each object is finite length
and the simple objects are the stable bundles. Finally,
$\mathrm{Hom}(\mathcal{F},\mathcal{G})=0$ if $\mathcal{F} \in \mathcal{C}_q$,
 $\mathcal{G} \in \mathcal{C}_{q'}$ and $q>q'$.\end{prop}

The relevance of this notion for us is then explained by the following results.

\begin{prop}[\cite{GL}] If $\Gamma$ is finite (resp. affine) then every 
indecomposable locally free sheaf is stable (resp. semistable).\end{prop}

From Serre duality and the fact that $\partial(\vec{\omega}) \leq 0$
if $g(\X_{\p,\l}) \leq 1$ and $\partial(\vec{\omega})=0$ if $g(\X_{\p,\l})=1$
 we deduce

\begin{cor}\label{C:HN1} Suppose that $\Gamma$ is finite. Then 
$$q \leq q' \Rightarrow \mathrm{Ext}^1(\mathcal{F},\mathcal{G})=0\;
\mathrm{for\;all}\; \mathcal{F} \in \mathcal{C}_q, \mathcal{G} \in 
\mathcal{C}_{q'}.$$
Suppose that $\Gamma$ is affine. Then
$$q < q' \Rightarrow \mathrm{Ext}^1(\mathcal{F},\mathcal{G})=0\;
\mathrm{for\;all}\; \mathcal{F} \in \mathcal{C}_q, \mathcal{G} \in 
\mathcal{C}_{q'}.$$
\end{cor}

\vspace{.1in}
\section{Ringel-Hall algebras}
\vspace{.2in}
\paragraph{\textbf{3.1.}} Let $k$ be a finite field with $q$ elements and
let $\mathcal{A}$ be a (small) abelian k-linear
category. We assume the following :
\begin{enumerate}
\item[i)] $\mathcal{A}$ is hereditary, i.e $Ext^2(V,W)=0$ for any objects
$V$ and $W$,
\item[ii)] $\mathcal{A}$ is Hom and Ext-finite, i.e for any objects 
$V,W$ the spaces $\mathrm{Hom}(V,W)$ and
$\mathrm{Ext}^1(V,W)$ are finite-dimensional.
\end{enumerate}
Let $Iso(\mathcal{A})$ denote the set of isomorphism classes of objects in
$\mathcal{A}$ and
let $K_0(\mathcal{A})$ be the Grothendieck group of $\mathcal{A}$.
The assignement
\begin{align*}
Iso(\mathcal{A}) \times Iso(\mathcal{A}) &\to \Z\\
\langle V, W \rangle &\mapsto \mathrm{dim\;Hom}(V,W)-\mathrm{dim\;Ext}^1(V,W)
\end{align*}
descends to a well-defined biadditive form
 $\langle \;,\;\rangle:\; K_0(\mathcal{A}) \otimes_{\Z}
K_0(\mathcal{A}) \to \Z$ called the \textit{Euler form}.
\vspace{.2in}
\paragraph{\textbf{3.2.}} Let $S,T,U \in Iso(\mathcal{A})$. Let 
$\mathcal{P}_{S,T}^U$ be the (finite) set of all exact sequences
$$0 \to T \stackrel{f}{\to} U \stackrel{g}{\to} S \to 0.$$
The group $\mathrm{Aut}(S) \times \mathrm{Aut}{T}$ acts faithfully on
$\mathcal{P}_{S,T}^U$ and we set
\begin{equation}\label{E:31}
\textbf{P}_{S,T}^U= \left| \mathcal{P}_{S,T}^U / \mathrm{Aut}(S) \times
\mathrm{Aut}(T)\right| \in \N.
\end{equation} 

Consider the ring $\C_q:=\C[v,v^{-1}]/(v^{2}-q)$. We endow the free
$\C_q$-module 
$$\mathbf{H}_{\mathcal{A}}:= \bigoplus_{S \in Iso(\mathcal{A})} \C_q[S]$$
with an associative algebra structure by setting
\begin{equation}\label{E:32}
[S] \cdot [T]:=v^{\langle S,T \rangle} \sum_{U} \mathbf{P}_{S,T}^U [U].
\end{equation}
Such algebras were first considered by Hall (\cite{Hall}, see Section 3.3)
when $\mathcal{A}$ is the category of representations of a discrete valuation
ring with finite residue field, and later in much more generality by Ringel
(see e.g \cite{Ri1}).
\vspace{.1in}
\section{Ringel-Hall algebras of cyclic quivers}
\vspace{.2in}
\paragraph{\textbf{4.1.}}
Let us consider the case of the quiver $\mathcal{C}_1$.
The isomorphism classes of indecomposable (nilpotent) representations of
$\mathcal{C}_1$ over the finite field $\mathbb{F}_q$ with $q$ elements
is in
natural bijection with $\N^*$ and we denote by $|n)$ the class of the
indecomposable representation of dimension $n$. The isomomorphism
classes of (nilpotent) representations is thus in bijection with the set
$\Pi$ of all partitions via the assignment $\lambda = (\lambda_1, \ldots
,\lambda_r) \mapsto |\lambda)=[\lambda_1) \oplus \cdots \oplus |\lambda_r)$.
Let $\mathbf{H}_1^{(q)}$ be the Hall algebra of $\mathcal{C}_1$ over
$\mathbb{F}_q$ and set $\mathbf{e}_r=v^{r(r-1)}|(1^r))$. It is known that
the structure constants for $\H_1^{(q)}$ are polynomials in $q$ and one can
consider $\H_1^{(q)}$ as a $\C_v$-algebra.
\vspace{.1in}
\paragraph{}Let $\Gamma=\C[y_1,y_2, \ldots]^{\mathfrak{S}_\infty}$
be Macdonald's ring of symmetric functions (see \cite{M}) and let 
$\{e_r\}_{r \in \N}$ denote the elementary symmetric polynomial. The 
following well-known result is due to P. Hall.
\begin{theo} The assignement ${e}_r \mapsto \mathbf{e}_r$ uniquely extends to
an isomorphism of algebras $\Upsilon:\Gamma\otimes_{\C}\C_v \stackrel{\sim}{\to} \mathbf{H}^{(q)}_1$.
\end{theo}
Under this isomorphism the element $|\lambda)$ corresponds to $q^{-n(\lambda)}
P_\lambda(q^{-1})$ where $n(\lambda)=\sum_i (i-1)\lambda_i$ and $P_\lambda(t)$
is the Hall-Littlewood polynomial (see \cite{M}). Let $p_r \in \Gamma$ denote
the power-sum symmetric function and set 
$$\mathbf{h}_r=\frac{[r]}{r}
\Upsilon(p_r).$$
It follows from \cite{M}, III.7, Ex. 2, that
\begin{equation}\label{E:4hr}
\mathbf{h}_r=\frac{[r]}{r}\sum_{\lambda, |\lambda|=r} n(l(\lambda)-1)|\lambda),
\end{equation}
where $n(l)=\prod_{i=1}^{l}(1-v^{2i})$.
\vspace{.1in}
\paragraph{\textbf{4.2.}} Let $p \in \N$ and let $\mathbf{H}^{(q)}_p$ be the 
Hall algebra of the 
cyclic quiver $\mathcal{C}_p$ over $\mathbb{F}_q$ (we will also denote it 
simply by $\mathbf{H}_p$ when there is no ambiguity in the ground field).
Denote by $\{\epsilon_i\}_{i \in \Z/p\Z}$ 
the canonical basis of $\Z^{\Z/p\Z}$. For each $i \in \Z/p\Z$ and $l \in \N$
define the \textit{cyclic segment} $[i;l)$ to be the image of the projection
to $\Z/p\Z$ of the segment $[i'-(l-1), i']$ for any $i' \in \Z,\; i' \equiv i\;
(\mathrm{mod}\;p)$. A \textit{cyclic multisegment} is a finite linear 
combination $\mathbf{m}=\sum_{i,l} a_{il}[i;l)$ with $a_{il} \in \N$. The
isomorphism classes of representations (resp.
indecomposable representations) of $\mathcal{C}_p$ over $k$ are in natural
bijection with the set of cyclic multisegments (resp. cyclic segments). For
$i \in \Z/p\Z$ we let $S_i=[[i;1)]$ be the class of the corresponding
simple object, More generally, for any $n \geq 1$ we denote by $S_i^{(n)}$ the
class of the semisimple module $[i;1)^{\oplus n}$. Finally, let 
$d(\mathbf{m})=\sum_{i,l} a_{il}(\epsilon_i+
\cdots \epsilon_{i-(l-1)}) \in \N^{\Z/p\Z}$ be the \textit{degree} of a 
multisegment $\mathbf{m}$ and set $\delta=(1,\ldots,1)$. 
The algebra $\mathbf{H}_p$ is graded and we denote by 
$\mathbf{H}_p[d]$ the
piece of degree $d$. It follows from 
\cite{Ri1} that the assignment
$E_i^{(n)} \mapsto v^{n(n-1)/2}[S_i^{(n)}]$ defines an embedding $\mathbf{U}_v^+
(\widehat{\mathfrak{sl}}_p) \hookrightarrow \mathbf{H}_p$. 
\begin{theo}[\cite{S1}] The embedding $\mathbf{U}_v^+
(\widehat{\mathfrak{sl}}_p) \hookrightarrow \mathbf{H}_p$ extends to
an isomorphism
$$\mathbf{H}_p \simeq 
\mathbf{U}^+_v(\widehat{\mathfrak{sl}}_p) \otimes_{\C_v} \mathcal{Z}$$
where $\mathcal{Z}=\C_v
[x_1,x_2,\cdots]$ is a central subalgebra and
where the element $x_i$ is homogeneous of degree $i\delta$.
\end{theo}
\vspace{.1in}
\paragraph{\textbf{Remark.}} The algebra structure $\bullet: \H_p \otimes \H_p
\to \H_p$ used in \cite{S1} differs slightly from the one defined by 
(\ref{E:32}). Namely, if $f \in \H_p[d]$ and $f' \in \H_p[d']$ then
$$f \bullet f'= v^{-2\sum_i d_i d'_i}f' \cdot f.$$
Note also that the opposite orientation of the cyclic quiver is used in 
\cite{S1}.
\vspace{.1in}

\paragraph{}To any $(\lambda_1,
\ldots , \lambda_r) \in \Pi$ we associate
$$f_\lambda =\sum_{j} [0;p\lambda_j) \in
\mathbf{H}_p$$
and we set
\begin{equation}\label{E:341}
\mathbf{H}_p^0:=\bigoplus_{\lambda \in \Pi} \C_v f_{\lambda}
\end{equation}
The assignment $|\lambda) \mapsto f_\lambda$ uniquely extends to an
algebra isomorphism $\Psi:\mathbf{H}_1^{(q)} \stackrel{\sim}{\to} 
\mathbf{H}^0_p$. In particular, $\mathbf{H}^0_p$ is a commutative 
subalgebra of $\H_p$ and is freely generated
by any of the three sets $\{\mathbf{e}_r\}_{r \in \N^*}$, 
$\{\mathbf{l}_r\}_{r \in \N^*}$ and 
$\{\mathbf{h}_r\}_{r \in \N^*}$ where
$$\mathbf{e}_r=\Psi(\mathbf{e}_r), \qquad
\mathbf{l}_r=f_{(r)},
\qquad \mathbf{h}_r=\Psi(\mathbf{h}_r).$$

\begin{lem}\label{L:31}
The multiplication maps $\U_v^+(\widehat{\mathfrak{sl}}_p) \otimes
\mathbf{H}^0_p \to \mathbf{H}_p$ and $\mathbf{H}^0_p \otimes 
\U_v^+(\widehat{\mathfrak{sl}}_p) \to \mathbf{H}_p$ induce isomorphisms of 
$\C_v$-modules.
\end{lem}
\noindent
\textit{Proof.} It is enough to prove both statements for the algebra structure
$\bullet$ used in \cite{S1}.
We first recall some notation and results from \cite{S1}.
By \cite{S1}, Proposition 2.2, there exists a canonical isomorphism $i:\Gamma
\otimes \C_v \stackrel{\sim}{\to} \mathcal{Z}$. Let
$\{s_\lambda\}_{\lambda \in \Pi}$ denote the
Schur basis of $\Gamma$ and set $\mathbf{s}_\lambda = i(s_\lambda)$. Let 
$\mathbf{b}_{\mathbf{m}}$ be the canonical basis of $\mathbf{H}_p$. To a
partition $\lambda \in \Pi$ we associate the multisegments $\mathbf{m}(\lambda)
=\sum_{i} [1-i, \lambda_i-i)$ and $\mathbf{m}^\lambda=\sum_{i,j}[i,\lambda_j)$.
We have the following (see \cite{S1}, Proposition 2.4)
\begin{equation}\label{E:342}
\sum_i E_i \bullet\mathbf{H}_p=\bigoplus_{\mathbf{m} \not\in 
\{ \mathbf{m}^\lambda,
\lambda \in \Pi\}} \C_v \mathbf{b}_{\mathbf{m}}.
\end{equation}
Let
$\Lambda^\infty$ be the level one Fock space of $\mathbf{H}_p$ considered in
\cite{S1} and let $|0\rangle$ be its vacuum vector. Let 
$\{\mathbf{b}^{\pm}_\lambda|\; \lambda \in \Pi\}$ be the Leclerc-Thibon 
canonical bases of $\Lambda^\infty$ (see \cite{LT}). Then (see \cite{S1},
Proposition 3.1 and \cite{LT}, Theorem 6.9)
\begin{equation}\label{E:343}
i)\mathbf{b}_{\mathbf{m}(\lambda)}|0\rangle=\mathbf{b}^+_\lambda, \qquad
ii)x \cdot |0\rangle=0 \Leftrightarrow x \in \bigoplus_{\mathbf{m} \not\in
\{\mathbf{m}(\lambda), \lambda \in \Pi\}} \C_v \mathbf{b}_\mathbf{m}
\end{equation}
\begin{equation}\label{E:344}
\mathbf{s}_\lambda|0\rangle = \mathbf{b}^-_{p\lambda}.
\end{equation}
Finally, it is easy to prove that
\begin{equation}\label{E:345}
\C_v[\mathbf{l}_1,\mathbf{l}_2, \ldots]|0\rangle=\bigoplus_\mu \C_v 
|p\mu\rangle.
\end{equation}
We now prove the lemma by induction. A direct computation shows that
$\mathbf{H}_p[\delta]=\mathbf{U}^+_v(\widehat{\mathfrak{sl}}_p)[\delta]
\oplus \C_v\mathbf{l}_1$. Assume that for any $r <k$ we have
$$\mathbf{H}_p[r\delta] = \big(\mathbf{U}^+_v(\widehat{\mathfrak{sl}}_p)
\otimes \C_v[\mathbf{l}_1,\ldots,\mathbf{l}_r]\big)[r\delta].$$ 
Observe that $\mathbf{b}^\lambda|0\rangle=\mathbf{b}^+_{n\lambda}$ and recall
that for any $\mu$, $\mathbf{b}^{\pm}_\mu \in |\mu\rangle + 
\bigoplus_{\nu <\mu}\C_v |\nu\rangle$.
From (\ref{E:343})i), (\ref{E:344})
 and (\ref{E:345}) it now follows that there exists $P_k \in \C_v[\mathbf{l}_1,
\ldots ,\mathbf{l}_k]$ and $x \in \bigoplus_{\mathbf{m} \neq 
\mathbf{m}^\lambda}\C_v \mathbf{b}_\mathbf{m}$ such that
$$\mathbf{s}_k|0\rangle=(P_k +x)|0\rangle.$$
Hence, from (\ref{E:343})ii) and (\ref{E:342}) we have
$$\mathbf{s}_k \in P_k + x + \bigoplus_{\mathbf{m} \not\in
\{\mathbf{m}(\lambda), \lambda \in \Pi\}} \C_v \mathbf{b}_\mathbf{m}
\subset P_k + \sum_i E_i \bullet\mathbf{H}_p.$$
Using the induction hypothesis, we obtain
$$\mathbf{s}_k \in \big(\mathbf{U}^+_v(\widehat{\mathfrak{sl}}_p)
\otimes \C_v[\mathbf{l}_1,\ldots,\mathbf{l}_r]\big)[k\delta].$$
But $\mathcal{Z}=\C_v[\mathbf{s}_1,\mathbf{s}_2, \ldots]$ and hence by 
Theorem~4.1 
the multiplication map 
$$\big(\mathbf{U}^+_v(\widehat{\mathfrak{sl}}_p)
\otimes \C_v[\mathbf{l}_1,\ldots,\mathbf{l}_r]\big)[k\delta] \to 
\mathbf{H}_p[k\delta]$$
is surjective. By graded dimension argument it is also injective. This
proves the first assertion of the Lemma. The second is proved in a similar
way.\qed
\vspace{.1in}

\paragraph{}The following result makes the link between $\U_v^+
(\widehat{\mathfrak{gl}}_p)$ and $\H_p$ explicit at the classical level.
\begin{lem}\label{L:4.2} There exists a unique isomorphism $\varphi:\;
\U^+(\widehat{\mathfrak{gl}}_p)
\stackrel{\sim}{\to} (\mathbf{H}_p)_{|v=1}$ extending the canonical embedding
$\U^-(\widehat{\mathfrak{sl}}_n) \hookrightarrow (\mathbf{H}_p)_{|v=1}$ such 
that $\varphi(ht^r)=\mathbf{h}_r$ for all $r\geq 1$.
\end{lem}
\noindent
\textit{Proof.} It follows from \cite{Ri2} that under the
canonical embedding $\U^-(\widehat{\mathfrak{sl}}_n) \hookrightarrow 
(\mathbf{H}_p)_{|v=1}$ we have, for $i=1, \ldots, n-1$ and $s \geq 0, r\geq 1$
$$e_it^s=\pm [i;1+ps), \qquad f_it^r=\pm [i+1;pr-1),
\qquad h_it^r=\pm([i;lr)-[i-1;lr)).$$
A direct computation then shows that the assignement 
$ht^r \mapsto \mathbf{h}_r$ extends to a well-defined algebra homomorphism.
The result now follows from Lemma~\ref{L:31}.\qed  
\vspace{.2in}
\paragraph{\textbf{4.3.}} Let $\U_v(\widehat{\mathfrak{sl}}_p)$ be the 
quantized enveloping algebra of the affine Kac-Moody algebra 
$\widehat{\mathfrak{sl}}_p$ and let $\U_v^{(D)}(\widehat{\mathfrak{sl}}_p)$ be 
Drinfeld's new realization (see \cite{Beck}).
In \cite{Beck}, Beck explicitely constructed an isomorphism
of $\C(v)$-algebras $i: \U_v(\widehat{\mathfrak{sl}}_p) \stackrel{\sim}{\to}
\U_v^{(D)}(\widehat{\mathfrak{sl}}_p)$. Let $\pi :
\U_v^{(D)}(\widehat{\mathfrak{sl}}_p) \to \U_v(L\mathfrak{sl}_p)$ be the
canonical projection (obtained by setting the central charge $c=0$). 
Set $\epsilon_0 = \pi \circ i (E_0)$ and let $\U_0^+ \subset
\U_v(L\mathfrak{sl}_p)$ be the subalgebra generated by $\epsilon_0$ and
$x^+_{i,0}$ for $i=1, \ldots p-1$. Let $\U^+ \subset \U_v(L\mathfrak{gl}_p)$
be the $\C_v$-subalgebra generated by $\U^+_0$ and $h_{0,l}$ for $l >0$. The
quantum analogue of Lemma~\ref{L:4.2} is the following conjecture :

\begin{conj} The assignement $\epsilon_0 \mapsto E_0, x^+_{i,0} \mapsto
E_i, h_{0,l} \mapsto \mathbf{h}_l$ extends to an isomorphism of 
$\C(v)$-algebras $\U^+\otimes \C(v) \stackrel{\sim}{\to} \H_{p}\otimes \C(v)$.
\end{conj}
\vspace{.1in}
\paragraph{\textbf{4.4.}} Let $\Gamma$ be a star Dynkin diagram
as in Section 1.4.
We will now define a quantum analog of the subalgebra $\U(\widehat{\n})
\subset \U(\g)$. To each type $A$ subdiagram $J_s$, $s=1,\ldots n$ we
associate the algebra $\mathbf{H}_{p_s}$ where $\phi_s: A_{p_s} 
\stackrel{\sim}{\to} J_s$ and we
denote by 
$E_0^{(s)}, E_{\phi_s(1)}, \ldots, E_{\phi_s(p_s-1)}$ and 
$\mathbf{h}_r^{(s)} \in \mathbf{H}_{p_s}$ the corresponding elements,
as defined in Section 4.2.
\vspace{.1in}
\paragraph{\textbf{Definition}.} 
Let $\U_v(\widehat{n})$ be the $\C_v$-algebra 
generated by $\mathbf{H}_{p_s}$ for $s=1, \ldots n$, $H_{\star,r}$,
 for $r \in \N^*$ and $E_{\star,k}$ for
$k \in \Z$, subject to the following set of relations :
\begin{enumerate}
\item[i)] For all $s$, $H_{\star,r}=\mathbf{h}_r^{(s)}$, and
$[\H_{p_s},\H_{p_{s'}}]=0$ for $s \neq s'$.
\item[ii)] We have
\begin{equation}\label{E:441}
[H_{\star,r},E_{\star,t}]=\frac{[2r]}{r}E_{\star,r+t},
\end{equation}
\begin{equation}\label{E:4415}
E_{\star,t_1+1}E_{\star,t_2}-v^2E_{\star,t_2}E_{\star, t_1+1}
= v^2E_{\star,t_1}E_{\star,t_2+1}-E_{\star,t_2+1}E_{\star,t_1},
\end{equation}
\item[iii)]
\begin{equation}\label{E:442}
a_{\star,j}=0 \Rightarrow [E_{\star,t},E_j]=0, \qquad 
vE_{\star,t}E_0^{(s)}-E_0^{(s)}E_{\star,t}=0,
\end{equation}
\begin{equation}\label{E:443}
p_s >2 \Rightarrow [E_0^{(s)}, vE_{\phi_s(1)}E_{\star,t}-E_{\star,t}
E_{\phi_s(1)}]=0,
\end{equation}
\item[iv)] For $s=1, \ldots, n$ and $r \geq 1$ set $E_{\phi_s(1)}^{(r)}=
\frac{r}{[r]}[E_{\phi_s(1)},H_{\star,r}]$. Then for any $r,r_1,r_2
\in \N$ and $t,t_1,t_2 \in \Z$,
\begin{equation}\label{E:444}
\mathrm{Sym}_{r_1,r_2} 
\big\{E^{(r_1)}_{\phi_s(1)} E^{(r_2)}_{\phi_s(1)}E_{\star,t}
-[2]E^{(r_1)}_{\phi_s(1)}E_{\star,t}E^{(r_2)}_{\phi_s(1)} +E_{\star,t}
E^{(r_1)}_{\phi_s(1)} E^{(r_2)}_{\phi_s(1)}\big\}=0,
\end{equation}
\begin{equation}\label{E:445}
\mathrm{Sym}_{t_1,t_2} \big\{E_{\star,t_1}E_{\star,t_2}E^{(r)}_{\phi_s(1)}
-[2]E_{\star,t_1}E^{(r)}_{\phi_s(1)}E_{\star,t_2}+E^{(r)}_{\phi_s(1)}
E_{\star,t_1}E_{\star,t_2}\big\}=0,
\end{equation}
\item[v)] For any $s=1, \ldots, n$ and $r,r_1,r_2 \in \N$, $t \in \Z$ we have
\begin{equation}\label{E:446}
E_{\star, t+1}E^{(r)}_{\phi_s(1)}-v^{-1}E^{(r)}_{\phi_s(1)}E_{\star,t+1}=
v^{-1}E_{\star,t}E^{(r+1)}_{\phi_s(1)}-E^{(r+1)}_{\phi_s(1)}E_{\star,t},
\end{equation}
\begin{equation}\label{E:447}
E^{(r_1+1)}_{\phi_s(1)}E^{(r_2)}_{\phi_s(1)}-v^{2}E^{(r_2)}_{\phi_s(1)}
E^{(r_1+1)}_{\phi_s(1)}=v^2E^{(r_2+1)}_{\phi_s(1)}E^{(r_1)}_{\phi_s(1)}
-E^{(r_1)}_{\phi_s(1)}E^{(r_2+1)}_{\phi_s(1)}.
\end{equation}
\end{enumerate}

Note that $E^{(r)}_{\phi_s(1)} \in \H_{p_s}$ by (\ref{E:4hr}). Moreover,
if Conjecture 4.3 is true then there is a natural homomorphism 
$\U_v(\widehat{\n}) \to \U_v(L\g)$.
\vspace{.1in}
\section{Main theorems}
\vspace{.2in}
\paragraph{\textbf{5.1.}} Fix a finite field $k$ with $q$ elements and
let $\mathbf{p}=(p_1, \ldots, p_n) \in \N^n$, $\l=(\lambda_1, \ldots,
\lambda_n) \in (\mathbb{P}^1(k))^n$ be as in Section 2.1. Let $\Gamma$
be the corresponding star Dynkin diagram, and let $\g$, $\widehat{\g}$,
$\widehat{\mathfrak{n}}$, etc... be associated to $\Gamma$ as in Section 1.
Observe that, from Proposition~2.4
there is a natural identification of $\Z$-modules
$K_0(Coh(\X_{\p,\l})) \stackrel{\sim}{\to}
\widehat{Q}$ given by the assignement
\begin{align*}
[S_i^{(s)}] &\mapsto \alpha_{\phi_s^{-1}(i)}\;\text{for\;} i=1, \ldots, 
p_s-1,\\
[S_0^{(s)}] &\mapsto \delta - \theta^{(s)},\\
[\O_{\X}(k\c)] &\mapsto \alpha_{\star} + k\delta.
\end{align*}
Define the symmetric bilinear form $(\,,\,): K_0(Coh(\X_{\p,\l}))\otimes
K_0(Coh(\X_{\p,\l})) \to \Z$ by $(M,N)=\langle M, N \rangle + \langle N, 
M \rangle$.
\begin{prop} The bilinear form $(\,,\,)$ on $K_0(Coh(\X_{\p,\l}))$ coincides
with the Cartan form on $\widehat{Q}$.
\end{prop}
 This is a direct consequence of Lemma~6.1, proved in Section 6.2.
\vspace{.2in}
\paragraph{\textbf{5.2.}}
Let $x$ be a closed point of 
$\mathbb{X}_{\mathbf{p},\l}$ and let 
$\mathcal{T}_x$ be the Serre subcategory
of $Coh(\mathbb{X}_{\mathbf{p},\l})$ consisting of (torsion) sheaves
supported at $x$. Note that there is a canonical embedding of Hall algebras
$$\H_{\mathcal{T}_x} \hookrightarrow 
\H_{Coh(\mathbb{X}_{\mathbf{p},\l})}.$$
If $x$ is an ordinary closed point of degree $pd$ then there is a canonical 
isomorphism
$$\Theta_x:\; \H_1^{(q^{d})} \stackrel{\sim}{\to} 
\H_{\mathcal{T}_x}.$$
For every $r \in \N$ we set 
$$\mathbf{h}_{r,x}=\begin{cases} 0 & \mathrm{if}\; r \not\equiv 0\;
(\mathrm{mod}\;{d})\\
\Theta_x(\mathbf{h}_{r/d})& \mathrm{if}\;r \equiv 0\;(\mathrm{mod}\;d)
\end{cases}$$
where $\mathbf{h}_{l} \in \H_1^{(q^{d})}$ is defined in Section 4.1. 

If $x=\sigma_s$ for some $s=1, \ldots, n$ then $x$ is of degree 
$\frac{p}{p_s}$ and there 
is a canonical isomorphism
$$\Theta_{\sigma_s}:\; \H_{p_s}^{(q)} \stackrel{\sim}{\to} 
\H_{\mathcal{T}_{\sigma_s}}.$$
For every
$r \in \N$ and $s=1, \ldots, n$ we set
$$\mathbf{h}_{r, \sigma_s} = \Theta_{\sigma_s}(\mathbf{h}_r)$$
where $\mathbf{h}_r$ is defined in Section 4.2. Finally we set
$$T_r=\bigoplus_{x} 
\mathbf{h}_{r,x} \in \H_{Coh(\mathbb{X}_{\mathbf{p},\l})},$$
where the sum ranges over all closed points of 
$\mathbb{X}_{\mathbf{p},\l}$. This sum is finite as $\mathbf{h}_{r,x}=0$
for all but finitely many $x$.

Recall that $S^{(s)}_i$ for $i=0, \ldots, p_s-1$ denotes the simple torsion 
sheaves supported at $\sigma_s$ (see Section 2.4.).
\vspace{.2in}
\paragraph{\textbf{5.3.}} Let $\mathbf{U}_{\p,\l}$ be the subalgebra of
$\H_{Coh(\mathbb{X}_{\mathbf{p},\lambda})}$ generated by $[S^{(s)}_i]$ for
$s=1, \ldots, n$ and $i=0, \ldots, p_s-1$, by $[\O_{\X}(k\c)]$ for $k \in \Z$
and by $T_r$ for $r \in \N^*$. The following is our first
 main result and will be proved in Section 6.
\begin{theo} The assignement $E_{\phi_s(i)} 
\mapsto [S^{(s)}_i]$, $E_0^{(s)} \mapsto [S_0^{(s)}]$, $H_{\star, r} 
\mapsto T_r$, $E_{\star,k} \mapsto [\O_{\X}(k\c)]$
extends to a homomorphism of algebras $\Phi:\U_v(\widehat{\n}) \otimes_{\C_v}
 \C_q\to \U_{\p,\l}$.
\end{theo}

Note that $\U_v(\widehat{\n})$ and $\U_{\p,\l}$ are both $\widehat{Q}$-graded
and that $\Phi$ is compatible with this grading.
\vspace{.2in}
\paragraph{\textbf{5.4.}} Following Ringel, we introduce a ``generic'' 
analogue of the Hall algebra $\U_{\p,\l}$. For each
finite field $k$ we denote by $\H^{(k)}$
the Hall algebra of the category of coherent sheaves on the weighted
projective line $\X_{\p,\l}$ defined over $k$. We will use the
notations $[S_i^{(s)}]^{(k)}$, $T_r^{(k)}$, etc... for the corresponding
elements. Consider the direct
product of rings $\H=\prod_k \H^{(k)}$. 
Let $\breve{\U}_{\p,\l}$ be the subalgebra of $\H$ generated by the
collection of elements
$$\breve{[S_i^{(s)}]}=\prod_k [S_i^{(s)}]^{(k)},$$
$$\breve{T_r}=\prod_k T_r^{(k)},$$
$$\breve{[\O_{\X}({l\c})]}=\prod_k [\O_{\X}(l\c)]^{(k)}.$$
Both $\H$ and $\breve{\U}_{\p,\l}$ are $\C_v$-modules
since each $\H^{(k)}$ is a $\C_v$-module. Moreover, the element $v \in \H$
does not satisfy any polynomial equations $P(v)=0$ and we can consider it as
a formal parameter. In particular, the algebra $\breve{\U}_{\p,\l} 
\otimes_{\C_v} \C(v)$ is well-defined.

Our second main result is the following.
\begin{theo} Assume that $\Gamma$ is finite or affine. The assignement
 $E_{\phi_s(i)} 
\mapsto \breve{[S^{(s)}_i]}$, $E_0^{(s)} \mapsto \breve{[S_0^{(s)}]}$,
$H_{\star, r} 
\mapsto \breve{T_r}$, $E_{\star,k} \mapsto \breve{[\O_{\X}({k\c})]}$
extends to an isomorphism of algebras 
$$\breve{\Phi}:\U_v(\widehat{\n}) \otimes_{\C_v} \C(v)
\stackrel{\sim}{\to} \breve{\U}_{\p,\l}\otimes_{\C_v} \C(v).$$
\end{theo}

This theorem is proved in Section 7 (for the finite type case) and in Section
8 (for the affine case).
\vspace{.2in}
\paragraph{\textbf{Remark.}} We believe that the map $\Phi$ in Theorem ~5.1
is in fact an isomorphism (for all $\Gamma$). When $\Gamma$ is finite 
or affine,
this would follow from the results of Section 7 and 8 together with the
flatness of the deformation $\U_v(\widehat{\n})$ of $\U(\widehat{\n})$.

\vspace{.1in}
\section{Computations of some Hall numbers}
\vspace{.1in}
In this section we check that the map $\Phi$ is well-defined, i.e that
 relations (\ref{E:441})-(\ref{E:447}) are satisfied in the Hall algebra.
This proves Theorem~5.1. For simplicity we write $\X$ for $\X_{\p,\l}$ and
$\O$ for the structure sheaf $\O_{\X}$.
\vspace{.2in}
\paragraph{\textbf{6.1.}} It is clear that if $S$ and $S'$ are torsion sheaves
with disjoint support then $[S][S']=[S \oplus S']=[S'][S]$ in 
$\H_{Coh(\X)}$. Moreover it follows from Sections 4.1 and 4.2 that
the subalgebra generated by $T_r$ for $r \geq 1$ is commutative. Hence 
relation i) in Section 4.4 is satisfied in $\U_{\p,\l}$.
\vspace{.2in}
\paragraph{\textbf{6.2.}} Recall that $\langle \,,\,\rangle$ denotes the 
Euler form on $K_0(Coh(\X))$. From (\ref{E:23}) and Theorem 2.5 one
easily deduces the following lemma. Set 
$$\delta=[\O(\c)]-[\O]=\sum_{i=0}^{p_1-1}[S^{(1)}_i]= \cdots =
\sum_{i=0}^{p_s-1}[S^{(s)}_i] \in K_0(Coh(\X)).$$
This corresponds to the imaginary root $\delta \in \widehat{Q}$ under
the identification 
$$K_0(Coh(\X)) \simeq \widehat{Q}.$$
\begin{lem}\label{L:665} The Euler form $\langle \,,\,\rangle$ is given by :
$$
\langle [\O], [\O] \rangle =1, \qquad \langle [\O], \delta \rangle =1,
\qquad \langle \delta, [\O] \rangle =-1,$$
$$\langle \delta,\delta \rangle =0, \qquad \langle \delta, [S^{(s)}_i] 
\rangle =0, \qquad \langle [S^{(s)}_i],\delta \rangle =0,
$$
and
$$\langle \O, [S^{(s)}_i] \rangle =\begin{cases} 1 & \mathrm{if}\; i=p_s\\
0 & \mathrm{if}\; i \neq p_s \end{cases}$$
$$\langle [S^{(s)}_i], \O \rangle =\begin{cases} -1 & \mathrm{if}\; i=1\\
0 & \mathrm{if}\; i \neq 1 \end{cases}$$
$$\langle [S^{(s)}_{i}], [S^{(s')}_{i'}] \rangle =\begin{cases} 1 & \mathrm{if}
\; s=s',\;i=i'\\
-1 & \mathrm{if}\; s=s',\;i\equiv i'+1\;(\mathrm{mod}\;p_s)\\
0 & \mathrm{otherwise} \end{cases}.$$
\end{lem}
\vspace{.2in}
\paragraph{\textbf{6.3.}} Let $\mathcal{C}_0$ be the Serre subcategory of
$Coh(\X)$ generated by the locally free sheaves $\O(k\c)$ for 
$k \in \N$, by the torsion sheaves supported at ordinary closed points, and by
the sheaves $S^{(s)}_{r\delta}$ for $s=1, \ldots, n$ and $r \geq 1$ defined by
the exact sequences
$$0 \to \O \to \O(r\c) \to S^{(s)}_{r\delta} \to 0$$
induced by multiplication by $x_s^{rp_s}$. Let $Coh(\mathbb{P}^1(k))$ 
be the category of coherent sheaves on $\mathbb{P}^1(k)$.
\begin{lem} The categories $\mathcal{C}_0$ and $Coh(\mathbb{P}^1(k))$ are 
equivalent.\end{lem}
\noindent
\textit{Proof.} By Serre's theorem and by Proposition 2.3, there is a
chain of equivalences
\begin{equation*}
\begin{split}
 Coh(\mathbb{P}^1(k)) &\simeq k[T,U]-modgr/k[T,U]-modgr_0 \\
&\simeq
S(\p,\l)^0-modgr/S(\p,\l)^0-modgr_0 \\
&\simeq \mathcal{C}_0
\end{split}
\end{equation*}
where $S(\p,\l)^0$ denotes the ($\Z$-graded) subring of $S(\p,\l)$ consisting
of elements of degree $k\c$ for some $k\geq 0$.\qed

The structure of the Hall algebra of $Coh(\mathbb{P}^1(k))$ has been
unraveled by Kapranov
\cite{K} and recovered by Baumann and Kassel in \cite{BK}. 
In particular,  relation (\ref{E:4415}) now 
follows from Lemma 6.2 and \cite{BK}, Lemma 16 i).

For every ordinary
closed point $x$ of degree $pd$ we set $$\Xi_x(s)=1 + 
\sum_{\beta \in Iso \mathcal{T}_x}[\beta]s^{\mathbf{d}(\beta)/p}
=1 + \sum_{r \geq 1}\Theta_x\circ \Upsilon(\xi_r)s^{rd},$$
where $\xi_r \in \Gamma$ denotes the complete symmetric function (we avoid the
usual notation $h_r$ by fear of confusion with $\mathbf{h}_r$). For every
exceptional point $\sigma_i$ we set
$$\Xi_{\sigma_i}(s)=1 + \sum_{r \geq 1} \Theta_{\sigma_i}\circ \Psi \circ 
\Upsilon(\xi_r) s^{r}.$$
Now define $\widehat{\xi}_r \in \H_{Coh(\X)}$ by the relation
$$1 + \sum_{r \geq 1} \widehat{\xi}_r s^r=\prod_{x \in \X_{\p,\l}} 
\Xi_{x}(s).$$
\begin{lem} We have
$$1 + \sum_{r \geq 1}\widehat{\xi}_r s^r=exp\big(\sum_{r \geq 1} 
\frac{T_r}{[r]} s^r\big).$$
\end{lem}
\noindent
\textit{Proof.} Classes of torsion sheaves supported at distinct
points of $\X$ commute in the Hall algebra $\H_{Coh(\X)}$. 
Hence,
\begin{equation*}
exp\big(\sum_{r \geq 1} 
\frac{T_r}{[r]} s^r\big)=\prod_{x \in \X} exp\big(\sum_{r \geq 1} 
\frac{\mathbf{h}_{r,x}}{[r]}\big).
\end{equation*}
Thus it is enough to show that for any closed point $x$ we have
$$\Xi_x(s)=exp\big(\sum_{r \geq 1} 
\frac{\mathbf{h}_{r,x}}{[r]}\big).$$
This is a consequence of the definitions of $\Xi_r$ and $\mathbf{h}_{r,x}$
and of the following identity in the ring $\Gamma$
 (see \cite{M}, I.2) :
$$\sum_{r \geq 0} \xi_r s^r =exp\big(\sum_{r \geq 1} \frac{p_r}{r}\big).$$
\qed

From Lemma 6.2 and
\cite{BK}, Lemma 19 we have, for any $r \geq 1$ and any $n \in \Z$,
\begin{equation}\label{E:61}
\widehat{\xi}_r [\O(n\c)]=[r+1][\O((n+r)\c)] + \sum_{s=0}^{r-1}[\O((n+s)\c)]
\widehat{\xi}_{r-s}.
\end{equation}
Relation (\ref{E:441}) is now a consequence of Lemma 6.3 and (\ref{E:61})
(see \cite{BK}, Proposition~25).
\vspace{.2in}
\paragraph{\textbf{6.4.}} We now check relations (\ref{E:442}) and
(\ref{E:443}). Note that the group $L(\mathbf{p})$ acts on 
$\H_{Coh(\X)}$. Thus from (\ref{E:299}) we can assume $t=0$. If
$a_{\star,j}=0$ and $j=\phi_s(k)$ we have, using Serre duality,
$$\mathrm{Ext}^1(S_k^{(s)},\O)=\mathrm{Ext}^1(\O,S_k^{(s)})=0$$
and
$$\mathrm{Hom}(\O, S_k^{(s)})=\mathrm{Ext}^1(S_k^{(s)},\O(\vec{\omega}))=
\mathrm{Ext}^1(S_{k+1}^{(s)},\O)=0.$$
Moreover, $\langle [\O],[S_k^{(s)}]\rangle=\langle [S_k^{(s)}], [\O]\rangle
=0$. The first equality in (\ref{E:442}) follows. Similarly, we have
$$\mathrm{Ext}^1(S^{(s)}_{p_s},\O)=\mathrm{Ext}^1(\O,S_{p_s}^{(s)})=0, \qquad
\mathrm{Hom}(\O,S_{p_s}^{(s)})=k.$$
Thus,
\begin{align*}
[S_{p_s}^{(s)}] \cdot [\O]&=v^2[\O \oplus S_{p_s}^{(s)}]\\
[\O]\cdot [S_{p_s}^{(s)}]&=v[\O \oplus S_{p_s}^{(s)}]
\end{align*}
and the second equality of (\ref{E:442}) follows. To prove (\ref{E:443}) we 
compute in an analogous fashion
\begin{align*}
[S_1^{(s)}]\cdot[\O]&=v^{-1} \big( [\O \oplus S_1^{(s)}] + [\O(\x_i)]\big)\\
[\O]\cdot[S_1^{(s)}]&=[\O \oplus S_1^{(s)}].
\end{align*}
Hence,
\begin{equation}\label{E:666}
v[S_1^{(s)}]\cdot[\O] - [\O]\cdot[S_1^{(s)}]=[\O(\x_i)].
\end{equation}
If $p_s >2$ then from (\ref{E:442}) we deduce
$$\big[ [S_{p_s}^{(s)}],[\O(\x_i)]\big]=
\big[[S_{p_s-1}^{(s)}],[\O]\big](\x_i)=0$$
which proves (\ref{E:443}).
\vspace{.2in}
\paragraph{\textbf{6.5.}}
 Let us set $[S^{(s)}_{1,r}]=\frac{r}{[r]}\big[[S_1^{(s)}],T_r\big]$.
An explicit computation in $\H_{p_s}$ shows that
\begin{equation*}
[S^{(s)}_{1,r}]=\sum_{\lambda, |\lambda|=r} n(l(\lambda))
\big[|\lambda)^{(s)} \oplus S_1^{(s)}\big] + \sum_{\lambda, |\lambda|=r}
n(l(\lambda)-1) \sum_{\nu \triangleleft
\lambda}\big[[1;p_s\nu +1) \oplus |\lambda \backslash \nu)^{(s)}\big],
\end{equation*}
where for simplicity we set $|\lambda)^{(s)}=\Theta_{\sigma_s}(f_\lambda)$
and where we write $\nu \triangleleft \lambda$ if $\nu \in \N$ is a part of
$\lambda$. The following lemma will be important for us.

\begin{lem}\label{L:666} For any $r_1, r_2 \in \N$ we have 
\begin{equation}\label{E:app1}
\big[ [S^{(s)}_{1,r_2}], T_{r_1}\big]=
\frac{[r_1]}{r_1}[S^{(s)}_{1,r_1+r_2}].
\end{equation}
\end{lem}
\noindent
\textit{Proof.} This is an explicit computation in $\H_p$, and 
is proved in Appendix~1.
\vspace{.2in}
\paragraph{\textbf{6.6.}} In this section we check relation (\ref{E:444}).
From (\ref{E:666}) and the relation $v [\O(\x_i)] \cdot [S^{(s)}_1] - 
[S^{(s)}_1]\cdot[\O(\x_i)]=0$ we deduce
\begin{equation*}\label{E:667}
\begin{split}
0&=v (v [S_1^{(s)}] \cdot [\O] - [\O] \cdot [S_1^{(s)}]) \cdot [S_1^{(s)}]
-[S_1^{(s)}] \cdot ( v [s_1^{(s)}]\cdot[\O]-[\O]\cdot [S_1^{(s)}])\\
&=[S_1^{(s)}] ^2 \cdot [\O]-[2][S_1^{(s)}]\cdot [\O] \cdot [S^{(s)}_1]
+ [\O]\cdot [S_1^{(s)}]^2
\end{split}
\end{equation*}
proving (\ref{E:444}) when $r_1=r_2=t=0$. By twisting by $t\vec{c}$ also obtain
the case $r_1=r_2=0$ and $t \in \Z$. Now let us apply $ad\; T_{r_1}$
to $(\ref{E:667})$. We get, using (\ref{E:441})
\begin{equation*}
\begin{split}
0=&\frac{[r_1]}{r_1} \mathrm{Sym}_{r_1, 0}\big\{ [S_{1,r_1}^{(s)}]
\cdot [S^{(s)}_{1,0}] \cdot [\O] -[2] [S^{(s)}_{1,r_1}]\cdot [\O] \cdot
[S^{(s)}_{1,0}] + [\O] \cdot [S_{1,r_1}^{(s)}] \cdot [S^{(s)}_{1,0}]\big\}\\
&\qquad-\frac{[2r_1]}{r_1} \big\{ [S_1^{(s)}]^2 \cdot [\O(r_1\vec{c})]-
[2] [S_{1}^{(s)}] \cdot [\O(r_1\vec{c})] \cdot [S_1^{(s)}]
+ [\O(r_1\vec{c})] \cdot [S_1^{(s)}]^2\big\},
\end{split}
\end{equation*}
Hence,
\begin{equation}\label{E:668}
\mathrm{Sym}_{r_1, 0}\big\{ [S_{1,r_1}^{(s)}]
\cdot [S^{(s)}_{1,0}] \cdot [\O] -[2] [S^{(s)}_{1,r_1}]\cdot [\O] \cdot
[S^{(s)}_{1,0}] + [\O] \cdot [S_{1,r_1}^{(s)}] \cdot [S^{(s)}_{1,0}]\big\}
=0,
\end{equation}
proving (\ref{E:444}) when $r_2=0$. Finally, we apply $ad\;T_{r_2}$ to
(\ref{E:668}). Using Lemma~\ref{L:666}, we get
\begin{equation*}
\begin{split}
0=&-\frac{[r_2]}{r_2} \mathrm{Sym}_{r_1+r_2, 0}
 \big\{[S_{1,r_1+r_2}^{(s)}]
\cdot [S^{(s)}_{1,0}] \cdot [\O] -[2] [S^{(s)}_{1,r_1+r_2}]\cdot [\O] \cdot
[S^{(s)}_{1,0}] \\
&\qquad \qquad \qquad \qquad \qquad \qquad \qquad \qquad \qquad
+ [\O] \cdot [S_{1,r_1+r_2}^{(s)}] \cdot 
[S^{(s)}_{1,0}]\big\}\\
&-\frac{[r_2]}{r_2} \mathrm{Sym}_{r_1, r_2}
\big\{ [S_{1,r_1}^{(s)}]
\cdot [S^{(s)}_{1,r_2}] \cdot [\O] -[2] [S^{(s)}_{1,r_1}]\cdot [\O] \cdot
[S^{(s)}_{1,r_2}] \\
&\qquad\qquad\qquad\qquad\qquad\qquad\qquad\qquad\qquad
+ [\O] \cdot [S_{1,r_1}^{(s)}] \cdot 
[S^{(s)}_{1,r_2}]\big\}\\
&+\frac{[2r_2]}{r_2} \mathrm{Sym}_{r_1, 0}
\big\{ [S_{1,r_1}^{(s)}]
\cdot [S^{(s)}_{1,0}] \cdot [\O(r_2\vec{c})] -[2] [S^{(s)}_{1,r_1}]\cdot 
[\O(r_2\vec{c})] \cdot
[S^{(s)}_{1,0}] \\
&\qquad \qquad \qquad \qquad \qquad \qquad \qquad \qquad \qquad
+ [\O(r_2\vec{c})] \cdot [S_{1,r_1}^{(s)}] \cdot 
[S^{(s)}_{1,0}]\big\}.
\end{split}
\end{equation*}
Using (\ref{E:668}) (with $r_1+r_2$) we deduce that (\ref{E:444}) holds for
all $r_1, r_2$ and $t$.
\vspace{.2in}
\paragraph{\textbf{6.7.}} We now prove relation (\ref{E:445}). Reasoning as in
Section 6.6 and using Lemma 6.5 we see that it is enough to deal with the
case $t_2=r=0$. Twisting (\ref{E:666}) by $n \vec{c}$ we obtain
\begin{equation}\label{E:669}
v [S_1^{(s)}] \cdot [\O(t\vec{c})]-[\O(t\vec{c})]\cdot [S_1^{(s)}]=[\O(t\vec{c}
+ \vec{x}_s)].
\end{equation}
Hence relation (\ref{E:445}) is equivalent to the relation
\begin{equation}\label{E:670}
\begin{split}
0=&[\O(t\vec{c})]\cdot [\O(\vec{x}_s)]-
v^{-1}[\O(\vec{x}_s)]\cdot[\O(t\vec{c})]\\
&+[\O]\cdot[\O(t\vec{c}+\vec{x}_s)] - 
v^{-1}[\O(t\vec{c}+\x_s)]\cdot [\O].
\end{split}
\end{equation}
An explicit computation using Lemma 6.1 gives
\begin{align}\label{E:671}
 [\O(\x_s)]\cdot [\O(t\c)]&=v^t [\O(t\c) \oplus \O(\x_s)],\\
[\O]\cdot[\O(t\c+\x_s)]&=v^{1+t}[\O \oplus \O(t\c+\x_s)].
\end{align}

The computation of the other two products is given by the following lemma. If
$i,j \geq 1$ let
us denote by $P_{i,j}$ the set of pairs of homogeneous polynomials 
$(P, Q) \in k[T,U]$ of degrees $i$ and $j$ respectively such that
$P$ and $Q$ are relatively prime and such that $T$ does not divide $P$.

\begin{lem} We have
\begin{equation*}
\begin{split}
[\O(t\c)]\cdot[\O(\x_s)]=&v^{1-t}\big\{ v^{2t} [\O(\x_s) \oplus \O(t\c)]\\
&
+ \sum_{k=1}^{t-1} \frac{1}{q-1}|P_{k,t-k-1}|[\O(k\c + \x_s) \bigoplus 
\O((t-k)\c)]\big\},
\end{split}
\end{equation*}
\begin{equation*}
\begin{split}
[\O(t\c+\x_s)]\cdot[\O]=&v^{-t}\big\{ v^{2t+2} [\O \oplus \O(t\c+\x_s)]
+(v^2-1)v^{2t} [\O(\x_s) \oplus \O(t\c)]\\
&\qquad\qquad\qquad+ 
\sum_{k=1}^{t-1} \frac{1}{q-1}|P_{t-k,k}|[\O(k\c + \x_s) \bigoplus 
\O((t-k)\c)]\big\}.
\end{split}
\end{equation*}
\end{lem}
\noindent
\textit{Proof.} A reasonning similar to that of \cite{BK}, Proposition 5, 
shows that
the coefficient $c_{\vec{l}_1, \vec{l}_2}$
of $[\O(\vec{l}_1 \oplus \O(\vec{l}_2)]$ in
the product $[\O(t\c)] \cdot [\O(\x_s)]$ is equal to 
$\frac{1}{q-1}|Q_{\vec{l}_1,\vec{l}_2}|$ where $Q_{\vec{l}_1,\vec{l}_2}$
is the set of pairs of relatively prime elements $P,Q \in S(\p,\l)$ of 
respective degrees $\vec{l}_1-\x_s$ and $\vec{l_2}-\x_s$. Comparing
images in the Grothendieck group, it is
clear that $c_{\vec{l}_1, \vec{l_2}}=0$ if $\vec{l}_1 + \vec{l}_2 
\neq t\c + \x_s$. Moreover, observe that if $\vec{l}=k\c + t_1\x_1 + \cdots
+ t_n \x_n$ with $0 \leq t_i <p_i$ then
$$S(\p,\l)[\vec{l}]=x_1^{t_1} \cdots x_n^{t_n}S(\p,\l)[k\c].$$
From there one easily deduces that $c_{\vec{l}_1,\vec{l}_2}=0$ unless
$\{\vec{l}_1,\vec{l}_2\}=\{k\c + \x_s, (n-k)\c\}$ for some $0 \leq k \leq t-1$.
Finally, $S(\p,\l)[k\c] =k[x_s^{p_s},x_r^{p_r}]$ for any $r=1, \ldots, n$,
$r \neq s$. This provides a bijection between $P_{k, t-k-1}$ and 
$Q_{k\c, (t-k)\c-\x_s}$ and proves the first relation in the Lemma. The second
relation is proved in a similar manner. \qed

From (\ref{E:671}) and Lemma 6.5 we see that (\ref{E:445}) is 
equivalent to the relation
$$\sum_{k=1}^{t-1} (v|P_{k,t-k-1}|-v^{-1}|P_{t-k,k}|)[\O(k\c+\x_s) \oplus
\O((t-k)\c)]=0,$$
which in turn is a consequence of the following lemma :

\begin{lem} If $1 \leq k \leq t-1$ then $P_{t-k,k}=v^2P_{k,n-k-1}$.
\end{lem}
\noindent
\textit{Proof.} Let us denote for simplicity by $S_d=(k[T,U])[t]$ the space
of homogeneous polynomials of degree $d$, and let us write $(P,Q)=1$ if
$P$ and $Q$ are relatively prime. Then
\begin{equation*}
\begin{split}
|P_{d_1,d_2}|=&|\{(P,Q) \in S_{d_1} \times S_{d_2}\;|(P,Q)=1, (P,T)=1\}|\\
=&|\{(P,Q) \in S_{d_1}\times S_{d_2}\;|(P,Q)=1\}|\\
&\qquad-|\{(P,Q) \in S_{d_1} \times
S_{d_2}\;|(P,Q)=1, P=TP'\}|.
\end{split}
\end{equation*} 
In a similar way,
\begin{equation*}
\begin{split}
|\{(P,Q) &\in S_{d_1} \times S_{d_2}\;|(P,Q)=1, P=TP'\}|\\
=&|\{(P',Q) \in S_{d_1-1} \times S_{d_2}\;|(P',Q)=1, (T,Q)=1\}|\\
=&P_{d_2,d_1-1}.
\end{split}
\end{equation*}
From \cite{BK}, Lemma 9, we have
$$|\{(P,Q) \in S_{d_1}\times S_{d_2}\;|(P,Q)=1\}|=(q^2-1)(q-1)q^{d_1+d_2-1}$$
if $d_1, d_2 \geq 1$.
Hence if $d_1, d_2 \geq 1$ then
\begin{equation}\label{E:672}
|P_{d_1,d_2}|=(q^2-1)(q-1)q^{d_1+d_2-1}-|P_{d_2,d_1-1}|.
\end{equation}
The Lemma now easily follows from (\ref{E:672}) by induction on $min(d_1,d_2)$.
\qed
\vspace{.2in}
\paragraph{\textbf{6.8.}} We here prove relation (\ref{E:446}).
Using the action of $L(\p)$ and Lemma~\ref{L:666} again we see that it is 
enough to consider the case $r=t=0$. Using the notations of Lemma~\ref{L:666}
we have 
$$[S^{(s)}_{1,1}]= [[1;p_s+1)] + (1-v^2)[|1)^{(s)} \oplus S^{(s)}_1].$$
Relation (\ref{E:446}) can now be checked through the following explicit 
computations :
$$
[\O] \cdot [[1;p_s+1)]= v[\O \oplus [1;p_s+1)]
$$
$$
[\O] \cdot [|1)^{(s)} \oplus S^{(s)}_1] =v [\O \oplus 
|1)^{(s)} \oplus S^{(s)}_1]
$$
$$
[\O(\c)] \cdot [S^{(s)}_1]= [\O(\c) \oplus S_1^{(s)}]
$$
$$
[S_1^{(s)}]\cdot [\O(\c)]=v^{-1}[\O(\c) \oplus S_1^{(s)}] + 
v^{-1}[\O(\c + \x_s)]
$$
$$
[[1;p_s+1)]\cdot [\O]=[\O \oplus [1;p_s+1)] + (1-v^{-2})[\O(\x_s) + |1)^{(s)}]
+v^{-2}[\O(\c+\x_s)]
$$
$$
[|1)^{(s)}\oplus S_1^{(s)}] \cdot [\O]=[\O \oplus |1)^{(s)} \oplus S_1^{(s)}]
+ v^{-2}[\O(\x_s) \oplus |1)^{(s)}] + v^{-2}[\O(\c) \oplus S_1^{(s)}].
$$
\qed
\vspace{.2in}
\paragraph{\textbf{6.9.}} The last relation (\ref{E:447}) is also proved
via a direct computation in $\H_{p_s}$ similar to Lemma~\ref{L:666}
(again, it is enough to assume $r_2=0$), (see Appendix 1).
\vspace{.1in}
\section{The finite type case}
\vspace{.1in}
In this section we assume that $\Gamma$ is finite and prove Theorem~5.2, i.e
we prove that the map $\breve{\Phi}$ in Theorem~5.2 is injective.
\vspace{.2in}
\paragraph{\textbf{7.1.}} Observe that the algebra $\U(\widehat{\n})$ does
not have finite-dimensional weight spaces. We introduce an exhaustive 
filtration as follows. For each $m \in \Z$, let $\widehat{\n}_{\geq m}$
be the Lie subalgebra of $\widehat{\n}$ generated by $\widehat{\n}_s$ for
$s=1, \ldots, n$, $h_{\star, r}$ for $r \geq 1$ and $e_{\star, k}$ for
$k \geq m$. We put $\U^{\geq m}(\widehat{\n})=\U(\widehat{\n}_{\geq m})$.
It is clear that $\U^{\geq m+1}(\widehat{\n}) \subset 
\U^{\geq m}(\widehat{\n})$ and that 
$\U(\widehat{\n})=\bigcup_{m} \U^{\geq m}(\widehat{\n})$.
Note that the assignement $e_{\star, t} \mapsto e_{\star, t+1}$
and the identity on all other generators induces an automorphism $\sigma$ of 
$\U(\widehat{\n})$ which restricts to isomorphisms 
$\sigma: \U^{\geq m}(\widehat{\n}) \stackrel{\sim}{\to} 
\U^{\geq m+1}(\widehat{\n})$. In particular, for any
$\alpha' \in \Z\delta \oplus 
\bigoplus_{j \neq \star} \Z \alpha_j$, we have
\begin{equation}\label{E:dimeq1}
\mathrm{dim}\;\U^{\geq m}(\widehat{\n})[\alpha' + d\alpha_\star]
=\mathrm{dim}\;\U^{\geq m+1}(\widehat{\n})[\alpha' + d\alpha_\star+d\delta].
\end{equation}

\vspace{.1in}

\begin{lem} For any $m \in \Z$ the algebra $\U^{\geq m}(\widehat{\n})$
has finite-dimensional weight spaces.\end{lem}
\noindent
\textit{Proof.} By (\ref{E:dimeq1}) it is enough to assume $m=0$.
Recall the notations of Section 1.4. We have $\widehat{\n}_{\geq 0}=
\widehat{\mathfrak{t}} \oplus \widehat{\mathfrak{f}}_{\geq 0}$
where $\widehat{\mathfrak{f}}_{\geq 0}=\mathfrak{l}[t]$. Hence, by the
Poincarr\'e-Birkhoff-Witt theorem, we have
$$\U^{\geq 0}(\widehat{\n})=\U(\widehat{\mathfrak{t}}) \otimes
\U(\widehat{\mathfrak{f}}_{\geq 0}).$$
Both $\U(\widehat{\mathfrak{t}})$ and $\U(\widehat{\mathfrak{f}}_{\geq 0})$
have finite-dimensional weight spaces. Now, if $\alpha$ is a root of
$\widehat{\n}_{\geq 0}$ then $\alpha=\alpha' + d\delta$
with $\alpha' \in \bigoplus_{i} \Z \alpha_i$ and $d \geq 0$. Note that
there are only finitely many weights $\beta$ of $\widehat{\mathfrak{t}}$
with $\beta=\beta' + d'\delta$ and $d' \leq d$. Hence,
$$\mathrm{dim}\;\U^{\geq 0}(\widehat{\n})[\alpha] = 
\sum_{\substack{\beta=\beta' + d'\delta \\ d' \leq d}} \mathrm{dim}\;
\U(\widehat{\mathfrak{t}})[\beta] \mathrm{dim}\;
\U(\widehat{\mathfrak{f}}_{\geq 0})[\alpha-\beta] < \infty.$$
\qed 

\vspace{.1in}

In a similar manner, we consider the filtration of 
$\U_v(\widehat{\n})$ (resp. of $\U^{(k)}_{\p,\l}$) by subalgebras
$\U^{\geq m}_v(\widehat{\n})$ (resp. $\U_{\p,\l}^{\geq m,(k)}$) 
generated by
$\H_{p_i}$ for $i=1, \ldots, n$ and by $E_{\star, t}$ for $t \geq m$
(resp. generated by the classes of simple torsion sheaves $[S_i^{(s)}]$,
by $T_r$ for $r \geq 1$ and by $[\O(t\c)]$ for $t \geq m$).

\vspace{.1in}

Let us view $\C$ as $\C_v$-module via the evaluation $v \mapsto 1$. It follows
from Proposition~1.1, Theorem~4.2 and Lemma~4.2 that
$$\U_v(\widehat{\n}) \otimes_{\C_v} \C \simeq \U(\widehat{\n}),$$
and that for any $m \in \Z$,
$$\U_v^{\geq m}(\widehat{\n}) 
\otimes_{\C_v} \C \simeq \U^{\geq m}(\widehat{\n}).$$
Hence, for any $\alpha \in \widehat{Q}$ we have
\begin{equation}\label{E:dim71}
\mathrm{dim}_{\C(v)} \U^{\geq m}(\widehat{\n})[\alpha]\otimes_{\C_v} \C(v) 
\leq 
\mathrm{dim}_\C \U^{\geq m}(\widehat{\n})[\alpha].
\end{equation}
For any $k$ the map $\Phi^{(k)}: \U_v(\widehat{\n}) 
\twoheadrightarrow
\U_{\p,\l}^{(k)}$ restricts to a map $\Phi^{(k)} : \U^{\geq m}_v(\widehat{\n})
\twoheadrightarrow \U^{\geq m,(k)}_{\p,\l}$ and, varying $k$, 
we thus obtain a morphism
\begin{equation}\label{E:dim72}
\breve{\Phi}: \U^{\geq m}_v(\widehat{\n})\otimes_{\C_v} \C(v)
\twoheadrightarrow \breve{\U}^{\geq m}_{\p,\l}\otimes_{\C_v}
\C(v).
\end{equation}
Combining (\ref{E:dim71}) with (\ref{E:dim72}), we see that the injectivity
of $\breve{\Phi}$ will follow from the system of inequalities (for all
$\alpha \in \widehat{Q}$ and all $m \in \Z$)
$$\mathrm{dim}\;\breve{\U}^{\geq m}_{\p,\l}[\alpha] \otimes_{\C_v} \C(v)
\geq \mathrm{dim}_\C\;\U^{\geq m}(\widehat{\n})[\alpha],$$
which in turn will follow from the system of inequalities
\begin{equation}\label{E:dim73}
\mathrm{dim}_\C\;\U^{\geq m,(k)}_{\p,\l}[\alpha] 
\geq \mathrm{dim}_\C\;\U^{\geq m}
(\widehat{\n})[\alpha]
\end{equation}
for all finite fields $k$. Finally, observe that the twist by $\c$ induces
isomorphisms $\U^{\geq m,(k)}_{\p,\l} \stackrel{\sim}{\to}
\U^{\geq m+1,(k)}_{\p,\l}$. In particular, for each
$\alpha' \in \Z\delta \oplus 
\bigoplus_{j \neq \star} \Z \alpha_j$, we have
\begin{equation}\label{E:dimeq2}
\mathrm{dim}\;\U^{\geq m,(k)}_{\p,\l}[\alpha' + d\alpha_\star]
=\mathrm{dim}\;\U^{\geq m+1,(k)}_{\p,\l}[\alpha' + d\alpha_\star+d\delta].
\end{equation}
Hence, in view of (\ref{E:dimeq1}), the system of inequalities (\ref{E:dim73})
for all $m \in \Z$ follow from the same system for $m=0$. This set of 
inequalities (and therefore Theorem~5.2) are proved in the rest of this 
section.
\vspace{.2in}
\paragraph{\textbf{7.2.}} From now on we fix a finite field $k$ and 
drop the superscript $(k)$ in the notation. Consider the sheaf
$$\mathbf{T}=\bigoplus_{0 \leq \vec{x} \leq \c} \O(\x).$$
It is proved in \cite{GL} that $\mathbf{T}$ is a tilting sheaf (in the sense
of, e.g \cite{Happ}), and that $\mathrm{End}(\mathbf{T})$ is a canonical
algebra $\Lambda(\p)$ of type $(\p)$ (see \cite{GL}), which is
itself a tilted algebra
of a tame hereditary algebra $\Sigma$ of (extended) Dynkin type 
$\widehat{\Gamma}$. In particular, there is a chain of equivalences
$$D^b(Coh(\X)) \simeq D^b(mod(\Lambda^{op})) \simeq D^b(mod(\Sigma))$$
of (bounded) derived categories. This allows one to translate the well-known
classification of indecomposable objects in $D^b(mod(\Sigma))$ and
$mod(\Sigma)$ into the classification of indecomposable sheaves in
$Coh(\X)$, as was done in \cite{L}, $\S 5.8$ (see also \cite{LM}, 5.4.1.)

\begin{prop}[\cite{L}]\label{P:777} The map $\mathcal{F} \to [\mathcal{F}] \in
K_0(coh(\X)) =\widehat{Q}$ induces a bijection between the set
of (isomorphism classes) of indecomposable vector bundles and the root system
of $\widehat{\mathfrak{f}}$.
\end{prop}

\vspace{.1in}

We begin with the following lemma.
\begin{lem}\label{L:72}
 The algebra $\U_{\p,\l}$ contains the classes $[\O(\vec{x})]$
of all the line bundles in $Coh(\X_{\p,\l})$. Moreover, the subalgebra
$\U^{\geq 0}_{\p,\l}$ contains the class of all line bundles of the form
$\O(\x)$ with $\x \geq 0$.\end{lem}
\noindent
\textit{Proof.} From (\ref{E:299}), the set of generators of 
$\U_{\p,\l}$ is closed under
twisting by $k\c \in L(\p)$ for $k \in \Z$. Thus the set of $\x \in L(\p)$
such that $[\O(\x)] \in \U_{\p,\l}$ is closed under twisting by
$k\c$. But from (\ref{E:666}) we deduce that this set is also invariant under
twisting by $\x_i$ for $i=1, \ldots, n$. Hence,
$[\O(\x)] \in \U_{\p,\l}$ for all $\x \in L(\p)$ as desired.\qed 

\vspace{.1in}

\begin{prop}\label{E:P72}
The algebra $\U_{\p,\l}$ contains the classes of all 
indecomposable locally free sheaves.\end{prop}
\noindent
\textit{Proof.} By Lemma~\ref{L:72}, $\U_{\p,\l}$ contains all line bundles.
We argue by induction on the rank. Assume that $\U_{\p,\l}$ contains all
indecomposable vector bundles of rank $r$, and let $\mathcal{F}$ be an
indecomposable vector bundle of rank $r+1$. By \cite{GL}, Proposition~2.6,
there exists a line bundle $\mathcal{L}$, a vector bundle $\mathcal{G}$ of
rank $r$ and an exact sequence
$$0 \to \mathcal{L} \to \mathcal{F} \to \mathcal{G} \to 0.$$
Let us decompose $\mathcal{G}=\mathcal{G}_1 \oplus \cdots \oplus \mathcal{G}_s$
into indecomposable vector bundles, ordered such that $\mu(\mathcal{G}_1)
\leq \mu(\mathcal{G}_2) \cdots \leq \mu(\mathcal{G}_s)$. By assumption,
$[\mathcal{G}_i]\in \U_{\p,\l}$ for each $i$. Furthermore, from
Proposition~\ref{P:HN1} and Corollary~\ref{C:HN1}
it follows that $\mathrm{Ext}^1(\mathcal{G}_i,
\mathcal{G}_j)=0$ if $i \leq j$. Thus
$$c_{\mathcal{G}}[\mathcal{G}]=[\mathcal{G}_1] \cdot [\mathcal{G}_2] \cdots
[\mathcal{G}_s] \in \U_{\p,\l}$$
for some nonzero $c_{\mathcal{G}} \in \C_v$, and 
$[\mathcal{G}] \in \U_{\p,\l}$.
Now let us write
$$[\mathcal{L}] \cdot [\mathcal{G}]=c_{\mathcal{F}} [\mathcal{F}] + R, 
\qquad\text{for\;some\;} 0 \neq c_{\mathcal{F}} \in \C_v.$$
Since there is at most one indecomposable vector bundle with a given class
in $K_0(Coh(\X))$ it follows that $R$ is a sum of terms of the form
$c_{\mathcal{H}}[\mathcal{H}]$ where $H$ is a direct sum of indecomposable
sheaves of rank at most $r$. By the same reasoning as for $\mathcal{G}$
we conclude that $[\mathcal{H}] \in \U_{\p,\l}$. Thus $[\mathcal{F}] \in
\U_{\p,\l}$ and the induction step is complete. \qed
\vspace{.2in}
\paragraph{\textbf{7.3.}} In order to prove (\ref{E:dim73}), we need a 
refinement of Proposition~\ref{P:777} and Proposition~\ref{E:P72}.
The following Lemma is proved in exactly
the same manner as Lemma~\ref{L:72}.
\begin{lem}\label{L:73}
 The subalgebra
$\U^{\geq 0}_{\p,\l}$ contains the classes of all line bundles of the form
$\O(\x)$ with $\x \geq 0$.\end{lem}\qed

\vspace{.1in}

Now let us consider the full subcategories
$$X_0=\{\mathcal{F} \in Coh(\X)\;| 
\mathrm{Ext}^1(\mathbf{T},\mathcal{F})=0\},$$
$$X_1=\{\mathcal{F} \in Coh(\X)\;| 
\mathrm{Hom}(\mathbf{T},\mathcal{F})=0\}.$$
Then $X_0 \cap X_1=\{0\}$, $X_0$ and $X_1$ are both stable under extension, 
$X_1$ is stable under subobjects while $X_0$ is stable under quotients. 
Furthermore, from \cite{L}, Cor.~4.7 we see that every indecomposable object
is either in $X_0$ or in $X_1$.

\vspace{.1in}

Recall the notations from Section 1.4. Let $\Delta_{\mathfrak{l}}$ denote the
set of roots of $\mathfrak{l}$, and let 
$\widehat{\Delta}_{\mathfrak{f}}^{\geq 0}=\Delta_{\mathfrak{l}} 
\oplus \N \delta$
be the set of roots of $\widehat{\mathfrak{f}}_{\geq 0}$.

\begin{lem}\label{E:L1} The map $\mathcal{F} \mapsto [\mathcal{F}] 
\in K_0(Coh(\X_{\p,\l}))$ induces a bijection between the set of indecomposable
vector bundles in $X_0$ and $\widehat{\Delta}^{\geq 0}_{\mathfrak{f}}$.
\end{lem}
\noindent
\textit{Proof.} Let $\mathcal{F}$ be an indecomposable vector bundle.
Since $X_0 \cap X_1 = \{0\}$ and since every indecomposable vector bundle
is either in $X_0$ or in $X_1$, it is enough to show that $\mathcal{F} \not\in
X_1$ if and only if $[\mathcal{F}] \in 
\widehat{\Delta}^{\geq 0}_{\mathfrak{f}}$.
Let us write $[\mathcal{F}]=\alpha' + d \delta + r [\O]$. It follows from
Lemma~\ref{L:665} that
$$\langle \O, \mathcal{F} \rangle= r + d, \qquad \langle \O(\c), \mathcal{F}
\rangle=d.$$
In particular, if $d \geq 0$ then 
$\mathrm{dim}\;\mathrm{Hom}(\mathbf{T},\mathcal{F}) \geq \mathrm{dim}\;
\mathrm{Hom}(\O,\mathcal{F}) \geq \langle \O, \mathcal{F}\rangle >0$, 
from which we deduce that $\mathcal{F}
\in X_0$. On the other hand, if $d <0$ then 
$\mathrm{dim}\;\mathrm{Ext}^1(\mathbf{T},\mathcal{F}) \geq \mathrm{dim}\;
\mathrm{Ext}^1(\O(\c),\mathcal{F})\geq - \langle \O(\c), \mathcal{F} \rangle 
>0$, from which we deduce that $\mathcal{F}
\in X_1$. The Lemma is proved.\qed

\begin{prop}\label{E:P73}
The algebra $\U^{\geq 0}_{\p,\l}$ contains the classes of all
indecomposable vector bundles in $X_0$.\end{prop}
\noindent
\textit{Proof.} We argue again by induction on the rank. It is easy to
see that a line bundle $\O(\x)$ is in $X_0$ if and only if $\x \geq 0$.
From Lemma~\ref{L:72} we deduce the statement of the Proposition for line
bundles. Assume that $\U^{\geq 0}_{\p,\l}$ contains all indecomposable
vector bundles of $X_0$ of rank at most $r$, and let $\mathcal{F} \in X_0$
be indecomposable of rank $r+1$. Since $\mathrm{Hom}(\mathbf{T},\mathcal{F})
\neq 0$ there exists a line bundle $\O(\x)$ with $\x \geq 0$, a vector
bundle $\mathcal{G}$ and an exact sequence
$$0 \to \mathcal{O}(\x) \to \mathcal{F} \to \mathcal{G} \to 0.$$
Since $\mathcal{F} \in X_0$, the vector bundle $\mathcal{G}$ decomposes
as a sum $\mathcal{G}=\mathcal{G}_1 \oplus \cdots \oplus \mathcal{G}_s$
of indecomposable vector bundles $\mathcal{G}_i \in X_0$ of rank at most $r$.
Arguing as in Proposition~\ref{E:P72}, we see that 
$\mathcal{G} \in \U^{\geq 0}_{\p,\l}$. Now let us write
$$[\mathcal{O}(\x)] \cdot [\mathcal{G}]=c_{\mathcal{F}} [\mathcal{F}] + R, 
\qquad\text{for\;some\;} 0 \neq c_{\mathcal{F}} \in \C_v.$$
Since there is at most one indecomposable vector bundle with a given class
in $K_0(Coh(\X))$ it follows that $R$ is a sum of terms of the form
$c_{\mathcal{H}}[\mathcal{H}]$ where $H$ is a direct sum of indecomposable
sheaves of rank at most $r$. Observe that these sheaves are in $X_0$ since
$X_0$ is stable under extension. By the same reasoning as for $\mathcal{G}$
we conclude that $[\mathcal{H}] \in \U^{\geq 0}_{\p,\l}$. 
Thus $[\mathcal{F}] \in \U^{\geq 0}_{\p,\l}$ and the induction is
complete. \qed
\vspace{.2in}
\paragraph{\textbf{7.4.}} The following result is clear from the discussions
of Section 4.2 and Section~4.3.

\begin{lem}\label{L:789}
The Hall algebra $\U^{\infty}_{\p,\l}$ of the category of torsion
sheaves on $\X$ is isomorphic to the subalgebra of 
$\U_v(\widehat{\n})$
generated by $\H_{p_s}$ for $s=1, \ldots, n$ and $H_{\star,r}$ for $r\geq 1$.
\end{lem}
\vspace{.1in}
In particular, for any $\alpha \in \widehat{Q}$, 
we have the equality of dimensions
\begin{equation}\label{E:dimeq3}
\mathrm{dim}(\U^{\infty}_{\p,\l}[\alpha])=
\mathrm{dim}(\U(\widehat{\mathfrak{t}})[\alpha]).
\end{equation}
\vspace{.2in}
\paragraph{\textbf{7.5.}}
We are now in position to prove (\ref{E:dim73}).
Let $\mathcal{F}_1, \ldots, \mathcal{F}_s$ be any indecomposable vector
bundles in $X_0$, ordered so that $\mu(\mathcal{F}_1) \leq
 \cdots \leq \mu(\mathcal{F}_s)$. Then
$$[\mathcal{F}_1] \cdot [\mathcal{F}_2] \cdots [\mathcal{F}_s] = c
[\mathcal{F}_1 \oplus \cdots \oplus \mathcal{F}_s],\qquad
\text{for\;some\;} 0 \neq c \in \C_v.$$
Hence 
$[\mathcal{F}_1 \oplus \cdots \oplus \mathcal{F}_s] \in \U^{\geq 0}_{\p,\l}.$
In this way we obtain a linearly independent set. In particular,
if $\dot{\U}^{\geq 0}_{\p,\l}$ is the subalgebra of $\U^{\geq 0}_{\p,\l}$
consisting of all classes of vector bundles then it follows from 
Lemma~\ref{E:L1}, Proposition~\ref{E:P73} and the PBW theorem that, for any 
$\alpha \in \widehat{Q}$,  we have
$$\mathrm{dim}(\dot{\U}^{\geq 0}_{\p,\l}[\alpha])=
\mathrm{dim}(\U(\widehat{\mathfrak{f}}_{\geq 0})[\alpha]).$$
On the other hand, since the category of torsion sheaves coincides with
$\mathcal{C}_\infty$, it follows from Corollary~7.1 that the multiplication map
$$\dot{\U}^{\geq 0}_{\p,\l} \otimes \U^{\infty}_{\p,\l} \to 
\U^{\geq 0}_{\p,\l}$$
is injective. Thus, for any $\alpha \in \widehat{Q}$, we have
\begin{equation*}
\begin{split}
\mathrm{dim}\; \U^{\geq 0}_{\p,\l}[\alpha] &\geq 
\mathrm{\dim}(\dot{\U}^{\geq 0}_{\p,\l} \otimes \U^{\infty}_{\p,\l})[\alpha]\\
&=\mathrm{dim}(\U(\widehat{\mathfrak{f}}_{\geq 0}) 
\otimes \U(\widehat{\mathfrak{t}})
[\alpha]\\
&=\mathrm{dim}(\U^{\geq 0}(\widehat{\n})[\alpha],
\end{split}
\end{equation*}
proving (\ref{E:dim73}) as desired.
\vspace{.1in}

\paragraph{\textbf{Remark.}} The set 
$$\mathbf{B}=\{[\mathcal{F}]\;|\; \mathcal{F}\;
\text{is\;a\;vector\;bundle\;or\;a\;torsion\;sheaf\;supported\;at\;
exceptional\;points}\}$$
can be considered as a natural analog of the PBW basis of the quantum 
Kac-Moody algebra $\U_q^+(\mathfrak{g})$ constructed using quivers (when
$\mathfrak{g}$ is of finite type). However, in the present case, $\mathbf{B}$
is not a basis of $\U_v(\widehat{\n})$ but rather only of the (quantum
analog of the) subalgebra generated by $\widehat{\n}^+_s$ for $s=1, \ldots, n$
and $e_{\star,t}$ for $t \in \Z$. 

\section{The affine case}

In this section we prove Theorem~5.2 when $\Gamma$ is an affine 
Dynkin diagram (necessarily of type $D_4^{(1)}$, $E_6^{(1)}$, $E_7^{(1)}$
or $E_8^{(1)}$), which we now assume is the case.
\vspace{.2in}
\paragraph{\textbf{8.1.}} Let $k$ be any field. The classification of
indecomposable objects in $Coh(\X)$ is obtained in \cite{LM}. The main 
tool used there is the theory of mutations (as in e.g., \cite{Bondal}), which
we now briefly recall. Let $\mathcal{C}$ be any hereditary k-linear category
with finite-dimensional morphism and extension spaces. We assume that 
$\mathcal{C}$ has Serre duality, i.e that there exists an equivalence
$\tau: \mathcal{C} \to \mathcal{C}$ and an isomorphism of bifunctors
$$(\mathrm{Ext}^1(X,Y))^* \simeq \mathrm{Hom}(Y,\tau X).$$
Let $\mathbf{A}=\{A_1, \ldots, A_p\}$ be a family of pairwise nonisomorphic
indecomposable objects in $\mathcal{C}$ such that $\mathrm{End}(A_i)=k$,
$\tau(A_i)\simeq A_{i+1\;\mathrm{mod}\;p}$ and such that 
$\mathrm{Hom}(A_i,A_j)=0$ if $i \neq j$. We denote by $\mathcal{A}$ the 
additive full subcategory of $\mathcal{C}$ generated by $\mathbf{A}$.

\vspace{.1in}

We define the functor 
$L_{\mathbf{A}}: \mathcal{C} \to \mathcal{C}$ of \textit{left mutation}
as follows. For any object $F$ of $\mathcal{C}$ we set
$\Sigma_{\mathbf{A}}F=\bigoplus_j \mathrm{Hom}(A_j,F) \otimes_k A_j$ and we
denote by $\sigma_F: \Sigma_{\mathbf{A}}F \to F$ the natural map. Then
$L_{\mathbf{A}}(F)$ is defined by the exact sequence
$$\Sigma_{\mathbf{A}}F \to F \to L_{\mathbf{A}}F \to 0.$$
The \textit{right mutation} $R_{\mathbf{A}}: \mathcal{C} \to \mathcal{C}$
is defined in a similar
fashion. Let $F$ be an object of $\mathcal{C}$ and set 
$\Omega_{\mathbf{A}} F=\bigoplus_j \mathrm{Ext}^1(F,A_j)^* \otimes_k A_j$.
The functors $\mathrm{Hom}(\Omega_{\mathbf{A}} F,-)$ and $\mathrm{Ext}^1(F, -)$ from
$\mathcal{A}$ to $\mathrm{mod}\;k$ are isomorphic. This gives rise to an
exact sequence (the universal extension)
$$\mu_F:\;0 \to \Omega_{\mathbf{A}} F \to R_{\mathbf{A}}F \to F \to 0$$
such that the connecting homomorphism $\mathrm{Hom}(\Omega_{\mathbf{A}} F,A)\to 
\mathrm{Ext}^1(F,A)$ is bijective for all $A \in \mathcal{A}$. Note that
$\mu_F$ depends on a choice of isomorphism $\mathrm{Hom}(\Omega_{\mathbf{A}} F,-)
\simeq \mathrm{Ext}^1(F, -)$, but $R_{\mathbf{A}}F$ doesn't. Moreover, the assignement
$F \mapsto R_{\mathbf{A}}F$ is not functorial in general.

\vspace{.1in}

Finally, define full subcategories $\mathcal{A}^{\Vdash}$ and 
$\mathcal{A}^{\vdash}$ of $\mathcal{C}$ by the following conditions :
\begin{align*}
\mathcal{A}^{\Vdash}&=\{F \in \mathcal{C}\;|\; \mathrm{Hom}(F,A_j)=0\;
\mathrm{for\;all\;}j\},\\
\mathcal{A}^{\vdash}&=\{F \in \mathcal{A}^{\Vdash}\;|\;\sigma_F: 
\Sigma_{\mathbf{A}}F \to F\;\mathrm{is\;a\;monomorphism\;}\}.
\end{align*}

\begin{theo}[\cite{LM}] The functor of left mutations $L_{\mathbf{A}}$ induces
an equivalence of categories 
$L: \mathcal{A}^{\vdash} \to \mathcal{A}^{\Vdash}$, with inverse given by
the right mutation $R_{\mathbf{A}}$. Moreover, $L_{\mathbf{A}}$ and 
$R_{\mathbf{A}}$ both commute with the equivalence $\tau$.\end{theo}

In particular, the restriction of $R_{\mathbf{A}}$ to $\mathcal{A}^{\Vdash}$
is a functor. 

\vspace{.1in}

Using the above methods, with $\mathcal{C}=Coh(\X)$ and $\mathbf{A}=
\{\O(i\vec{\omega})\}$, Lenzing and Meltzer proved the following result.
Recall that $\mathcal{C}_q$ denotes the full subcategory consisting of
the zero sheaf together with semistable sheaves of slope $q$. In particular, 
$\mathcal{C}_{\infty}$ is the subcategory of all torsion sheaves.

\begin{theo}[\cite{LM}]\label{T:81} The following holds
\begin{enumerate}
\item[i)] The categories $\mathcal{C}_q$ for $q \in \mathbb{Q} \cup \{\infty\}$
are all canonically equivalent.
\item[ii)] For $q >0$, left and right mutations with respect to 
$\mathbf{A}=\{\O(i\vec{\omega})\}$ induce inverse equivalences 
$$\mathcal{C}_q \xrightarrow{R_{\mathbf{A}}} 
\mathcal{C}_{\frac{q}{q+1}}\xrightarrow{L_{\mathbf{A}}} \mathcal{C}_q.$$
\end{enumerate}
\end{theo}

In particular, the categorie $\mathcal{C}_q$ are described, like 
$\mathcal{C}_{\infty}$, in terms of categories of representations of cyclic quivers
(see Section 2.4).
\vspace{.2in}

\paragraph{\textbf{8.2.}} Recall that the symmetric bilinear form $(\,,\,)$
on $K_0(Coh(\X))$ coincides with the Cartan form on $\widehat{Q}$
(Proposition~5.1). This, together with Corollary 1.1
and \cite{LM}, Theorem~4.6. implies the following result. We denote by 
$\widehat{\Delta}_{\n}$ the set of roots of $\widehat{\n}$.

\begin{prop}\label{P:81}
There is an indecomposable sheaf $\mathcal{F} \in Coh(\X)$
such that $[\mathcal{F}]=\alpha \in \widehat{Q}$ if and only if $\alpha \in
\widehat{\Delta}_{\n}$. Such a sheaf is unique (up to isomorphism) if and only
if $\alpha$ is a real root. Moreover, $\mathcal{F}$ is a torsion sheaf
(resp. a locally free sheaf) if and only if $[\mathcal{F}]$ is a root of
$\widehat{\mathfrak{t}}$ (resp. a root of $\widehat{\mathfrak{f}}$).
\end{prop}

\vspace{.1in}

Combining Theorem~\ref{T:81} and Proposition~\ref{P:81} we obtain the following
description of the structure of the set $\widehat{\Delta}_{\n}$.
Consider linear forms $\mathbf{d} : \widehat{Q} \to \Z$ and $\mathbf{r} :
\widehat{Q} \to \Z$ corresponding to the degree and rank forms on 
$K_0(Coh(\X))$. For $\alpha \in \widehat{Q}$ set 
$\mu(\alpha)=\frac{\mathbf{d}(\alpha)}{\mathbf{r}(\alpha)} \in \mathbb{Q} 
\cup \{\infty\}$ and for $q \in \mathbb{Q} \cup \{\infty\}$ put
$\widehat{\Delta}_q=\{\alpha \in \widehat{\Delta}_\n\;|\;\mu(\alpha)=q\}$
and $\widehat{\n}_q=
\bigoplus_{\alpha \in \widehat{\Delta}_q} \widehat{\n}[\alpha]$. 

\begin{cor}\label{cor:wei} The following holds :
\begin{enumerate}
\item[i)] Each
$\widehat{\n}_q$ is a Lie subalgebra,
 $\widehat{\n}=\bigoplus_{q}\widehat{\n}_q$ (as graded vector spaces) and 
$\widehat{\n}_\infty=\widehat{\mathfrak{t}}$,
\item[ii)] For each $q$ there is a 
canonical isomorphism $u_q:\widehat{\Delta}_q \simeq \widehat{\Delta}_\infty$ which lifts 
to a vector space isomorphism 
$\widehat{\n}_q \simeq \widehat{\n}_{\infty}$.
\item[iii)] The automorphism $\mathcal{F} \mapsto \mathcal{F}(\vec{\omega})$ of
$Coh(\X)$ induces, for each $q$ an automorphism $\tau$ of 
$\widehat{\Delta}_q$ of order $p$.
\end{enumerate}
\end{cor}
\noindent
\textit{Proof.} Assertions i) and iii) are immediate, as is the first
claim of ii). Observe that $u_q$ is induced by equivalences of categories, and
thus compatible with the Euler form. In particular, $u_q$ maps real roots
to real roots and imaginary roots to imaginary roots. But all real roots of
$\widehat{\n}$ have multiplicity one and all imaginary roots have multiplicity
$|\Gamma|$. Therefore, $u_q$ lifts (in a noncanonical way) to a vector
space isomorphism $\widehat{\n}_q \simeq \widehat{\n}_{\infty}$.  \qed
\vspace{.1in}
\paragraph{\textbf{Remarks.}}It is possible to show that in fact
$\widehat{\n}_q \simeq \widehat{\n}_\infty$ as Lie algebras. We won't need this
result.
\vspace{.1in}
\paragraph{}We will call a root $\alpha \in \widehat{\Delta}_q$ \textit{simple} if it
corresponds, under the equivalence $\widehat{\Delta}_q \simeq 
\widehat{\Delta}_\infty$, to the class of a simple torsion sheaf 
$[S_i^{(s)}]$ (supported at an exceptional point). Note that all simple
roots are real. In addition, the simple roots $\alpha \in \widehat{\Delta}_q$
decompose as a union of $n$ $\tau$-orbits of size $p_1, \ldots, p_n$
respectively, and the set of simple sheaves in $\mathcal{C}_q$ is naturally
in bijection with $\Z/p_1\Z \sqcup \cdots \sqcup \Z/p_n\Z$. On the other hand,
there exists a unique minimal imaginary root $\delta_q \in \widehat{\Delta}_q$
of which all other imaginary roots are multiples. 
\vspace{.2in}

\paragraph{\textbf{8.3.}} We now introduce some exhaustive filtrations on
$\U(\widehat{\n})$, $\U_v(\widehat{\n})$ and $\U^{(k)}_{\p,\l}$ 
respectively. Let
$\widehat{\n}_{\geq 0}$ be the subalgebra of $\widehat{\n}$ generated by
$\widehat{\mathfrak{t}}$ and by the root vectors corresponding to the
simple roots $\alpha_*,\tau(\alpha_*), \ldots, \tau^{p-1}(\alpha^*)$ in 
$\widehat{\Delta}_0$. We set $\U^{\geq 0}(\widehat{\n})=
\U(\widehat{\n}_{\geq 0})$.

\begin{lem} The algebra $\U^{\geq 0}(\widehat{\n})$ has finite-dimensional
weight spaces.\end{lem}
\noindent
\textit{Proof.} For $m \in \Z$, let us denote by 
$\widehat{\mathfrak{u}}_{\geq m}$ the subalgebra of $\widehat{\n}$ generated
by $\widehat{\mathfrak{t}}$ and $e_{\star,t}$ for $t \geq m$. This defines
an exhaustive filtration of $\widehat{\n}$, and the proof of Lemma~7.1 shows
that $\U(\widehat{\mathfrak{u}}_{\geq m})$ has finite-dimensional weight spaces
for any $m$. Now choose $m \ll 0$ such that $\widehat{\mathfrak{u}}_{\geq m}$
contains the root vectors of weight
$\alpha_*, \tau(\alpha_*), \ldots, \tau^{p-1}(\alpha_*)$. Then, for any $\alpha
\in \widehat{Q}$ we have
$$\mathrm{dim}\;\U^{\geq 0}(\widehat{\n})[\alpha] \leq \mathrm{dim}\;
\U(\widehat{\mathfrak{u}}_{\geq m})[\alpha] < \infty.$$
\qed

The following result is proved in the same way as Lemma~7.2.
\begin{lem} The algebra $\U^{(k)}_{\p,\l}$ contains the classes of all line 
bundles $\O(\x)$.\end{lem}

Let $E(\x)\in \U_v(\widehat{\n})$ be the element satisfying $\Phi(E(\x))=
[\O(x)]$ provided by the proof of Lemma~7.2. Let $\U^{\geq 0,(k)}_{\p,\l}$
(resp. $\U^{\geq 0}_v(\widehat{\n})$) be generated by the classes of the 
simple 
torsion sheaves $[S_i^{(s)}]$, by $T_r$ for $r \geq 1$ and by the classes
of the line bundles $[\O(i\vec{\omega})]$ for $i=0, \ldots, p-1$ (resp. by
$\H_{p_i}$ for $i=1, \ldots, n$ and by $E(i\vec{\omega})$ for $i=0, \ldots, 
p-1$). It is clear that $\Phi(\U^{\geq 0}_v(\widehat{\n})) = 
\U^{\geq 0, (k)}_{\p,\l}$.

Let us denote by $\sigma$ the automorphism of $\U(\widehat{\n})$
(resp. of $\U_v(\widehat{\n})$) defined by $\sigma(e_{\star,t})=e_{\star,t+1}$
and the identity on all other generators (resp. by $\sigma(E_{\star,t})=
E_{\star, t+1}$ and the identity on all other generators). Let us also denote 
by the same letter the automorphism of $\U_{\p,\l}^{(k)}$ defined by
$\sigma([\mathcal{F}])=[\mathcal{F}(\c)]$ for any $\mathcal{F}\in 
Coh(\X)$. It is clear that $\sigma$ commutes with the morphism
$\Phi^{(k)}: \U_v(\widehat{\n}) \to \U^{(k)}_{\p,\l}$. Finally, for $m \in \Z$
we put
$$\U^{\geq m}(\widehat{\n})=\sigma^m(\U^{\geq 0}(\widehat{\n})), \quad
\U^{\geq m}_v(\widehat{\n})=\sigma^m(\U^{\geq 0}_v(\widehat{\n})), \quad
\U^{\geq m, (k)}_{\p,\l}=\sigma^m (\U^{\geq 0, (k)}_{\p,\l}).$$
Arguing in the same way as in Section~7.1, we see that Theorem~5.2 will follow
from the set of inequalities
\begin{equation}\label{E:dim81}
\mathrm{dim}_\C\;\U^{\geq 0,(k)}_{\p,\l}[\alpha] 
\geq \mathrm{dim}_\C\;\U^{\geq 0}
(\widehat{\n})[\alpha]
\end{equation}
for all finite fields $k$ and all $\alpha \in \widehat{Q}$. This set of
inequalities is proved in the rest of this Section.
\vspace{.2in}

\paragraph{\textbf{8.4.}} From now on we fix the finite field $k$ and drop the
superscript $(k)$ in all notations. The aim of this paragraph is to prove the
following result.

\begin{prop}
$\U_{\p,\l}^{\geq 0}$ contains the classes of all simple sheaves in 
$\mathcal{C}_q$ for $q>0$.\end{prop}

\vspace{.1in}

We start with some preliminary technical results. Define a sequence $F_m$ for 
$m \geq 1$ (the \textit{Farey sequence}) inductively as follows. 
We set $F_1=\{\frac{0}{1}, \frac{1}{0}\}$ and if $F_m=\{\frac{0}{1}, 
\frac{a_1}{b_1}, \ldots, \frac{a_r}{b_r}, \frac{1}{0}\}$ then
$$F_{m+1}=\big\{\frac{0}{1},\frac{a_1}{b_1+1}, \frac{a_1}{b_1}, \ldots,
\frac{a_i}{b_i}, \frac{a_i+a_{i+1}}{b_i+b_{i+1}},\frac{a_{i+1}}{b_{i+1}},
\ldots, \frac{1}{0}\big\}.$$
In the above, we treat the fraction $\frac{1}{0}$ as representing $\infty$.
The following facts are well-known about the sequences $F_m$ :
\begin{enumerate}
\item[i)] Each new term $\frac{a_i+a_{i+1}}{b_i+b_{i+1}}$ appearing in
$F_{m+1}$ is a reduced fraction.
\item[ii)] If $\frac{a_i}{b_i}$ and $\frac{a_{i+1}}{b_{i+1}}$ are consecutive
entries in some $F_m$ then $b_ia_{i+1}-b_{i+1}a_i=1$.
\item[iii)] Every positive rational number belongs to $F_m$ for $m \gg 0$.
\end{enumerate}

\vspace{.1in}

\begin{prop} Let $\frac{a}{b}$ and $\frac{c}{d}$ be consecutive entries
in some $F_m$. Let $\U=\{U_i\}_{i \in I}$ be simple sheaves in 
$\mathcal{C}_{\frac{a}{b}}$ forming a single $\tau$-orbit of size $p$.
Let $\mathcal{U}$ be the additive closure of $\U$. Then
\begin{enumerate}
\item[i)] For any $i \in I$, $U_i$ is a vector bundle of rank $b$ and of degree
$a$,
\item[ii)] $\mathcal{C}_{\frac{a+c}{b+d}} \subset \mathcal{U}^{\vdash}$
and the right mutation $R_\U$ restricts to an equivalence
$\mathcal{C}_{\frac{c}{d}} \stackrel{\sim}{\to} \mathcal{C}_{\frac{a+c}{b+d}}$.
\end{enumerate}
Let us denote by $\{V_j\}_{j \in J}$ the collection of all simple sheaves
in $\mathcal{C}_{\frac{c}{d}}$. Then for any $j \in J$, 
\begin{enumerate}
\item[iii)]
$\sum_i \mathrm{dim\;Ext}^1(V_j,U_i)=\frac{p}{c_j},\qquad
\mathrm{dim\;Ext}^1(V_j,U_i) \leq 1\quad\mathrm{for\;all}\;i,$
\end{enumerate}
where $c_j \in \{p_1, \ldots, p_n\}$ is the size of the $\tau$-orbit of $V_j$.
\end{prop}
\noindent
\textit{Proof.} We will give a proof by induction. The result is easily checked
for $F_1$ (see \cite{LM}, Proposition 3.9 for ii)). Let us assume that 
properties i) and iii) are proved for all consecutive entries of $F_l$. We will
prove ii) for entries in $F_l$ and i) and iii) for entries in $F_{l+1}$.

Let $\frac{a}{b}$ and $\frac{c}{d}$ be consecutive entries in
$F_l$, and let $\U=\{U_i\}_{i \in I}$ and $\{V_j\}_{j \in J}$ be as above.

\begin{claim} We have $\mathcal{C}_{\frac{a+c}{b+d}}
\subset \mathcal{U}^{\vdash}$.\end{claim}
\noindent
\textit{Proof.} Let $\mathcal{F}$ be an indecomposable sheaf in 
$\mathcal{C}_{\frac{a+c}{b+d}}$. Since $\frac{a}{b}
< \frac{a+c}{b+d}$ we have $\mathcal{F} \in 
\mathcal{U}^{\Vdash}$. Consider the morphism $\sigma_\mathcal{F}:
\Sigma_{\U}\mathcal{F} \to \mathcal{F}$ and let $\mathcal{G}$ be the image
of $\sigma_\mathcal{F}$. 
Thus $\sigma_\mathcal{G}: \Sigma_{\U}\mathcal{G} =\Sigma_\U\mathcal{F} \to 
\mathcal{G}$ is an epimorphism. We may assume that $\mathcal{G}$ is 
indecomposable. We claim that $\mathrm{Hom}(\mathcal{G}, U_i) \neq 0$
for some $i \in I$. Indeed, if not then by Serre duality 
$\mathrm{Ext}^1(U_i,\mathcal{G}) =0$ for all $i$ and using
the Riemann-Roch theorem we obtain
\begin{equation}\label{E:C81}
\begin{split}
\sum_i \mathrm{dim\;Hom}(U_i,\mathcal{G})&=
\sum_i \langle U_i, \mathcal{G} \rangle\\
&=\mathbf{r}(U_i)\mathbf{d}(\mathcal{G})-
\mathbf{d}(U_i)\mathbf{r}(\mathcal{G})\\
&=b\mathbf{d}(\mathcal{G})-a\mathbf{r}(\mathcal{G}).
\end{split}
\end{equation}
On the other hand, we have
\begin{equation}\label{E:C82}
\sum_i \mathrm{dim\;Hom}(U_i,\mathcal{G})=
\frac{\mathbf{r}(\Sigma_\U\mathcal{G})}{b} \leq 
\frac{\mathbf{r}(\mathcal{G})}{b}.
\end{equation}
Combining (\ref{E:C81}) and (\ref{E:C82}) we deduce that $\mu(\mathcal{G})
\geq \frac{a}{b} + \frac{1}{b^2}$. But since $\mathcal{G} 
\hookrightarrow \mathcal{F}$ we have 
$$\mu(\mathcal{G}) \leq \mu(\mathcal{F})
=\frac{a+c}{b+d}=\frac{a}{b}+
\frac{bc-ad}{bd}=\frac{a}{b}+
\frac{1}{b(b+d)},$$
a contradiction. Thus there exists a nonzero map $\gamma :G \to U_i$. 
But then $\gamma$ is automatically an isomorphism, and $\sigma_\mathcal{G}$
is also an isomorphism.
This completes the proof
of Claim~8.1.\qed

\vspace{.1in}

Now, since $V_j$ is of rank $d$ and of degree $c$, the 
Riemann-Roch theorem implies that
\begin{equation}\label{E:C88}
\mathbf{d}(R_\U V_j)=c+(bc-ad)a=a+c,
\end{equation}
\begin{equation}\label{E:C89}
\mathbf{r}(R_\U V_j)=d+(bc-ad)b_{i}=b+d.
\end{equation}
In particular, $\mu(R_\U V_j)=\frac{a+c}{b+d}$. Hence
from the above claim we deduce that $\{R_\U V_j\}_{j \in J}$ forms a complete
set of simple sheaves in $\mathcal{C}_{\frac{a+c}{b+d}}$
and ii) follows. Note that i) for $\frac{a+c}{b+d}$ is now a
consequence of (\ref{E:C88}) and (\ref{E:C89}).

It remains to prove iii) for the pairs $\{\frac{a}{b},
\frac{a+c}{b+d}\}$ and $\{\frac{a+c}{b+d},
\frac{c}{d}\}$. This is taken care of by the following claims.

\begin{claim} For all $j \in J$ we have
\begin{equation}\label{E:dimL81}
\sum_i \mathrm{dim\;Ext}^1(R_{\U}V_j,U_i)=\frac{p}{c_j},\qquad
\mathrm{dim\;Ext}^1(R_{\U}V_j,U_i) \leq 1\quad\mathrm{for\;all\;}i.
\end{equation}
\end{claim}
\noindent
\textit{Proof.}
 By definition, we have
\begin{equation*}
[R_{\U}V_j]=[V_j]+\sum_k \mathrm{dim\;Ext}^1(V_j,U_k)[U_k]
=[V_j]-\sum_k\langle V_j,U_k\rangle [U_k]
\end{equation*}
since $V_j \in \mathcal{C}_{q'}$, $U_k \in \mathcal{C}_q$ and $q<q'$, so that
$\mathrm{Hom}(V_j,U_k)=0$. Observe that $R_\U V_j$ is indecomposable by 
Theorem~8.1. Hence $R_\U V_j \in \mathcal{C}_{q''}$ for some $q<q^{''}<q'$.
Therefore,
\begin{equation*}
\mathrm{dim\;Ext}^1(R_\U V_j,U_i)=-\langle R_{\U}V_j, U_i\rangle
=-\langle V_j,U_i\rangle + \sum_k \langle V_j, U_k \rangle \langle U_k, 
U_i\rangle.
\end{equation*}
But by Theorem~8.2 and Lemma~6.1, we have
\begin{equation}\label{E:888}
\langle U_k, U_i \rangle =\begin{cases} 1 & \mathrm{if}
\; U_k=U_i\\
-1 & \mathrm{if}\ \tau(U_k)=U_i\\
0 & \mathrm{otherwise} \end{cases}.
\end{equation}
Hence $\mathrm{dim\;Ext}^1(R_\U V_j,U_i)=\langle V_j, \tau^{-1}U_i\rangle$,
and (\ref{E:dimL81}) is proved.\qed

\vspace{.1in}

Similarly, let $K \subset J$ be such that $\{V_k\;|k \in K\}$ forms a single
$\tau$-orbit of size $p$.
\begin{claim} For all $j \in J$ we have
\begin{equation}\label{E:dimL82}
\sum_k \mathrm{dim\;Ext}^1(V_j,R_{\U}V_k)=\frac{p}{c_j},\qquad
\mathrm{dim\;Ext}^1(V_j,R_{\U}V_k) \leq 1\quad\mathrm{for\;all\;}k.
\end{equation}
\end{claim}
\noindent
\textit{Proof.} As in the previous lemma,
\begin{equation*}
\mathrm{dim\;Ext}^1(V_j,R_\U V_k)=-\langle V_j, R_\U V_k\rangle
=-\langle V_j,R_\U V_k\rangle + \sum_i \langle V_j, U_i \rangle \langle V_k, 
U_i\rangle.
\end{equation*}
Since the $\tau$-orbit of $V_k$ is of size $p$, the induction hypothesis
on $\mathrm{Ext}^1(V_k, U_i)$ implies that there exists $i(k) \in I$ such that
$$\langle V_k,U_{i(k)} \rangle=1, \qquad \langle V_k,U_{i'} \rangle=0
\;\text{for\;all\;}i' \neq i(k).$$
As a consequence,
$$\mathrm{dim\;Ext}^1(V_j,R_\U V_k)=-\langle V_j,R_\U V_k\rangle +
\langle V_j, U_{i(k)}\rangle$$
Finally, it is easy to see that the assignement $k \mapsto i(k)$ is a bijection
between $K$ and $I$. Using (\ref{E:888}) we deduce equation (\ref{E:dimL82}).

This concludes the proof of Claim~8.3 and of Proposition~8.3. \qed

\vspace{.1in}

We are now in position to prove Proposition~8.2. We argue again by induction.
By definition, $\U^{\geq 0}_{\p,\l}$ contains the classes of all the simple
sheaves in $\mathcal{C}_{\infty}$ and the classes of simple sheaves
$\O(i\vec{\omega})$. Assume that $\U^{\geq 0}_{\p,\l}$ contains the classes
of all simple sheaves of $\mathcal{C}_q$ if $q >0$ and $q$ appears in $F_l$.
Let $\frac{a}{b}$ and $\frac{c}{d}$ be two consecutive entries in $F_l$.
Pick a collection $\U=\{U_i\}_{i \in I}\subset \U^{\geq 0}_{\p,\l}$ 
of simple sheaves in $\mathcal{C}_{\frac{a}{b}}$ forming a single $\tau$-orbit
of size $p$. Let $V$ be any simple sheaf of $\mathcal{C}_{\frac{c}{d}}$.
From Proposition~8.2, iii), we deduce that there exists $I_0
\subset I$ such that $\Sigma_\U V= \bigoplus_{i \in I_0} U_{i}$. Recall that
$R_\U V$ is the universal extension of $V$ with respect to $\mathcal{U}$.
Moreover, it follows from Proposition~8.3, iii) and Serre duality that
$$\mathrm{dim\;Hom\;}(U_i,R_{\U} V)=\langle U_i, R_\U V\rangle
=-\langle R_{\U} V, \tau U_i \rangle \in \{0,1\}.$$
We conclude that for any subset $J \subset I_0$ there exists a unique
indecomposable object $O_J$ fitting in an exact sequence
$$0 \to \bigoplus_{i \in J} U_i \to O_J \to V \to 0$$
(specifically, $O_J \simeq R_\U V / \bigoplus_{j \not\in J} U_j$). 
\begin{claim} $[O_J] \in \U^{\geq 0}_{\p,\l}$
for any $J \subset I_0$.\end{claim}
\noindent
\textit{Proof.} We argue by induction. Assume that 
$[O_J] \in \U^{\geq 0}_{\p,\l}$ for any subset $J$ of size at most $k$, and
take $J\subset I_0$ with $|J|=k$. Observe that for any proper subset 
$J \subset I$
the class $[\bigoplus_{j \in J} U_j]$ is in $\U^{\geq 0}_{\p,\l}$ (this can
easily be deduced from the fact that the Serre subcategory of $Coh(\X)$
generated by $\U$ is equivalent to the category of representations of the 
cyclic quiver $A_{p-1}^{(1)}$). Thus, 
\begin{equation}\label{E:C84}
\U^{\geq 0}_{\p,\l} \owns [\bigoplus_{j \in J} U_j] \cdot [V]= c_{J}[O_J] 
+ \sum_{K \subset J, K \neq J}
c_K [O_K \oplus \bigoplus_{l \not\in K}U_l]
\end{equation}
for some nonzero $c_J, c_K \in \C_v$. By the induction hypothesis, $[O_K]
\in \U^{\geq 0}_{\p,\l}$ for all $K$, and thus 
$[O_K \oplus \bigoplus_{l \not\in K}U_l] = c [O_K] \cdot 
[\bigoplus_{l \not\in K}U_l] \in \U^{\geq 0}_{\p,\l}$. We conclude that
$[O_J] \in \U^{\geq 0}_{\p,\l}$, proving the claim. \qed

\vspace{.1in}

In particular, taking $J=I_0$ we see that
$[R_\U V] \in \U^{\geq 0}_{\p,\l}$. 
Finally, by Proposition~8.2 ii), 
$R_\U V$ is a simple sheaf in $\mathcal{C}_{\frac{a+c}{b+d}}$ and every
simple sheaf in $\mathcal{C}_{\frac{a+c}{b+d}}$ arises in this way. Thus
$\U^{\geq 0}_{\p,\l}$ contains the classes of all simple sheaves in
$\mathcal{C}_q$ for $q>0$ arising in $F_{l+1}$. This completes the induction
and the proof of Proposition~8.3. \qed 
\vspace{.2in}

\paragraph{\textbf{8.4.}} In this paragraph, we construct certain elements
in $\U^{\geq 0}_{\p,\l}$ corresponding to imaginary root vectors of $\widehat{\n}$.
We first introduce some useful notation. For any $q \in \mathbb{Q}$ we denote by
$\mathcal{C}_{>q}$ (resp. $\mathcal{C}_{<q}$) the additive subcategory of $Coh(\X)$
generated by $\mathcal{C}_{q'}$ for $q' >q$ (resp. $q' <q$). Also, if $\mathcal{F}$
is any coherent sheaf and $u \in \H_{Coh(\X)}$ then we let $c_\mathcal{F}(u)$ be
the coefficient of $[\mathcal{F}]$ in $u$. Finally, we
will say that a sheaf $\mathcal{F}\in \mathcal{C}_q$ is a \textit{$\lambda$-sheaf}
if $\mathcal{F}=\mathcal{X}_1 \oplus \cdots \oplus \mathcal{X}_l$ where 
$\mathcal{X}_1, \ldots \mathcal{X}_l$ are indecomposable sheaves in
$\mathcal{C}_q$ corresponding, under the equivalence $\mathcal{C}_q \simeq 
\mathcal{C}_{\infty}$ to simple torsion sheaves of degree $\lambda_1 \delta,
\ldots, \lambda_l \delta$ respectively.

\vspace{.1in}

Define elements $T^{(q)}_r \in \U^{\geq 0}_{\p,\l}$ for $q \in \mathbb{Q}$, $q>0$
and $r \geq 1$ as follows. Set $T^{\infty}_r=T_r$. Assume that $T_r^{(q)}$ has been 
defined for all $q$ belonging to the sequence $F_l$. Let $\frac{a}{b}$ and 
$\frac{c}{d}$ be consecutive entries in $F_l$, and let $\U=\{U_i\}_{i \in \Z/p\Z}$ be
any collection of simple sheaves in $\mathcal{C}_{\frac{a}{b}}$ forming a single
$\tau$-orbit. We put
$$T_r^{(\frac{a+c}{b+d})}=[U_1]^r \cdot [U_2]^r \cdots [U_p]^r \cdot 
T_r^{(\frac{c}{d})}.$$
This definition depends on the choice (and ordering) of the $U_i$ but this will
not be important for us.
For any $q>0$ and any partition
$\lambda=(\lambda_1 \geq \cdots \geq \lambda_l)$ we define
$T_\lambda^{(q)}=T^{(q)}_{\lambda_1} \cdots T^{(q)}_{\lambda_l}$.

\begin{claim}\label{C:909} We have $c_{\mathcal{F}}(T_\lambda^{(q)}) \neq 0$ for any 
$\lambda$-sheaf $\mathcal{F}$ in $\mathcal{C}_q$.
\end{claim}
\noindent
\textit{Proof.} Let us first consider the case when $\lambda=\lambda_1$, and
argue by induction. The statement is clear for $q \in F_1$. Assume that it holds
for all $q \in F_l$ and let $\frac{a}{b}$, $\frac{c}{d}$ and
 $\{U_i\}_{i \in \Z/p\Z}$ have the same meaning as above. Let $\mathcal{F}
\in \mathcal{C}_{\frac{c}{d}}$ be a $\lambda$-sheaf. 
From Proposition~8.3, iii) we see that 
$\Omega_{\U}(\mathcal{F})=U_1 ^{\oplus r} \oplus \cdots \oplus U_p^{\oplus r}$.
It remains to notice that
\begin{equation*}
c_{[U_1 ^{\oplus r} \oplus \cdots \oplus U_p^{\oplus r}]}
([U_1]^r \cdot [U_2]^r \cdots [U_p]^r) \neq 0.
\end{equation*} 
This proves the claim for $\lambda$ of length one. The general case easily follows.
\qed

\vspace{.1in}

By Section 2.4 and Section 4.2, the subalgebra $\mathbf{C}_q$ of 
$\U^{\geq 0}_{\p,\l}$ generated by the simple sheaves in $\mathcal{C}_q$ is 
isomorphic to $\U^+_v(\widehat{\mathfrak{sl}}_{p_1}) \otimes \cdots \otimes 
\U^+_v(\widehat{\mathfrak{sl}}_{p_n})$. Let $\mathcal{B}$ be a basis of the latter
 algebra consisiting of weight vectors, and let $\mathcal{B}^{(q)}$ be the 
corresponding basis of $\mathbf{C}_q$.

\begin{lem}\label{L:finito} For any $q >0$ and any sheaf $\mathcal{F}$ appearing in
$T_r^{(q)}$ we have $[\mathcal{F}]=r\delta_q \in K_0(Coh(\X)) \simeq 
\widehat{Q}$. \end{lem}
\noindent
\textit{Proof.} We argue again by induction. Assume the statement in the
Lemma for all $q$ appearing in $F_l$ and let $\frac{a}{b}, \frac{c}{d}$ be
consecutive entries in $F_l$ and let $\U=\{U_i\}_{i \in \Z/p\Z}$ be any
collection of simple sheaves in $\mathcal{C}_{\frac{a}{b}}$ forming a single
$\tau$-orbit. It is enough to observe that
$$\delta_{\frac{a+c}{b+d}}=[U_1]+ \cdots + [U_p] + \delta_{\frac{c}{d}}$$
since the right mutation $R_{\U}$ defines an equivalence $\mathcal{C}_{\frac{c}{d}} \simeq \mathcal{C}_{\frac{a+c}{b+d}}$. \qed

\begin{prop}\label{P:887} The elements $\{b^{(q)} \cdot T_\lambda^{(q)}\;|b^{(q)} \in 
\mathcal{B}^{(q)},\; \lambda \in \Pi\}$ are linearly independent.\end{prop}
\noindent
\textit{Proof.} Let $Coh^{ss}(\X)$ denote the set of semistable coherent sheaves
on $\X$. Set
\begin{align*}
I'_q&=\bigoplus_{\substack{\mu(\mathcal{F})=q\\\mathcal{F} \not\in Coh^{ss}(\X)}}
\C[\mathcal{F}] \subset \H_{Coh(\X)},\\
I^{ss}_q&=\bigoplus_{\substack{\mu(\mathcal{F})=q\\\mathcal{F} \in Coh^{ss}(\X)}}
\C[\mathcal{F}] \subset \H_{Coh(\X)}.
\end{align*}
If $[\mathcal{F}] \in I'_q$ then $\mathcal{F}=\mathcal{F}^- \oplus \mathcal{F}^0 
\oplus \mathcal{F}^+$ where $\mathcal{F}^- \in \mathcal{C}_{<q}$ and
$\mathcal{F}^+ \in \mathcal{C}_{>q}$ are both nonzero and where $\mathcal{F}^0
\in \mathcal{C}_q$.

\begin{lem}\label{L:99} $I'_q \oplus I_q^{ss}$ is a subalgebra of $\H_{Coh(\X)}$, and
$I'_q$ is an ideal in $I'_q \oplus I_q^{ss}$.\end{lem}
\noindent
\textit{Proof.} The first statement is obvious. For the second, observe that if
$[\mathcal{G}] \in I_q^{ss}$ and $[\mathcal{F}]=[\mathcal{F}^- \oplus \mathcal{F}^0
\oplus \mathcal{F}^+] \in I'_q$ then by the HN filtration there can be no
monomorphism $\mathcal{F}^- \oplus \mathcal{F}^0 \oplus \mathcal{F}^+ 
\hookrightarrow \mathcal{G}$
or epimorphism $\mathcal{G} \twoheadrightarrow  
\mathcal{F}^- \oplus \mathcal{F}^0 \oplus \mathcal{F}^+$.\qed

\vspace{.1in}

We will write $\pi_{ss}$ for the projection.
 $I^{ss}_q \oplus I_q' \to I^{ss}_q$. Note that $T_r^{(q)} \in I_q^{ss} \oplus
I_q'$. By Lemma~\ref{L:99} we have, for any partition $\lambda=(\lambda_1 \geq
\cdots \geq \lambda_l)$,
\begin{equation}\label{E:851}
\pi_{ss}(T_\lambda^{(q)})= \pi_{ss}(T^{(q)}_{\lambda_1}) \cdots
\pi_{ss}(T^{(q)}_{\lambda_l}).
\end{equation}

From this and Claim~\ref{C:909} we deduce the following result.
\begin{claim}\label{C:987}
 Let $\mu >\lambda$ and let $\mathcal{F}$ be any $\mu$-sheaf. Then
$c_{\mathcal{F}}(T_{\lambda}^{(q)}) =0$.
\end{claim}

Now, to prove that the set $\{b^{(q)} \cdot T^{(q)}_{\lambda}\}$ is linearly
independent, it is enough to consider their image under the projection $\pi_{ss}$.
Let us call a sheaf $\mathcal{G} \in \mathcal{C}_q$ 
\textit{exceptional} if it corresponds, under the equivalence $\mathcal{C}_q \simeq
\mathcal{C}_{\infty}$, to a torsion sheaf supported at exceptional points. In
particular, if $\mathcal{B}^{(q)} \owns b^{(q)}=\sum_i c_i [\mathcal{G}_i]$ then
all $\mathcal{G}_i$ are exceptional. Observe that if $\mathcal{G}$ is exceptional
and if $\mathcal{F}$ is a $\lambda$-sheaf then (up to a scalar)
$[\mathcal{G}] \cdot [\mathcal{F}]
=c[\mathcal{F} \oplus \mathcal{G}]$. 
Proposition~\ref{P:887} now easily follows from Claim~\ref{C:987} and from the
fact that the elements of $\mathcal{B}^{(q)}$ are linearly independent.\qed
\vspace{.2in}

\paragraph{\textbf{8.6.}} Let $\mathcal{I}$ denote the collection of all tuples
$$\sigma=\{q_1, \ldots, q_l; d_1, \ldots, d_l\;|\; q_i, d_i \in \mathbb{Q}^+,
\; q_1 > q_2 > \cdots > q_l\}.$$
We introduce a total order on $\mathcal{I}$ as follows. $\sigma=\{q_i; d_i\}
\succ \sigma'=\{q'_i; d'_i\}$ if there exists $j \geq 0$ such that
$q_1=q'_1, d_1=d'_1, \ldots q_j=q'_j, d_j=d'_j$ and
$q_{j+1} > q'_{j+1}$ or $q_{j+1}=q'_{j+1}$ and $d_i>d'_{i+1}$.

\vspace{.1in}

Let us associate to any coherent sheaf $\mathcal{F}$ the element $\sigma(\mathcal{F})
\in \mathcal{I}$ defined as
$$\sigma(\mathcal{F})=\{q_1, \ldots, q_l;\mathbf{d}(\mathcal{F}_1), \ldots, 
\mathbf{d}(\mathcal{F}_l)\}$$
where $\mathcal{F}=\mathcal{F}_1 \oplus \cdots \oplus \mathcal{F}_l$ with
$\mathcal{F}_i \in \mathcal{C}_{q_i}$ and $q_1 > \cdots > q_l$.

\begin{prop}\label{P:final} Fix $q_1 > \cdots > q_l >0$ and for $1 \leq i \leq l$ let
$b^{(q_i)} \in \mathcal{B}^{(q_i)}$ and $\lambda_i \in \Pi$. We have
$$(b^{(q_1)} \cdot T^{(q_1)}_{\lambda_1}) \cdots 
(b^{(q_l)} \cdot T^{(q_l)}_{\lambda_l}) \in 
\pi_{ss}(b^{(q_1)} \cdot T^{(q_1)}_{\lambda_1}) \cdots 
\pi_{ss}(b^{(q_l)} \cdot T^{(q_l)}_{\lambda_l}) \oplus 
\bigoplus_{\sigma(\mathcal{F})\succ \gamma} \C[\mathcal{F}],$$
where $\gamma=\{q_1, \ldots, q_l; \mathbf{d}(b^{(q_1)}\cdot T_{\lambda_1}^{(q_1)}),
\ldots, \mathbf{d}(b^{(q_l)}\cdot T_{\lambda_l}^{(q_l)})\}$.\end{prop}
\noindent
\textit{Proof.} We give a proof by induction. The statement is clear for $l=1$.
Assume that the Proposition is true for a product of $l-1$ terms, and let us set
$$A= (b^{(q_1)} \cdot T^{(q_1)}_{\lambda_1}) \cdots 
(b^{(q_{l-1})} \cdot T^{(q_{l-1})}_{\lambda_{l-1}}).$$
The next statement is a simple consequence of the HN filtration and of 
the definition of $\sigma$ :

\begin{lem} Let $\mathcal{F}$ and $\mathcal{G}$ be any coherent sheaves in
$\mathcal{C}_{>0}$. Then 
$$[\mathcal{F}] \cdot [\mathcal{G}] \in 
\bigoplus_{\sigma(\mathcal{H}) \succ \sigma(\mathcal{F})} \C[\mathcal{H}].$$
\end{lem}

\vspace{.1in}

In particular, we deduce that
$$A \cdot (b^{(q_l)} \cdot T_{\lambda_l}^{(q_l)}) \in 
\pi_{ss}(b^{(q_1)} \cdot T^{(q_1)}_{\lambda_1}) \cdots 
\pi_{ss}(b^{(q_{l-1})} \cdot T^{(q_{l-1})}_{\lambda_{l-1}})
\cdot (b^{(q_{l})} \cdot T^{(q_{l})}_{\lambda_l}) \oplus 
\bigoplus_{\sigma(\mathcal{F})\succ \gamma} \C[\mathcal{F}].$$
But 
$$b^{(q_l)}\cdot T^{(q_l)}_{\lambda_l} \in \pi_{ss}(b^{(q_l)}\cdot 
T^{(q_l)}_{\lambda_l}) \oplus \bigoplus_{\sigma(\mathcal{H}) \succ \gamma'}\C
[\mathcal{H}],$$
where $\gamma'=\{q_l; \mathbf{d}(b^{(q_l)}\cdot T^{(q_l)}_{\lambda_l})\}$.
The Proposition easily follows.\qed

\vspace{.1in}

\begin{cor}\label{Co:final} The elements of $\U^{\geq 0}_{\p,\l}$
$$\{(b^{(q_1)} \cdot T^{(q_1)}_{\lambda_1}) \cdots 
(b^{(q_l)} \cdot T^{(q_l)}_{\lambda_l})\;|q_1 > \cdots > q_l, b^{(q_i)} \in 
\mathcal{B}^{(q_i)}, \lambda_i \in \Pi\}$$
of $\U^{\geq 0}_{\p,\l}$ are linearly independent.\end{cor}
\noindent
\textit{Proof.} By Proposition~\ref{P:final} it is enough to prove that, for any
\textit{fixed} $q_1 > \cdots > q_l$, the set
$$\{\pi_{ss}(b^{(q_1)} \cdot T^{(q_1)}_{\lambda_1}) \cdots 
\pi_{ss}(b^{(q_{l})} \cdot T^{(q_{l})}_{\lambda_{l}})\;|\;b^{(q_i)} \in 
\mathcal{B}^{(q_i)}, \lambda_i \in \Pi\}$$
is linearly independent. But by the HN filtration again, we have (up to a scalar)
\begin{equation*}
\begin{split}
\pi_{ss}&(b^{(q_1)} \cdot T^{(q_1)}_{\lambda_1}) \cdots 
\pi_{ss}(b^{(q_{l})} \cdot T^{(q_{l})}_{\lambda_{l}})\\
&=\sum_{\mathcal{H}_i \in \mathcal{C}_{q_i}}c_{\mathcal{H}_1}
(\pi_{ss}(b^{(q_1)} \cdot T^{(q_1)}_{\lambda_1})) \cdots
c_{\mathcal{H}_l}
(\pi_{ss}(b^{(q_l)} \cdot T^{(q_l)}_{\lambda_l}))[\mathcal{H}_1 \oplus
\cdots \oplus \mathcal{H}_l].
\end{split}
\end{equation*}
 The corollary is now a consequence of Proposition~\ref{P:887}. \qed

\vspace{.2in}

We are finally ready to prove (\ref{E:dim81}). Recall the notations of 
Section~8.2. Let us set $\widehat{\n}_{>0}=\bigoplus_{q>0} \widehat{\n}_q$ and let
$\widehat{\n}_0$ be the subalgebra generated by the root vectors of weight
$\alpha_*, \tau(\alpha_*), \ldots, \tau^{p-1}(\alpha_*)$. Similarly, let us denote
by $\U_{\p,\l}^{(q)} \subset \U^{\geq 0}_{\p,\l}$ the vector space spanned
by the elements $\{b^{(q)} \cdot T_{\lambda}^{(q)}\;|\; b \in \mathcal{B}^{(q)},
\lambda \in \Pi\}$, and by $\U^0_{\p,\l}$ the subalgebra generated by
the line bundles $[\O], [\O(\vec{\omega})], \ldots, [\O((p-1)\vec{\omega})]$.

\vspace{.1in}

By the PBW theorem we have, for any $\alpha \in \widehat{Q}$,
\begin{equation*}
\begin{split}
\mathrm{dim\;}\U^{\geq 0}(\widehat{\n})[\alpha] \leq & 
\mathrm{dim\;}\U(\widehat{\n}_0 \oplus \widehat{\n}_{>0}) [\alpha]\\
=& \sum_{\substack{\alpha=\alpha_0+ \cdots + \alpha_l \\ q_1> \cdots >q_l >0}}
\mathrm{dim\;}\U(\widehat{\n}_0)[\alpha_0] \cdot\mathrm{dim\;}\U(\widehat{\n}_{q_1})
[\alpha_1] 
\cdots \mathrm{dim\;}\U(\widehat{\n}_{q_l})[\alpha_l].
\end{split}
\end{equation*}
On the other hand, from Corollary~\ref{Co:final} we deduce
$$\mathrm{dim}\;\U^{\geq 0}_{\p,\l}[\alpha]
\geq \sum_{\substack{\alpha=\alpha_0+ \cdots + \alpha_l \\ q_1> \cdots >q_l >0}}
\mathrm{dim\;}\U^0_{\p,\l}[\alpha_0] \cdot \mathrm{dim\;}\U^{(q_1)}_{\p,\l}[\alpha_1]
\cdots \mathrm{dim\;}\U^{(q_l)}_{\p,\l}[\alpha_l].$$
To conclude, observe that, by Corollary~\ref{cor:wei} and Section~7.4, 
we have for any $\beta$
$$\mathrm{dim\;}\U(\widehat{\n}_q)[\beta]=\mathrm{dim\;}\U(\widehat{\mathfrak{t}})
[\beta_\infty]=\mathrm{dim\;}\U^{(\infty)}_{\p,\l}[\beta_\infty],$$
where $\beta_\infty$ corresponds to $\beta$ under the identification $\widehat{\Delta}_q
\simeq \widehat{\Delta}_{\infty}$, and that by Lemma~\ref{L:finito} and
Proposition~\ref{P:887} we have
$$\mathrm{dim\;}\U^{(\infty)}_{\p,\l}[\beta_\infty]=\mathrm{dim\;}\U^{(q)}_{\p,\l}
[\beta].$$

\vspace{.2in}

\section{Conjectures}

\vspace{.1in}

\paragraph{\textbf{9.1.}}We first recall the Kac conjectures for representations of quivers. Let $Q=(I,\Omega)$
be an arbitrary quiver with no edge loops. Let $\g$ be the Kac-Moody algebra associated to
$Q$, let $\h \subset \g$ be a Cartan subalgebra and let $\Delta^+ \subset \h^*$ denote the set of
positive roots. We view the dimension $\mathbf{dim}\;M \in \mathbb{N}^I$ of some object
$M \in Rep_k(Q)$ as an element of $\h^*$ via the map $\mathbb{N}^I \to \h^*, (n_i)_i \mapsto 
\sum_i n_i \alpha_i$.

\begin{theo}[Kac] The following hold~:\\
\begin{enumerate}
\item[i)] There exists an indecomposable object in $Rep_k(Q)$ of dimension $\alpha$
if and only if $\alpha \in \Delta^+$. Moreover, the set of such indecomposables
forms a $1-\frac{(\alpha,\alpha)}{2}$-parameter family with a unique component 
of maximal dimension.
\item[ii)] For every root $\alpha$ there exists a polynomial $P_{\alpha}(v) \in 
\Z[v]$ with leading term $v^{1-\frac{(\alpha,\alpha)}{2}}$ such that the
number of absolutely indecomposable representations of dimension $\alpha$
over the finite field $ \mathbb{F}_q$ is $P_\alpha(q)$.
\end{enumerate}
\end{theo}

\begin{conj}[Kac] We have $P_{\alpha}(v) \in \N[v]$ and $P_{\alpha}(0)
=\mathrm{dim}\;\g_{\alpha}$.
\end{conj}

This conjecture is trivial for finite type Dynkin quivers, proved for affine
quivers and remains open in all other cases (see however \cite{CB-VdB} for
a recent important progress).

\vspace{.2in}

\paragraph{\textbf{9.2.}} The results of this paper naturally lead us to the following
variant of the above conjectures. Let $\X_{\p,\l}$ be a weighted
projective line and let $\mathfrak{Lg}$ be the corresponding loop Kac-Moody
algebra. Recall the identification $K_0(Coh(\X_{\p,\l})) \simeq 
\h^*_{\mathfrak{Lg}}$ of Proposition~5.1.

\begin{conj} The following hold~:
\begin{enumerate}
\item[i)] There exists an indecomposable coherent sheaf $\mathcal{F} \in 
Coh(\X_{\p,\l})$ if and only if $[\mathcal{F}] \in \h^*_{\mathfrak{Lg}}$
is a positive root. Moreover, the set of such indecomposables
forms a $1-\frac{(\alpha,\alpha)}{2}$-parameter family with a unique component 
of maximal dimension.

\item[ii)]  For every root $\alpha$ there exists a polynomial $P_{\alpha}(v) \in 
\N[v]$ with leading term $v^{1-\frac{(\alpha,\alpha)}{2}}$ such that the
number of absolutely indecomposable coherent sheaves of dimension $\alpha$
over the finite field $ \mathbb{F}_q$ is $P_\alpha(q)$, and we have
 $P_{\alpha}(0)
=\mathrm{dim}\;\mathfrak{Lg}_{\alpha}$.
\end{enumerate}
\end{conj}

\paragraph{}Conjecture i) was first made by Crawley-Boevey in relation to the
Deligne-Simpson problem (\cite{CB}). It is easy to check (using the 
descriptions of $Coh(\X_{\p,\l})$ in Sections 7 and 8) that i) and ii) are
true when $\X_{\p,\l}$ is of genus $g \leq 1$.

\vspace{.1in}

\section{Appendix}

\vspace{.1in}

In this appendix we give the details of the calculations of Lemma~\ref{L:666}.
The computations of Section 6.9 are conducted in an exactly similar fashion
and we leave them to the reader. For simplicity we drop the index $(s)$ 
throughout in the notation (as in $p_s$ for instance).

\vspace{.2in}
\paragraph{\textit{Proof of Lemma~\ref{L:666}.}} Suppose that
(\ref{E:app1}) holds for $r_1=1$ and arbitrary $r_2$. Set $\widetilde{T}_r
=\frac{r}{[r]}T_r$. Then
$$[S_{1,r_2}]= [[S_1],\widetilde{T}_{r_2}]=(-1)^{r_2}ad^{r_2}
\widetilde{T}_{r_2} ([S_1]).$$
Since $[T_r,T_s]=0$ for all $r,s$ we deduce
\begin{equation*}
\begin{split}
[[S_{1,r_2}],\widetilde{T}_{r_1}]&=(-1)^{r_2}[ad^{r_2}\widetilde{T}_{1}([S_1]),
\widetilde{T}_{r_1}]\\
&=(-1)^{r_2}ad^{r_2}\widetilde{T}_{1}([[S_1],\widetilde{T}_{r_1}])\\
&=(-1)^{r_2}ad^{r_2}\widetilde{T}_1([S_{1,r_1}])\\
&=S_{1,r_1+r_2}.
\end{split}
\end{equation*}
Thus it is enough to prove (\ref{E:app1}) for $r_2=1$. Let us introduce some
notations. We label the boxes of a partition $\lambda$ according to their 
height. For instance, the partition $(32^21)$ is labelled
$$
\centerline{\epsfbox{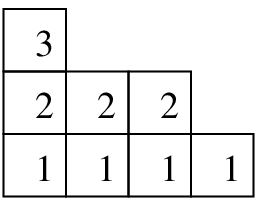}}
$$
Let us write $\lambda \subset \nu$ if $\nu$ is obtained from $\lambda$ by 
addition of a single box, and we write $\nu\backslash \lambda=\boxed{i}$ 
if the 
corresponding box is labelled by $i$. We also denote by $l(\lambda)_{\geq i}$
(resp. $l(\lambda)_{>i}$) the number of parts $\lambda_j$ of $\lambda=
(\lambda_1 \geq \lambda_2 \geq \cdots)$ with value $\lambda_j \geq i$ (resp.
$\lambda_j >i$. Finally, if $\nu \backslash \lambda =\boxed{i}$ we set
$$p(\nu,\lambda)=\frac{v^{2l(\lambda)_{\geq i}}-v^{2l(\lambda)_{>i}}}{v^2-1}.$$

\begin{lem}\label{L:A1} We have
$$\widetilde{T}_1 \widetilde{T}_r=\sum_{\lambda, |\lambda|=r}n(l(\lambda)-1)
\sum_{\nu \supset \lambda} p(\nu,\lambda) [|\nu)].$$
\end{lem}
\noindent
\textit{Proof.} It is enough to prove that
$$[|\lambda)] \cdot [|1)] = \sum_{\nu \supset \lambda} p(\nu,\lambda)[|\nu)].$$
By Section 4.2, we can assume $p=1$. Let $(M)\nu,x)$ be a representation of
$\mathcal{C}_1$ of type $\nu$. The space $N=\mathrm{Ker}\;x$ has a natural 
filtration
$N=N_0 \supset N_1 \supset \cdots$ where $N_i=N \cap \mathrm{Im}\;x^i$. We have
$\mathrm{dim}\;N_i=l(\nu)_{>i}$. Let 
$\lambda \subset \nu$ and let $L \subset N$ be a line. Then $M_\nu /L$ is of
type $\lambda$ if and only if $L \subset N_{i-1}$ but $L \not\subset N_i$ where
$\nu \backslash \lambda=\boxed{i}$. Note that the number of such lines is 
$p(\nu,\lambda)$. \qed

A similar computation shows that 
$$[|\nu)] \cdot [S_1] = v^{2l(\lambda)} [|\nu) \oplus S_1],$$
$$[S_1] \cdot [|\nu)]= [|\nu) \oplus S_1] + \sum_{\gamma \triangleleft \nu}
[|\nu \backslash \gamma) \oplus [1; p\gamma +1)].$$
Hence,
\begin{equation}\label{E:A2}
\widetilde{T}_1\widetilde{T}_rS_1=\sum_{\lambda} n(l(\lambda)-1) 
\sum_{\nu \supset \lambda} p(\nu,\lambda) v^{2l(\nu)}[|\nu) \oplus S_1],
\end{equation}
\begin{equation}\label{E:A3.5}
\begin{split}
S_1\widetilde{T}_1\widetilde{T}_r=&\sum_{\lambda} n(l(\lambda)-1) 
\sum_{\nu \supset \lambda} p(\nu,\lambda) [|\nu) \oplus S_1]\\
&+ \sum_{\lambda} n(l(\lambda)-1) \sum_{\nu \supset \lambda} p(\nu,\lambda)
\sum_{\gamma \triangleleft \nu}[|\nu \backslash \gamma) \oplus [1;p\gamma+1)].
\end{split}
\end{equation}

A reasoning close to that of Lemma~\ref{L:A1} shows that
$$
[|\lambda) \oplus S_1] \cdot T_1=\sum_{\nu \supset \lambda} p(\nu,\lambda)
[|\lambda) \oplus S_1] + v^{2l(\lambda)_{>1}}[|\lambda) \oplus [1;p+1)],$$
\begin{equation*}
\begin{split}
[|1) \oplus S_1] \cdot [|\lambda)]=& \sum_{\nu \supset \lambda}p(\nu,\lambda)
[|\nu) \oplus S_1]\\
&\qquad+\sum_{\gamma \triangleleft \lambda}
\sum_{\sigma \supset \lambda \backslash \gamma} p(\sigma,\lambda 
\backslash \gamma) [|\sigma) \oplus [1;p\gamma+1)],
\end{split}
\end{equation*}
from which one deduces
\begin{equation}\label{E:A3}
\begin{split}
\widetilde{T}_r S_1 \widetilde{T}_1=& \sum_{\lambda} n(l(\lambda)-1)
\sum_{\nu \supset \lambda} v^{2l(\lambda)}p(\nu,\lambda) [|\nu) \oplus S_1]\\
&+ \qquad\sum_{\lambda} n(l(\lambda)-1)v^{2(l(\lambda) + l(\lambda)_{>1})}
[|\lambda) \oplus [1;p+1)],
\end{split}
\end{equation}
\begin{equation}\label{E:A4}
\begin{split}
\widetilde{T}_1 S_1 \widetilde{T}_r=& \sum_{\lambda} v^2n(l(\lambda)-1)
\sum_{\nu \supset \lambda} p(\nu,\lambda) [|\nu) \oplus S_1]\\
&+ \qquad\sum_{\lambda} v^2n(l(\lambda)-1)
\sum_{\gamma \triangleleft \lambda}\sum_{\sigma \supset \lambda \backslash 
\gamma} p(\sigma,\lambda \backslash \gamma) [|\sigma) \oplus [1;p\gamma+1)].
\end{split}
\end{equation}

We are now in position to prove Lemma~\ref{L:666}. We have
\begin{equation}\label{E:A5}
[S_{1,r},\widetilde{T}_1]=S_1\widetilde{T}_r \widetilde{T}_1 +
\widetilde{T}_r \widetilde{T}_1S_1 -\widetilde{T}_r S_1\widetilde{T}_1
-\widetilde{T}_1 S_1\widetilde{T}_r.
\end{equation}
The coefficient $c_\nu$ of $[|\nu) \oplus S_1]$ in (\ref{E:A5}) is
\begin{equation}\label{E:A6}
\begin{split}
c_{\nu} =& (1-v^2) \sum_{\lambda \subset \nu} n(l(\lambda)-1)p(\nu,\lambda)
+
\sum_{\lambda \subset \nu} n(l(\lambda)-1)p(\nu,\lambda)(v^{2l(\nu)}-
v^{2l(\lambda)})\\
=&(1-v^2) 
\sum_{\underset{\nu \backslash \lambda \neq 1}{\lambda \subset \nu}}
n(l(\nu)-1)p(\nu,\lambda) \\
&\qquad\qquad+ \sum_{\underset{\nu \backslash \lambda= 1}{\lambda \subset \nu}}
n(l(\nu)-2)(1-v^2)(1-v^{2(l(\nu)-1)})p(\nu,\lambda)\\
=&(1-v^2)n(l(\nu)-1)\sum_{\lambda \subset \nu} p(\nu,\lambda)\\
=&n(l(\nu))
\end{split}
\end{equation}
since 
$$\sum_{\lambda \subset \nu} p(\nu,\lambda)=\frac{v^{2l(\nu)}-1}{v^2-1}.$$

Similarly, the coefficient $c_{\nu,\gamma}$ of $[|\nu) \oplus [1;p\gamma+1)]$
in (\ref{E:A5}) is, if $\gamma >1$,
\begin{equation}\label{E:A7}
c_{\nu,\gamma}=\sum_{\lambda \subset \nu \cup \{\gamma\}} n(l(\lambda)-1)
p(\nu \cup \{\gamma\},\lambda)-v^2 \sum_{\sigma \subset \nu} n(l(\sigma)) 
p(\nu,\sigma).
\end{equation}
Observe that
$$ p(\nu,\sigma)=
\begin{cases}
p(\nu \cup \{\gamma\},\sigma \cup \{\gamma\} & \mathrm{if}\; \nu \backslash
\sigma > \boxed{\gamma},\\
v^{-2}p(\nu \cup \{\gamma\},\sigma \cup \{\gamma\} & \mathrm{if}\; 
\nu \backslash
\sigma < \boxed{\gamma},\\
p(\nu \cup \{\gamma\},\sigma \cup \{\gamma\} -v^{2(l(\nu)_{\geq \gamma})}&
 \mathrm{if}\; \nu \backslash
\sigma =\boxed{\gamma}
\end{cases}$$
substituting in (\ref{E:A7}) yields
\begin{equation}\label{E:A8}
\begin{split}
c_{\nu,\gamma}&=n(l(\nu)) \big\{(1-v^2) 
\sum_{\underset{\nu \backslash \sigma \geq \gamma}{\sigma \subset \nu}}
p(\nu \cup \{\gamma\},\sigma \cup \{\gamma\}) + v^{2(l(\nu)_{\geq \gamma}+1)}
\big\}\\
&=n(l(\nu)) \big\{ (1-v^2) \frac{v^{2(l(\nu)_{\geq \gamma}+1)}-1}{v^2-1}+
 v^{2(l(\nu)_{\geq \gamma}+1)}
\big\}\\
&=n(l(\nu)).
\end{split}
\end{equation}
Similar computation shows that (\ref{E:A8}) holds for $\gamma=1$ also.
Equations (\ref{E:A6}) and (\ref{E:A8}) together prove Lemma~\ref{L:666}.\qed

\small{

\vspace{4mm}
Olivier Schiffmann,\texttt{schiffma@dma.ens.fr},\\
Yale University, 10 Hillhouse Avenue, 06520-8283
New Haven, CT, USA, and\\
DMA, \'Ecole Normale Sup\'erieure, 45 rue d'Ulm, 75230 Paris Cedex 05-FRANCE

\begin{thebibliography}{99}
\bibitem[BK]{BK}
Baumann, P., Kassel, C., \emph{The Hall algebra of the category of coherent 
sheaves on the projective line.}, J. Reine Angew. Math. \textbf{533}
207-233, (2001). 
\bibitem[Be]{Beck}
Beck J., \emph{Braid group action and quantum affine algebras.}, 
Comm. Math. Phys. \textbf{165} no. 3, 555-568, (1994).
\bibitem[Bo]{Bondal}
Bondal A.I., \emph{Representations of associative algebras and coherent 
sheaves}, Izv. Akad. Nauk. SSSR Ser. Math. \textbf{53}, 25-44, (1989).
\bibitem[CB]{CB}
Crawley-Boevey B., \emph{On the Deligne-Simpson problem and weighted projective
lines}, in preparation (Talk at ICRA X 2002).
\bibitem[CB-VdB]{CB-VdB}
Crawley-Boevey B., Van den Bergh M., \emph{Absolutely indecomposable
representations and Kac-Moody algebras}, preprint RA/0106009.
\bibitem[E]{E}
Enriquez B., \emph{PBW and duality theorems for quantum groups and quantum 
current algebras}, preprint QA/9904113.
\bibitem[GL]{GL}
Geigle W., Lenzing H., \emph{A class of weighted projective curves arising in 
the representation theory of finite-dimensional algebras}, Singularities, 
representations of algebras and vector bundles, Springer LNM \textbf{1273},
 265-297, (1987).
\bibitem[Hal]{Hall}
Hall P., \emph{The algebra of partitions}, Proc. 4th Canadian Math. Congress, 
147-159, (1959).
\bibitem[Hap1]{Happ}
Happel, D. \emph{Triangulated categories in the representation theory of
finite dimensional algebras}, London Math. Soc. Lect. Notes \textbf{119}
(1988).
\bibitem[Hap2]{Happ2}
Happel, D., 
\emph{A characterization of hereditary categories with tilting object}, 
Invent. Math. \textbf{144} no. 2, 381-398, (2001). 
\bibitem[Kac]{Kac}
Kac V., \emph{Infinite dimensional Lie algebras}, Cambridge Univ. Press, (1990).
\bibitem[Kap]{K}
Kapranov, M,. \emph{Eisenstein series and quantum affine algebras},
Algebraic geometry 7,  
J. Math. Sci. \textbf{84} no. 5, 1311-1360, (1997).
\bibitem[La]{Lam}
Lamotke, K., \emph{Regular solids and isolated singularities.}
Adv. Lect. in Math., Braunschweig, (1986).
\bibitem[LT]{LT}
Leclerc B., Thibon J-Y., \emph{Littlewood-Richardson coefficients
and Kazhdan-Lusztig polynomials.}, Combinatorial methods in representation
theory (Kyoto, 1998), 155-220, Adv. Stud. Pure Math. \textbf{28}. 
\bibitem[Le]{L}
Lenzing H.,\emph{Curve singularities arising from the representation theory
of tame hereditary algebras}, Representation theory, I (Ottawa, Ont., 1984),
 Lecture Notes in Math. \textbf{1177}, 199-231.
\bibitem[LM]{LM}
Lenzing H., Meltzer H., \emph{Sheaves on a weighted projective line of 
genus one, and representations of a tubular algebra},
Representations of algebras (Ottawa, ON, 1992), 313-337, CMS Conf. Proc., 
\textbf{14}. 
\bibitem[LP]{LP}
Lin Y., Peng L., \emph{2-Extended affine Lie algebras and tubular algebras},
to appear in Adv. Math..
\bibitem[Lu]{Lu}
Lusztig G., \emph{Introduction to quantum groups}, Progress in Mathematics, 
\textbf{110}. Birkhauser, (1993).
\bibitem[MRY]{MRY}
Moody R., Rao S., Yokonuma Y., \emph{Toroidal Lie algebras and vertex 
representations}, Geom. Dedicat. \textbf{35}, 283-307 (1990).
\bibitem[N]{Nak}
Nakajima, H., \emph{Quiver varieties and Kac-Moody algebras.},
Duke Math. J. \textbf{91}, no. 3, 515-560, (1998).
\bibitem[M]{M}
Macdonald I., \emph{Symmetric functions and Hall polynomials}, Oxford press 
(1979).
\bibitem[PX]{PX}
Peng, L., Xiao, J., \emph{Triangulated categories and Kac-Moody algebras.},
Invent. Math. \textbf{140}, no. 3, 563-603, (2000).
\bibitem[Re]{Re}
Reineke M., \emph{The Harder-Narasimhan system in quantum groups and 
cohomology of quiver moduli}, math.QA/0204059.
\bibitem[Ri1]{Ri1}
Ringel C., \emph{Hall algebras and quantum groups}, 
Invent. Math. \textbf{101} no. 3, 583-591, (1990).
\bibitem[Ri2]{Ri2}
Ringel C., \emph{Hall algebras}, Topics in algebra, Part 1 (Warsaw 1988), PWN,
Warsaw, 433-447, (1990).
\bibitem[Ri3]{Ri3}
Ringel C., \emph{Tame algebras and integral quadratic forms}, LNM \textbf{1099},
(1984).
\bibitem[S]{S1}
Schiffmann O., \emph{The Hall algebra of a cyclic quiver and canonical bases of Fock spaces.}, 
Internat. Math. Res. Notices, no. \textbf{8}, 413-440, (2000).
\end{thebibliography}
\end{document}